\documentclass[aos,MSNbibl,seceqn,dvips]{arximspdf}
\usepackage{graphicx}

\doi{10.1214/12-AOS1019} %kopijuoti is PTS
\volume{40}
\issue{4}
\pubyear{2012}
\firstpage{2128}
\lastpage{2161}

\makeatletter
\def\mid{\vert}

\newcommand{\llvert}{\vert}
\newtheorem{theorem}{Theorem}[section]
\newproclaim{defin}[theorem]{Definition}
\newproclaim{example}{Example}
\newtheorem{corollary}[theorem]{Corollary}
\newtheorem{proposition}[theorem]{Proposition}

\newproclaim{excont}{Example}
\makeatother

\begin{document}
\begin{frontmatter}

\title{General theory for interactions in sufficient cause models with
dichotomous exposures}
\runtitle{Sufficent cause interactions}

\begin{aug}
\author[a]{\fnms{Tyler J.} \snm{VanderWeele}\corref{}\ead[label=e1]{tvanderw@hsph.harvard.edu}\ead[label=u1,url]{http://www.hsph.harvard.edu/faculty/tyler-vanderweele/}\thanksref{t1}}
\and
\author[b]{\fnms{Thomas S.} \snm{Richardson}\ead[label=e2]{thomasr@uw.edu}\thanksref{t2}}
\thankstext{t1}{Supported by the National Institutes of Health (R01 ES017876).}
\thankstext{t2}{Supported by
NSF (CRI 0855230), the National
Institutes of Health (R01 AI032475) and the Institute of Advanced Studies,
University of Bologna.}
\runauthor{T. J. VanderWeele and T. S. Richardson}
\affiliation{Harvard University and University of Washington}
\address[a]{Department of Epidemiology\\
Harvard School of Public Health\\
677 Huntington Avenue\\
Boston, Massachusetts 02115\\
USA\\
\printead{e1}\\
\printead{u1}}

\address[b]{Department of Statistics\\
University of Washington\\
Box 354322\\
Seattle, Washington 98195\\
USA\\
\printead{e2}}
\end{aug}

% HISTORY:
\received{\smonth{4} \syear{2010}}
\revised{\smonth{5} \syear{2012}}

% ABSTRACT
%
\begin{abstract}
The
sufficient-component
cause framework assumes the existence
of sets of
sufficient causes that
bring about an event.
For a binary outcome and an
arbitrary number of
binary causes any
set of potential outcomes can be
replicated by positing a
set of sufficient causes; typically this
representation is not unique. A
sufficient cause
interaction is said to be
present if within all
representations
there exists a sufficient cause in
which two or more
particular causes are all
present. A singular interaction is said
to
be present if for some subset of
individuals there is a unique minimal
sufficient cause. Empirical and
counterfactual conditions
are given for sufficient cause interactions and singular
interactions between an arbitrary
number of causes.
Conditions
are given for cases in which
none, some or
all of a given set of causes
affect the outcome monotonically.
The relations
between
these results, interactions in
linear statistical models and
Pearl's probability of
causation are
discussed.
\end{abstract}

% KEYWORDS
%

\end{frontmatter}

%s1 #&#
%s1 ###
\section{Introduction}

Rothman's sufficient-component cause model~\cite{Rothman1976}
postulates a
set of different causal mechanisms, each sufficient to bring about the
outcome under consideration. Rothman refers to these hypothesized causal
mechanisms as ``sufficient causes,''
conceiving of them as minimal sets of actions, events or states of nature
which together initiate a process resulting in the outcome.
%likely be many different sufficient causes i.e. many different causal
%mechanisms by which the outcome could come about.

Thus each sufficient cause is hypothesized to consist of a set of
``component causes.''  Whenever all components of a particular sufficient
cause are present, the outcome occurs; within every sufficient cause, each
component would be necessary for that sufficient cause to lead to the
outcome.  Models of this kind have a long history: a simple version is
considered by Cayley~\cite{cayley1853}; it also corresponds to the INUS
model introduced by Mackie~\cite{Mackie1965} in the philosophical
literature; see also~\cite{blisstoxicity1939} for an early application.
Much recent work has sought to relate the model to other causal
modeling frameworks \cite%
{Greenland2002,Flanders2006,VanderWeele2006,VanderWeele2007a,VanderWeele2008a}.

In traditional sufficient-component cause [SCC] models, the outcome and
all the component
causes are events, or equivalently, binary random variables. An SCC model
with $k$ component causes implies a set of $2^k$ potential outcomes.
Conversely, in Section~\ref{secsccframe} we show that for any given
list of
potential outcomes there is at least one SCC model which represents this
set. However, in general there may be many such SCC models.

One question concerns whether, given a set of potential outcomes
implied by
some (unknown) SCC model, one may infer that two component causes are
present within some sufficient cause in the unknown SCC model. In general,
it is possible that two SCC models both imply the same set of potential
outcomes, yet although~$A$ and $B$ occur together in some sufficient
component cause in the first model, $A$ and $B$ are not present
together in
any sufficient component cause in the second. In \cite
{VanderWeele2008a} two
sufficient component causes are said to form a ``sufficient cause
interaction'' (or to be ``irreducible'') if they are both present
within at
least one sufficient cause in \textit{every} SCC model for a given set of
potential outcomes. Of course, in general, the distribution of potential
outcomes for a given population is also unknown, though it is constrained
(marginally) by the observed data from a randomized experiment. In \cite
{VanderWeele2008a} empirical conditions are given which are sufficient to
ensure that for any set of potential outcomes compatible with experimental
data, all compatible SCC models will contain a sufficient cause
involving $A$
and~$B$. These results were an improvement upon earlier empirical tests for
the existence of a two-way interaction in an SCC model~\cite{Rothman1998},
which required the assumption of monotonicity; see also~\cite{Koopman1981,Greenland1988,Novick2004,Aickin2002,VanderWeele2007c}.
The new results are able
to establish the existence of an interaction in situations where monotonicity
does not hold.
%%%%%%%
%%%%%%%
% \ If so, then there
%exists some causal mechanism in which both component causes are needed
%for
%the mechanism to operate and it is then appropriate to speak of
%synergism. \
%Rothman and Greenland give an empirical condition for
%synergism between two binary variables in the sufficient-component
%cause
%framework under the assumption that both causes affect the outcome
%monotonically. \ VanderWeele and Robins~\cite{VanderWeele2008a}
%provide a
%formal definition of a sufficient cause interaction, corresponding to
%the
%concept of synergism, and derive an empirical condition for a 2-way
%sufficient cause interaction without making use of monotonicity
%assumptions.
In this paper we develop empirical conditions that are sufficient for the
existence of a sufficient cause containing a given subset of an arbitrary
number of variables, both with and without monotonicity assumptions.

As illustrative motivation for the theoretical development, we will consider
data presented in a study by~\cite{Taylor1998}, summarized in Table~\ref
{tabdata},
from a case-control study of bladder cancer examining possible three-way
interaction between smoking ($1={}$present), and genetic variants on NAT2 ($0=R$,
$1=S$ genotype) and NAT1 ($1$~for the *10 allele) for Caucasian individuals.
We return to this example at the end of this paper to examine the
evidence for a
sufficient cause containing all three: smoking, the S genotype on NAT2 and
the *10 allele on NAT1.

%
%t1 #&#
%t1 ###
\begin{table}
\caption{Case-control data from a study of bladder cancer \protect\cite{Taylor1998}}\label{tabdata}
\begin{tabular*}{\textwidth}{@{\extracolsep{\fill}}lccccc@{}}
\hline
&  &  & \multicolumn{1}{c}{\textbf{Cases}} & \multicolumn{1}{c}{\textbf{Controls}} &
\\
\textbf{Smoking}& \textbf{NAT2}&\textbf{NAT1*10} &\multicolumn{1}{c}{$\bolds{(n=215)}$} &\multicolumn{1}{c}{$\bolds{(n=191)}$}
& \multicolumn{1}{c@{}}{\textbf{Odds ratio (95\% CI)}}\\
\hline
0 & 0 & 0  & \phantom{0}6  &13 & 1\phantom{0000000000.}\\
0 & 0 & 1  & \phantom{0}8  & 16 & 1.1 (0.3, 3.9)\phantom{0} \\
0 & 1 & 0  &16  &31 & 1.1 (0.4, 3.5)\phantom{0} \\
0 & 1 & 1  & \phantom{0}6  &10 & 1.3  (0.3, 5.3)\phantom{0} \\
1 & 0 & 0  & 42  &32 & 2.8  (1.1, 8.3)\phantom{0} \\
1 & 0 & 1  & 41  &26 & 3.4  (1.2, 10.1) \\
1 & 1 & 0 & 61  &51 & 2.6  (0.9, 7.3)\phantom{0} \\
1 & 1 & 1  & 35  &12 & 6.3 (2.0, 20.3)\\
\hline
\end{tabular*}
\end{table}

The remainder of this paper is organized as follows:  Section~\ref{secsccframe} presents the sufficient-component cause framework as
formalized by VanderWeele and Robins~\cite{VanderWeele2008a}.
%the section also gives a theorem that
%demonstrates that for a binary outcome and an arbitrary number of
%binary
%causes, given any given set of potential outcomes it is possible to
%construct a set of sufficient causes which replicates the potential
%outcome
%responses. \
Section~\ref{secirreducible} describes general $n$-way irreducible
interactions (aka ``sufficient cause interactions'') and characterizes these
in terms of potential outcomes. Section~\ref{sectestsnway} derives
empirical conditions for the existence of irreducible interactions both with
and without monotonicity assumptions. Section~\ref{secsingular} describes
``singular'' interactions which arise in genetic contexts, provides a
characterization, derives empirical conditions that are sufficient for their
existence and relates this notion to Pearl's probability of causation.
Section~\ref{seclinearlink} discusses the relation between singular
and sufficient cause
interactions and linear statistical models. Section~\ref{secinterpretation}
provides some comments regarding stronger interpretations of sufficient
cause models, and returns to the data presented in Table~\ref{tabdata}.
Finally Section~\ref{secconclude} offers some possible
extensions to the present work.

%s2 #&#
%s2 ###
\section{Notation and basic concepts}
\label{secsccframe}

%In this section we present the sufficient-component cause framework as
%formalized by VanderWeele and Robins~\cite{VanderWeele2008a}. \ We
%also give
%a result relating the sufficient-component cause framework to the
%potential
%outcomes framework in the case of $n$ binary causes; the result
%generalizes
%that presented in VanderWeele and Robins~\cite{VanderWeele2008a} for
%the
%case of $n=2$. We conclude the section by providing definitions
%concerning
%sufficient cause interactions and by providing a formalization of the
%idea
%of testing for the joint presence of two or more causes in the same
%causal
%mechanism.

We will use the following notation: An \textit{event} is a binary random
variable taking values in $\{0,1\}$. We use uppercase roman to indicate
events ($X$), boldface to indicate sets of events ($\mathbf{C}$), and
lowercase to indicate specific values both for single random variables
($X=x$%
), and, with slight abuse of notation, for sets $\{\mathbf{C=c}\} \equiv
\{\forall i, (\mathbf{C})_i=(\mathbf{c})_i \}$ and $\{\mathbf{a\leq b}\}
\equiv\{\forall i, (\mathbf{a})_i\leq(\mathbf{b})_i \}$; $\mathbf{1}$~and $%
\mathbf{0}$ are vectors of $1$'s and $0$'s; the cardinality of a set is
denoted $|\mathbf{C}|$. We use fraktur (${\mathfrak{B}}$) to denote
collections of sets of events.

The complement of some event $X$ is denoted by $\overline{X} \equiv
1-X$. A
\textit{literal} event associated with $X$, is either $X$ or $\overline{X}$.
For a given set of events $\mathbf{C}$, $\mathbb{L}(\mathbf{C})$ is the
associated set of literal events
\[
\mathbb{L}(\mathbf{C}) \equiv\mathbf{C} \cup\{ \overline{X} \mid X \in
\mathbf{C}\}.
\]
For a literal $L\in\mathbb{L}(\mathbf{C})$, and an assignment $\mathbf{c}$
to $\mathbf{C}$, $(L)_{\mathbf{c}}$ denotes the value assigned to $L$
by $%
\mathbf{c}$. The \textit{conjunction} of a set of literal events
$\mathbf{B} =
\{F_1,\ldots,F_m\} \subseteq\mathbb{L}(\mathbf{C})$ is defined as
\[
\bigwedge (\mathbf{B}) \equiv\prod
_{i=1}^m F_i = \min\{F_1,
\ldots,F_m\};
\]
note that $\bigwedge (\mathbf{B}) =1$ if and only if for
\textit{all} $%
i $, $F_i=1$. We also define $B_1 \wedge B_2 \equiv\bigwedge \{B_1,B_2\}$.
We will use ${\mathbb{I}}(A)$ to denote the indicator function
for event $A$. There is a simple correspondence between conjunctions of
literals and indicator functions: let $\mathbf{B}=\{X_1,\ldots,X_s\}$
and $%
\mathbf{C}=\{Y_1,\ldots,Y_t\}$, then
%
%e1 #&#
%e2.1 ###
\begin{equation}
\label{eqindicator} \bigwedge \bigl(\{X_1,
\ldots,X_s,\overline{Y}_{ 1},\ldots,%
\overline{Y}_{ t}\} \bigr) = 1\quad \Leftrightarrow\quad {\mathbb{I}} \bigl(
\{%
\mathbf{B} = \mathbf{1}, \mathbf{C} = \mathbf{0}\} \bigr)=1.
\end{equation}
Similarly, the \textit{set of literals corresponding to an assignment $%
\mathbf{c}$ to $\mathbf{C}$} is defined
\[
\mathbf{B}^{[\mathbf{c}]} \equiv\bigl\{L \mid L\in{\mathbb{L}}(\mathbf{C}),
(L)_{\mathbf{c}}=1\bigr\}
\]
so that $\bigwedge (\mathbf{B}^{[\mathbf{c}]}) = {\mathbb{I%
}}(\mathbf{C}=\mathbf{c})$; note that $|\mathbf{B}^{[\mathbf
{c}]}|=|\mathbf{C%
}|$. The \textit{disjunction} of a set of binary random variables is defined
as
\[
\bigvee \bigl(\{Z_1,\ldots, Z_p\}
\bigr) \equiv\max\{Z_1,\ldots, Z_p\};
\]
note that $\bigvee (\{Z_1,\ldots, Z_p\}) =1$ if and only if for
\textit{some} $j$, $Z_j=1$. Similarly we let $B_1 \vee B_2 \equiv\bigvee\{B_1,B_2\}$.
Given a collection of sets of literals ${\mathfrak{%
B}} = \{\mathbf{C}_1,\ldots, \mathbf{C}_q\}$, we define
\[
\bigvee \bigwedge({\mathfrak{B}}) \equiv%
\bigvee _i \Bigl(\bigwedge (
\mathbf{C}%
_i) \Bigr).
\]
We use $\dot{\mathbb{P}}(\mathbb{L}(\mathbf{C}))$ to denote the set of
subsets of $\mathbb{L}(\mathbf{C})$ that do not contain both $X$ and $%
\overline{X}$ for any $X\in\mathbf{C}$; more formally,
\[
\dot{\mathbb{P}}\bigl(\mathbb{L}(\mathbf{C})\bigr) \equiv \bigl\{ \mathbf{B }
\llvert \mathbf{B} \subset\mathbb{L}(\mathbf{C}) \mbox{ for all } X \in
\mathbf{C}, \{X,\overline{X}\} \nsubseteq\mathbf{B}  \bigr\} .
\]
Note that if $\mathbf{B} \in\dot{\mathbb{P}}(\mathbb{L}(\mathbf{C}))$,
and $%
|\mathbf{B}|=|\mathbf{C}|$, so that for all $C\in\mathbf{C}$, exactly one
of $C$ or $\overline{C}$ is in $\mathbf{B}$, then an assignment of
values $%
\mathbf{b}$ to $\mathbf{B}$ induces a unique assignment $\mathbf{c}$ to
$%
\mathbf{C}$ and vice versa. %
%will be written as $X_{1}\ldotsX_{m}$. \ The disjunctive or OR
%operator, $\tbigvee$, is defined by $A\tbigvee B=A+B-AB$ so that $A
%B=1$ if and only if either $A=1$ or $B=1$. \ \

%s2.1 #&#
%s2.1 ###
\subsection{Potential outcomes models}
\label{secpotoutcomes}
Consider a potential outcome model~\cite{Neyman1923,Rubin1974,Rubin1978}
with $s$ binary factors, $X_{1},\ldots,X_{s}$,
which represent hypothetical interventions or causes, and let $D$
denote some
binary outcome of interest.  We use $\Omega$ to denote the sample
space of
individuals in the population and use $\omega$ for a particular sample
point. Let $D_{x_{1},\ldots,x_{s}}(\omega)$ denote the counterfactual value
of $%
D $ for individual $\omega$ if the cause $X_{j}$ were set to the
value $%
x_{j} $ for $j=1,\ldots,s$. The potential outcomes framework we employ makes
two assumptions: first, that for a given individual these counterfactual
variables are deterministic; second, in asserting that the
counterfactual $%
D_{x_{1},\ldots,x_{s}}(\omega)$ is well defined, it is implicitly assumed that
the value that $D$ would take on for individual\vadjust{\goodbreak} $\omega$ is determined
solely by the values that $X_{1},\ldots,X_{s}$ are assigned for this
individual, and not the assignments made to these variables for other
individuals $\omega^\prime$. This latter assumption is often called ``no
interference''~\cite{Cox1958}, or the stable unit treatment value
assumption (SUTVA)
\cite{Rubin1990}. An
example of a situation where this assumption might fail is a vaccine trial
where there is ``herd'' immunity.

We will use $D_{x_{1},\ldots, x_{s}}(\omega)$, $D_{X_{1}=x_{1},\ldots
,X_{s}=x_{s}}(\omega)$, $D_{\mathbf{c}}$ and $D_{\mathbf{C=c}}(\omega)$,
with $\mathbf{C}=\{X_{1},\ldots,X_{s}\}$ interchangeably.  In this setting
there will be $2^{s}$ potential outcomes for each individual $\omega$ in
the population, one potential outcome for each possible value of $%
(X_{1},\ldots,X_{s})$; we use ${\mathcal{D}}(\mathbf{C};\omega)$ to
denote the
set of all such potential outcomes for an individual, and ${\mathcal
{D}}(%
\mathbf{C};\Omega)$ for the population. Note that if $G=g(\mathbf{C})$ is
some deterministic function of $\mathbf{C}$, then $G_{\mathbf
{C=c}}(\omega
)=g(\mathbf{c)}$, and hence is constant; thus our usage is consistent with
the definition of $(L)_{\mathbf{c}}$ in the previous section.

%
%t2 #&#
%t2 ###
\begin{table}
\tabcolsep=0pt
\caption{All potential outcomes and actual outcomes for three binary
causes, $%
X_{1}$, $X_{2} $ and $X_{3}$, in a~population with two individuals}\label{taboutcomes}
\begin{tabular*}{\textwidth}{@{\extracolsep{\fill}}lcccccccccc@{}}
\hline
\textbf{Individual} & $\bolds{D_{000}}$ & $\bolds{D_{001}}$ &
$\bolds{D_{010}}$ & $\bolds{D_{011}}$ & $\bolds{D_{100}}$ & $ \bolds{D_{101}}$ & $ \bolds{D_{110}}$ & $ \bolds{D_{111}}$
& $\bolds{(X_1,X_2,X_3)}$ & $\bolds{ D}$ \\
\hline
1 & 0 & 1 & 1 & 0 & 0 & 1 & 1 & 0 & $(1,0,1)$ & 1 \\
2 & 0 & 1 & 1 & 0 & 0 & 1 & 1 & 1 & $(0,0,0)$ & 0 \\
\hline
\end{tabular*}   \vspace*{-3pt}
\end{table}

The actual observed value of $D$ for individual $\omega$ will be
denoted by
$D(\omega)$ and similarly the actual value of $X_{1},\ldots,X_{s}$ by
$%
X_{1}(\omega),\ldots,X_{s}(\omega)$. Actual and counterfactual
outcomes are
linked by the \textit{consistency axiom} which requires
% Mathematically, it could be that $%
%D_{X_{1}(\omega),\ldots,X_{s}(\omega)}(\omega)\neq$ $D(\omega)$;
%however we
%will require the ``consistency" assumption
that
%
%e2 #&#
%e2.2 ###
\begin{equation}
\label{eqcons} D_{X_1=X_{1}(\omega),\ldots,X_s=X_{s}(\omega)}(\omega)= D(\omega),
\end{equation}
that is, that the value of $D$ which would have been observed if $%
X_{1},\ldots,X_{s} $ had been set to the values they actually took is
equal to
the value of $D$ which was in fact observed~\cite{Robins1986}. It
follows from this axiom that
$D_{X_{1}(\omega),\ldots,X_{s}(\omega)}(\omega)=D$ is observed, but
it is the
only potential outcome for individual $\omega$ that is observed. % is
%the
%potential outcome , the value
%of $D$ which would have been observed if $X_{1},,\ldots,,X_{s}$ had been
%set to
%what they in fact were.
% Since $X_{1},,\ldots,,X_{s}$ are binary, all of the
%potential outcomes for an individual $\omega$ can be listed in a
%vector
%with $2^{s}$ components and this vector we will denote by $\mathcal{D}%
%(\omega)$.

%
%ex1 #&#
\begin{example}
Consider a binary outcome $D$ with three binary causes of interest, $%
X_{1}$, $X_{2}$ and $X_{3}$.  Suppose that the population consists of two
individuals. The potential outcomes (left-hand side) and actual
outcomes (right-hand side) are
shown in Table~\ref{taboutcomes}. \label{ex2omegas}
\end{example}

We use the notation $A \perp\!\!\!\perp B \mid C$ to
indicate that $A$ is independent of $B$, conditional on $C$ in the population
distribution.

%s2.2 #&#
%s2.2 ###
\subsection{Definitions for sufficient cause models}
\label{secscms}

The following definitions generalize those in~\cite{VanderWeele2008a} to
sub-populations, $\varnothing\neq\Omega^* \subseteq\Omega$:\vadjust{\goodbreak} %
%%% TSR COMMENT: I don't believe the next statement:
%%%
%%% This is only true if we interpret the SCC model in the strong sense.
%%%%
% As
%discussed at the end of this section, the definitions of a sufficient
%cause
%and minimal sufficient cause correspond with the notion of a causal
%mechanism.

%
%de2.1 #&#
\begin{defin}[(Sufficient cause)]
\label{defsc} A subset $\mathbf{B}$ of the putative (binary) causes
${\mathbb{L}}(\mathbf{C})$ for $D$ forms a \textit{sufficient cause
for $D$ (relative to $\mathbf{C}$) in sub-population } $\Omega^*$ if
for all
$\mathbf{c}\in\{0,1\}^{|\mathbf{C}|}$ such that $(\bigwedge (\mathbf{B}))_{%
\mathbf{c}} = 1$, $D_{\mathbf{c}}(\omega)=1$ for all $\omega\in\Omega^*
\subseteq\Omega$. [We require that there exists a $\mathbf{c}%
^*$ such that $(\bigwedge (\mathbf{B}))_{\mathbf{c^*}} = 1$.]
\end{defin}

Observe that if $\mathbf{B}$ is a sufficient cause for $D$, then any
intervention setting the variables $\mathbf{C}$ to $\mathbf{c}$ with
$(\bigwedge (\mathbf{B}))_{%
\mathbf{c}} = 1$ will ensure
that $D_{\mathbf{c}}(\omega) =1$ for all $\omega\in\Omega^*$. We restrict
the definition to nonempty sets $\Omega^*$, to preclude every set
$\mathbf{B%
}$ being a sufficient cause in an empty sub-population. Likewise we require
that there exists some~$\mathbf{c}^*$ such that $(\bigwedge (\mathbf{B}))_{%
\mathbf{c^*}} = 1$ in order to avoid logically inconsistent conjunctions,
for example, $X_1 \wedge\overline{X}_1$, being classified (vacuously)
as a
sufficient cause. As a direct consequence, for any binary random
variable $X$, at most one of $X$ and $\overline{X}$ appear in any sufficient cause~$\mathbf{B}$.

%
%pr2.2 #&#
\begin{proposition}
\label{propsmallerset} In $\Omega^*$ if $\mathbf{B}$ is a sufficient cause
for $D$ relative to $\mathbf{C}$, then $\mathbf{B}$ is sufficient for
$D$ in
any set $\mathbf{C}^*$ with $\mathbf{B} \subseteq\mathbf{C}^* \subseteq
\mathbf{C}$.
\end{proposition}
$\mathbf{B}$ may be sufficient for $D$ relative to $\mathbf
{C}$ in
$\Omega^*$, but not relative to $\mathbf{C}^\prime\supset\mathbf{C}$.

%
%pr2.3 #&#
\begin{proposition}
\label{propsubset} If $\mathbf{B}$ is a sufficient cause for $D$ relative
to $\mathbf{C}$ in $\Omega^*$, then $\mathbf{B}$ is sufficient for $D$
relative to $\mathbf{C}$ in any subset $\varnothing\neq\Omega^{**}
\subseteq
\Omega^*$.
\end{proposition}
$\mathbf{B}$ may be sufficient for $D$ relative to $\mathbf
{C}$ in
$\Omega^*$, but not in $\Omega^\prime\supset\Omega^*$.

%
%de2.4 #&#
\begin{defin}[(Minimal sufficient cause)]
\label{defmsc} A set $\mathbf{B} \subset{\mathbb{L}}(\mathbf{C})$
forms a
\textit{minimal sufficient cause for $D$ (relative to $\mathbf{C}$) in
sub-population $\Omega^*$} if $\mathbf{B}$ constitutes a sufficient cause
for $D$ in $\Omega^*$, but no proper subset $\mathbf{B}^*\subset\mathbf{B}$
also forms a sufficient cause for $D$ in $\Omega^*$.
\end{defin}

Note that (in some $\Omega^*$) $\mathbf{B}$ may be a minimal sufficient
cause for $D$ relative to $\mathbf{C}$, but not relative to $\mathbf{C}^*
\subset\mathbf{C}$, so the analog of Proposition~\ref{propsmallerset} does
not hold. For individual $2$ in Table~\ref{taboutcomes} $\{
{X}_1,{X}_3\}$
is a minimal sufficient cause relative to $\{X_1,X_2,X_3\}$. However,
if we
suppose that for $\omega= 2$, $X_2$ is not caused by $X_1$ and $X_3$, so
for all $x_1$, $x_3$, $X_{2 X_1=x_1,X_3=x_3}(\omega= 2) =
X_2(\omega= 2)$, then $\{{X}_1,{X}_3\}$ is not a minimal sufficient cause
relative to $\{X_1, X_3\}$.
\begin{eqnarray*}
D_{X_1=0,X_3=1}(\omega= 2) &=& D_{X_1=0,X_2=0,X_3=1}(\omega= 2)=1,
\\
D_{X_1=1,X_3=1}(\omega= 2) &=& D_{X_1=1,X_2=0,X_3=1}(\omega= 2)= 1
\end{eqnarray*}
[since $X_2(\omega= 2)=0$]; hence $X_3$ is a sufficient cause of $D$
relative to $\{X_1,X_3\}$; hence $\{X_1,{X}_3\}$ is not minimal
relative to $%
\{X_1, X_3\}$ for $\omega= 2$.

Similarly, if $\mathbf{B}$ is a minimal sufficient cause for $D$
relative to
$\mathbf{C}$ in $\Omega^*$, it does not follow that $\mathbf{B}$ is a
\textit{minimal} sufficient cause for $D$ relative to $\mathbf{C}$ in
subsets $\Omega^{**}\subseteq\Omega^*$, so the analog to Proposition~\ref{propsubset}
does not hold. In particular, it may be the case that for
all $%
\omega\in\Omega^*$, $\mathbf{B}$ is not a minimal sufficient cause
for $D$
in $\{\omega\}$.

In the language of digital circuit theory~\cite{marcovitz2001}, sufficient
causes are termed ``implicants,'' and minimal sufficient causes are ``prime
implicants.''

%When a set of several sufficient causes for $D$ is such that $D$
%occurs if
%and only if one of the sufficient causes is realized then the set of
%sufficient causes is said to be determinative for $D$.

%
%de2.5 #&#
\begin{defin}[(Determinative set of sufficient causes)]
\label{defdsc} A set of sufficient causes for $D$, ${\mathfrak{B}} = \{
\mathbf{B}_{1},\ldots,\mathbf{B}_{n}\} \subseteq\dot{\mathbb
{P}}({\mathbb{L}%
}(\mathbf{C}))$, is said to be \textit{determinative for $D$ (relative
to $%
\mathbf{C}$) in sub-population $\Omega^*$} if for all $\omega\in\Omega^*$
and for all $\mathbf{c}$, $D_{\mathbf{c}}(\omega)=1$ if and only if $(\bigvee \bigwedge ({\mathfrak{B}}))_{\mathbf{c}}=1$.
\end{defin}

We will refer to a determinative set of sufficient causes for $D$ as a
\textit{sufficient cause model}.
Observe that in any sub-population $\Omega^*$ for which there exists a
determinative set of sufficient causes, the vectors of potential outcomes
for $D$ are identical, so ${\mathcal{D}}(\mathbf{C},\omega) = {\mathcal
{D}}(%
\mathbf{C},\omega^{\prime})$ for all $\omega, \omega^\prime\in\Omega^*$.

%
%de2.6 #&#
\begin{defin}[(Nonredundant set of sufficient causes)]
\label{defnrsc} A determinative set of sufficient causes ${\mathfrak{B}}$,
for $D$, is said to be \textit{nonredundant} (in $\Omega^*$,
relative to \textbf{C}) if there is no proper subset ${\mathfrak{B}}%
^* \subset{\mathfrak{B}}$ that is also determinative for $D$.
\end{defin}

Note that sufficient causes are conjunctions, while sets of sufficient
causes form disjunctions of conjunctions; minimality refers to the
components in a particular conjunction, that each component is required for
the conjunction to be sufficient for $D$; nonredundancy implies that each
conjunction is required for the disjunction of the set of conjunctions
to be
determinative. If for some set of sufficient causes $\mathfrak{B}
\subseteq
\dot{\mathbb{P}}(\mathbb{L}(\mathbf{C}))$, for all $X \in\mathbf{C}$, and
all $\mathbf{B} \in\mathfrak{B}$, either $X \in\mathbf{B}$ or
$\overline{X}
\in\mathbf{B}$, then $\mathfrak{B}$ is a nonredundant set of sufficient
causes.

%incremented
%environment only used for such continuations

%
%ex1 #&#
\begin{excont}[(Revisited)]%\hyperref[ex:2omegas]
The set $\mathfrak{B}_1 = \{ \{X_{1},X_{2}\}, \{X_{2},\overline{X}%
_3\}, \{\overline{X}_2,X_{3}\}\}$ forms a determinative set of sufficient
causes for the individual $\omega=2$, since
%
%e3 #&#
%e2.3 ###
\begin{equation}
\label{eqomega2a} D_{\mathbf{c}}(\omega=2)= \bigl(\strut(X_{1}
\wedge X_{2}) \vee (X_{2} \wedge
\overline{X}_3)\vee (\overline{X}%
_2
\wedge X_{3}) \bigr)_{\mathbf{c}}
\end{equation}
as does $\mathfrak{B}_2 = \{ \{X_{1},X_{3}\}, \{X_{2},\overline{X}_3\},
\{\overline{X}_2,X_{3}\}\}$:
%
%e4 #&#
%e2.4 ###
\begin{equation}
\label{eqomega2b} D_{\mathbf{c}}(\omega=2)= \bigl(\strut(X_{1}
\wedge X_{3}) \vee (X_{2} \wedge
\overline{X}_{3}) \vee (\overline{X}%
_{2}
\wedge X_{3}) \bigr)_{\mathbf{c}}.
\end{equation}
\end{excont}

%environment only used for such continuations
As this example shows, determinative sets of sufficient causes are not,
in general, unique.

%s2.3 #&#
%s2.3 ###
\subsection{Sufficient cause representations for a population}

As noted, if $\mathbf{B}$ is a sufficient cause for $D$ in $\Omega^*$, then
all the units in $\Omega^*$ will have $D=1$ for any assignment $\mathbf{c}$
to $\mathbf{C}$, such that $(\bigwedge (\mathbf{B}))_{%
\mathbf{c}} = 1$. In most realistic settings it is unlikely that any
set $%
\mathbf{B}$ will be sufficient to ensure $D=1$ in an entire population.\vadjust{\goodbreak}
Consequently, different sets of sufficient causes will be required within
different sub-populations. A sufficient cause representation is a set of
sub-populations, each with its own determinative sufficient cause
representation.

%A \textit{sufficient cause representation} for an individual's
%potential outcomes ${\mathcal D}(\mathbf{C};\omega)$ is a determinative
%set of sufficient causes $\mathfrak{B} \subseteq\mathbb{P}(\mathbb{L}(
%for all $\mathbf{c}\in\{0,1\}^{|\mathbf{C}|}$, $D_\mathbf{c}(\omega)
%= (\dnf({\mathfrak B}))_\mathbf{c}$.

%
%de2.7 #&#
\begin{defin}
\label{defpopsuffcause} A \textit{sufficient cause representation} $(
\mathbf{A}, {\mathfrak{B}})$ for ${\mathcal{D}}(\mathbf{C};\Omega)$ is an
ordered set $\mathbf{A} = \langle A_1, \ldots, A_p\rangle$ of binary random
variables, with $(A_i)_{\mathbf{c}} = A_i$ for all $i, \mathbf{c}$, and a
set ${\mathfrak{B}}=\langle\mathbf{B}_1,\ldots,\mathbf{B}_p\rangle$,
with $%
\mathbf{B}_i \in\dot{\mathbb{P}}(\mathbb{L}(\mathbf{C}))$, such that for
all $\omega$, $\mathbf{c}$, $D_{\mathbf{c}}(\omega)=1 \Leftrightarrow
\mbox{ for some }j, A_j(\omega)=1\mbox{ and }  (\bigwedge( \mathbf{B}_j )  )_{\mathbf{c}}=1$.
\end{defin}

%(1) Note I have allowed Ai's that are identically zero.
Note that the binary random variables $A_i$ and the sets
$\mathbf{B%
}_i$ are naturally paired via the orderings of $\mathbf{A}$ and
${\mathfrak{B%
}}$; we will refer to a pair $(A_i, \mathbf{B}_i)$ as \textit
{occurring} in
the representation. The requirement that $(A_i)_{\mathbf{c}} = A_i$ for
all $%
i, \mathbf{c}$ implies that $\mathbf{A}\cap\mathbf{C} = \varnothing$, and
further that the $A_i$ are unaffected by interventions on the $X_i$;
this is
in keeping with the interpretation of the $A_i$ as defining pre-existing
sub-populations with particular sets of potential outcomes for~$D$.

%Q I don't know if we really need the next two propositions ?
% I put them in, to try to tie this section to the previous ones.

%
%pr2.8 #&#
\begin{proposition}
If $(\mathbf{A},{\mathfrak{B}})$ is a sufficient cause representation
for ${%
\mathcal{D}}(\mathbf{C};\Omega)$, then $\mathbf{B}_i$ is a sufficient cause
of $D$ in the sub-population in which $A_i(\omega) = 1$.
\end{proposition}

%
%pr2.9 #&#
\begin{proposition}
If $(\mathbf{A}, {\mathfrak{B}})$ is a sufficient cause representation
for ${%
\mathcal{D}}(\mathbf{C};\Omega)$, then for all $\mathbf{A}^* \subseteq
\mathbf{A}$, if
\[
\varnothing\neq\Omega^{\mathbf{A}^*}_{\mathbf{A}\setminus\mathbf{A}^*} \equiv\bigl\{ \omega\mid
\mbox{ for all }A_i\in\mathbf{A}, A_i(\omega)=1 %
\mbox{ iff } A_i\in\mathbf{A}^*\bigr\},
\]
%
%A_i\notin\mathbf{A}^* \},
%%% Note I opted not to define this set in terms of wedge(A; A in
% A \notin\mathbf{A}) because I want to reserve wedge for causes; I
%want to avoid any ambiguity regarding the status of the Ai's solely as
% indicators for subpopulations.
then
\[
{\mathfrak{B}}^{\mathbf{A}^*} \equiv\bigl\{ \mathbf{B}_i \mid
\mathbf{B}_i \in{%
\mathfrak{B}}; A_i \in
\mathbf{A^*}\bigr\}
\]
forms a determinative set of sufficient causes (relative to $%
\mathbf{C}$) for $\Omega^{\mathbf{A}^*}_{\mathbf{A}\setminus
\mathbf{A}^*}$.
\end{proposition}

Note that $\Omega^{\mathbf{A}^*}_{\mathbf{A}\setminus\mathbf
{A}%
^*} $ consists of the sub-population in which $A_i(\omega)=1$ for all $A_i
\in\mathbf{A}^*$ and $A_j(\omega)=0$ for all $A_j \in\mathbf{A}%
\setminus\mathbf{A}^*$.

\begin{pf} Suppose for some $\omega\in\Omega^{\mathbf
{A}^*}_{%
\mathbf{A}\setminus\mathbf{A}^*}$, $\mathbf{B}_j \in{\mathfrak{B}}^{%
\mathbf{A}^*}$, and $\mathbf{c}$ we have\break $(\bigwedge %
( \mathbf{B}_j )  )_{\mathbf{c}}=1$. Since $\omega\in
\Omega^{\mathbf{A}^*}_{\mathbf{A}\setminus\mathbf{A}^*}$, $A_j(\omega) =1$.
It then follows from the definition of a sufficient cause representation
that $D_{\mathbf{c}}(\omega)=1$. Conversely, suppose $D_{\mathbf
{c}}(\omega
)=1$. As $(\mathbf{A}, {\mathfrak{B}})$ is a sufficient cause
representation, for some $j$, $A_j(\omega)=1$ and $ (%
\bigwedge ( \mathbf{B}_j)  )_{\mathbf{c}}=1$. Since,
by hypothesis, $\omega\in\Omega^{\mathbf{A}^*}_{\mathbf{A}\setminus
\mathbf{A}^*}$, it follows that $A_j \in\mathbf{A}^*$, hence $\mathbf{B}_j
\in{\mathfrak{B}}^{\mathbf{A}^*}$.
\end{pf}

%Thus the collection ${\mathfrak B}^{\mathbf{A}^*}$ of all the
%sufficient causes $\mathbf{B}_i$ corresponding to the $A_i$'s in an
%arbitrary subset $\mathbf{A}^*\subset\mathbf{A}$
% is determinative for $D$ in the sub-population $\Omega^{
% $A_j(\omega)=0$ for all $A_j \mathbf{A} \setminus\mathbf{A}^*$.\par

%
%th2.10 #&#
\begin{theorem}
\label{thmrepexists} For any ${\mathcal{D}}(\mathbf{C};\Omega)$, there
exists a sufficient cause representation $(\mathbf{A}, {\mathfrak{B}})$.
\end{theorem}

\begin{pf} Let $p=2^{|\mathbf{C}|}$, and define
${\mathfrak{B}}
\equiv\{\mathbf{B} \mid\mathbf{B}\subseteq\dot{\mathbb{P}}({\mathbb
{L}}(%
\mathbf{C})), |\mathbf{B}|=|\mathbf{C}|\}\equiv\langle\mathbf
{B}_1, \ldots,
\mathbf{B}_p\rangle$, ordered arbitrarily.\vadjust{\goodbreak} Further define $A_i(\omega
)\equiv
D_{\mathbf{B}_i=\mathbf{1}}(\omega)$.
%Further define $\mathbf{b}_i$ to be the unique assignment to $\bf C$
%such that ${\dbigwedge}(\mathbf{B}_i)_
%{\mathbf{b}_i}=1$, and define $A_i(\omega)\equiv D_{\mathbf{b}_i}(
Given an arbitrary $\mathbf{c}$, for some $j$, $\mathbf{B}^{[\mathbf{c}]}=
\mathbf{B}_j$, by construction of ${\mathfrak{B}}$. We then have
\[
D_{\mathbf{c}}(\omega)=1 \quad\Leftrightarrow \quad D_{\mathbf{B}_j=\mathbf{1}%
}(\omega)=1\quad
\Leftrightarrow\quad A_j(\omega)=1\quad\mbox{and}\quad \Bigl(%
\bigwedge ( \mathbf{B}_j) \Bigr)_{\mathbf{c}}=1
\]
as required. The last step follows since by definition $\mathbf
{B}^{[\mathbf{%
c}]}=\mathbf{1}$ if and only if $\mathbf{C}=\mathbf{c}$.
%Let $\mathbf{A}^*\subseteq\langle A_1,\ldots,A_p\rangle$. For all $
% A_i(\omega)=1 \Leftrightarrow D_{\mathbf{B}_i=\mathbf{1}}(\omega)=1.
%Hence for all $\omega\in\Omega^{\mathbf{A}^*}_{\mathbf{A}\setminus
%B}^{\mathbf{A}^*}))_\mathbf{c}=1$ iff
%for some $A_i\in\mathbf{A}^*$, $\dbigwedge(\mathbf{B}_i)_\mathbf{c}=
%For each conjunction $F_{1}^{i},\ldots,F_{n_{i}}^{i}$ we may define $A_{i}(
%)=1$ if $n_{i}=s$ and if $D_{x_{1},\ldots,x_{s}}(\omega)=1$ for
%$x_{1},,\ldots,,x_{s}$
%which satisfy $F_{1}^{i},\ldots,F_{s}^{i}=1$ and $A_{i}(\omega)=0$
%otherwise.
\end{pf}

\cite{VanderWeele2008a} prove this for the case of $|\mathbf{C}|=2$;
see also \cite%
{Flanders2006} and~\cite{VanderWeele2006} for discussion of the case
$|\mathbf{C}|=1$.

\setcounter{excont}{0}
\begin{excont}[(Revisited)]%{{\hyperref[ex:2omegas]
The construction given in the proof of Theorem~\ref{thmrepexists} would yield the following sets of sufficient causes to
represent ${\mathcal{D}}(\mathbf{C};\Omega)$ shown in Table~\ref{taboutcomes}:%
%
%e5 #&#
%e2.5 ###
\begin{eqnarray}\label{eq21}
{\mathfrak{B}} &=& \langle\mathbf{B}_1,\ldots, \mathbf{B}_8
\rangle
\nonumber
\\
&=& \bigl\langle\{ X_1,X_2,X_3\}, \{
X_1,X_2,\overline{X}_3\}, \{
X_1,%
\overline{X}_2,X_3\}, \{
X_1,\overline{X}_2,\overline{X}_3\},
\\
&&\hspace*{6pt} \{ \overline{X}_1,X_2,X_3\}, \{
\overline{X}_1,X_2,\overline{X} %%
_3
\}, \{ \overline{X}_1,\overline{X}_2,X_3\},
\{ \overline {X}_1,\overline{X}%
_2,
\overline{X}_3\} \bigr\rangle
\nonumber
\end{eqnarray}
$\mbox{with }A_1={\mathbb{I}}\bigl(\{\omega=2\}\bigr),
A_4=A_5=A_8=0, A_2=A_3=A_6=A_7=1.$
\end{excont}

\section{Irreducible conjunctions}
\label{secirreducible}

We saw in Example~\ref{ex2omegas} above with $\omega=2$ that an
individual's potential
outcomes may be such that there are two determinative sets of common
causes $%
\mathfrak{B}$ and $\mathfrak{B}^\prime$ and $\{X_1,X_2\}$ is in
$\mathfrak{B}
$, but not in $\mathfrak{B}^\prime$. However, certain conjunctions are such
that in every representation either the conjunction is present or it is
contained in some larger conjunction; such conjunctions are said to be
``irreducible.''

%Suppose that $F_{1},,\ldots,,F_{m}$ are such that each $F_{k}$ is
%either a member of the set of binary causes $\{X_{1},,\ldots,,X_{s}\}$ or
%is the
%complement of such a member then $F_{1},,\ldots,,F_{m}$ is said to exhibit
%nonredundant
%minimal sufficient cause representation for $D$ there exists within the
%representation a sufficient cause which contains $F_{1},,\ldots,,F_{m}$
%within
%its conjunction.

%
%de3.1 #&#
\begin{defin}
\label{defsci} $\mathbf{B} \in\dot{\mathbb{P}}({\mathbb{L}}(\mathbf{C}))$
is said to be %a \textit{sufficient cause interaction (or to be
\textit{irreducible for ${\mathcal{D}}(\mathbf{C},\Omega)$} if in every
representation $(\mathbf{A},{\mathfrak{B}})$ for ${\mathcal{D}}(\mathbf
{C}%
,\Omega)$, there exists $\mathbf{B}_i \in{\mathfrak{B}}$, with $\mathbf
{B}%
\subseteq\mathbf{B}_i$.
\end{defin}

\cite{VanderWeele2008a} also refer to irreducibility of $\mathbf{B}$
for ${%
\mathcal{D}}(\mathbf{C},\Omega)$ as a ``sufficient cause interaction''
between the components of $\mathbf{B}$. (Note, however, that if $\mathbf{B}$
is irreducible, this does not, in general, imply that $\mathbf{B}$ is
either a
minimal sufficient cause, or even a sufficient cause, only that there
is a
sufficient cause that contains $\mathbf{B}$.) It can be shown (via
Theorem~\ref{thmcharsci} below) that $\{X_{2},\overline{X}_{3}\}$ and $\{
\overline{%
X}_{2},X_{3}\}$ are irreducible for ${\mathcal{D}}(\mathbf{C};\Omega)$ in
Table~\ref{taboutcomes}. In Section~\ref{secinterpretation} we
provide an
interpretation of irreducibility in terms of the existence of
a mechanism involving the variables in $\mathbf{B}$. %
Using $\mathbf{C}_{1}\,\dot{\cup}\,\mathbf{C}_{2}$ to indicate the
disjoint union of $\mathbf{C}_{1}$ and $\mathbf{C}_{2}$, we now
characterize irreducibility:

%
%th3.2 #&#
\begin{theorem}
\label{thmcharsci} Let $\mathbf{C}=\mathbf{C}_1 \,\dot{\cup}\, \mathbf
{C}_2$, $%
\mathbf{B} \in\dot{\mathbb{P}}({\mathbb{L}}(\mathbf{C}_1))$, $|\mathbf
{B}%
| = |\mathbf{C}_1|$.\break Then $\mathbf{B}$ is irreducible for ${\mathcal
{D}}(%
\mathbf{C},\Omega)$ if and only if there exists $\omega^*\in\Omega$
and values $%
\mathbf{c}^*_2$ for $\mathbf{C}_2$ such that: \textup{(i)} $D_{\mathbf{B} =
\mathbf{1}, \mathbf{C}_2=\mathbf{c}_2^*}(\omega^*)=1$; \textup{(ii)}
for all $L \in\mathbf{B}$,\break $D_{\mathbf{B}\setminus\{L\} = \mathbf{1},
L =
0, \mathbf{C}_2=\mathbf{c}_2^*}(\omega^*)=0$.
\end{theorem}

Thus $\mathbf{B}$
is irreducible if and only if there exists an individual in $\Omega$ who
would have response $D=1$ if every literal in $\mathbf{B}$ is set to $1$,
but $D=0$ whenever one literal is set to $0$ and the rest to $1$ (in some
context $\mathbf{C}_{2}=\mathbf{c}_{2}^{\ast}$). Note that conditions~(i)
and~(ii) are equivalent to
%
%e6 #&#
%e3.1 ###
\begin{equation}
D_{\mathbf{B}=\mathbf{1},\mathbf{C}_{2}=\mathbf{c}_{2}^{\ast}}\bigl(\omega^{\ast}\bigr)-\sum
_{L\in\mathbf{B}}D_{\mathbf{B}\setminus\{L\}=\mathbf
{1},L=0,%
\mathbf{C}_{2}=\mathbf{c}_{2}^{\ast}}\bigl(\omega^{\ast}\bigr) > 0.
\label{eqcond}
\end{equation}

\begin{pf}
($\Rightarrow$) We adapt the proof of Theorem~\ref{thmrepexists} to show that if for all $\omega\in\Omega$ and assignments
$\mathbf{c}_2^*$ to $\mathbf{C}_2$, at least one of (i) or (ii) does not
hold, then there exists a representation $(\mathbf{A},{\mathfrak{B}})$
for ${%
\mathcal{D}}(\mathbf{C},\Omega)$ such that for all $\mathbf{B}_i \in{%
\mathfrak{B}}$, $\mathbf{B} \nsubseteq\mathbf{B}_i$. Define
\begin{eqnarray*}
{\mathfrak{B}}^{\dag} &\equiv& \bigl\langle\mathbf{B}^\dag_i
\bigr\rangle\equiv \bigl\{ \mathbf{B}^* \mid\mathbf{B}^* \in\dot{\mathbb{P}}
\bigl({\mathbb {L}}(%
\mathbf{C})\bigr), |\mathbf{B}^*|=|\mathbf{C}|,
\mathbf{B} \nsubseteq \mathbf{%
B}^* \bigr\},
\\
{\mathfrak{B}}^{\ddag} &\equiv& \bigl\langle\mathbf{B}^\ddag_i
\bigr\rangle\equiv \bigl\{ \mathbf{B}^* \mid\mathbf{B}^* \in\dot{\mathbb{P}}
\bigl({\mathbb {L}}(%
\mathbf{C})\bigr), |\mathbf{B}^*|=|\mathbf{C}|-1,
\mathbf{B}\setminus\mathbf{B}^* =\{L\}, L\in\mathbf{B} \bigr\},
\end{eqnarray*}
under arbitrary orderings. Thus ${\mathfrak{B}}^{\dag}$ is the set of
subsets of exactly $|\mathbf{C}|$ literals that do not include $\mathbf{B}$
as a subset, while ${\mathfrak{B}}^{\ddag}$ contains those subsets of
size $|%
\mathbf{C}|-1$ that contain all but one literals in $\mathbf{B}$.

\begin{longlist}[]
\item[] For $\mathbf{B}^\dag_i \in{\mathfrak{B}}^{\dag}$ define the
corresponding $A^\dag_i(\omega) \equiv D_{\mathbf{B}^\dag_i=\mathbf{1}%
}(\omega)$;

\item[] For $\mathbf{B}^\ddag_i \in{\mathfrak{B}}^{\ddag}$ define $%
A^\ddag_i(\omega) \equiv D_{\mathbf{B}^\ddag_i=\mathbf{1},L_i=0}(\omega
)D_{%
\mathbf{B}^\ddag_i=\mathbf{1},L_i=1}(\omega)$, where\break $\{L_i\}\equiv
\mathbf{B%
}\setminus\mathbf{B}^\ddag_i$.
\end{longlist}

The representation is given by $(\mathbf{A}, {\mathfrak{B}}) \equiv
(\mathbf{%
A}^\dag \cup \mathbf{A}^{\ddag}, {\mathfrak{B}}^{\dag} \cup {%
\mathfrak{B}}^{\ddag})$, where $\mathbf{A}^\dag\equiv\langle
A_i^\dag\rangle$, $\mathbf{A}^\ddag\equiv\langle A_i^\ddag\rangle$.
To see
this, first note that if for some $\omega$ and $\mathbf{c}$, there is a pair
$(A_j,\mathbf{B}_j)$ in $(\mathbf{A}, {\mathfrak{B}})$ such that $%
A_j(\omega) = 1$ and $(\bigwedge (\mathbf{B}_j))_{%
\mathbf{c}}=1$. Then by construction of $\mathbf{A}^\dag$ and $\mathbf
{A}%
^{\ddag}$ it follows that $D_{\mathbf{c}}(\omega) =1$. For the converse,
suppose that for some $\mathbf{c}$ and $\omega$, $D_{\mathbf{c}}(\omega)=1$.
There are two cases to consider:

\begin{longlist}
\item[$(\bigwedge (\mathbf{B}))_{\mathbf{c}} = 0$.]
In this case $\mathbf{B} \nsubseteq\mathbf{B}^{[\mathbf{c}]}$,
so for some $j$, $\mathbf{B}^\dag_j = \mathbf{B}^{[\mathbf{c}]}$, hence
$%
A_j^\dag(\omega) \equiv D_{\mathbf{B}^\dag_j = \mathbf{1}}(\omega) = D_{
\mathbf{c}}(\omega) =1$, as required.

\item[$(\bigwedge (\mathbf{B}))_{\mathbf{c}} = 1$.] %
Let $\mathbf{c}$ be partitioned as $(\mathbf{c}_1,\mathbf{c}_2)$.
Since (i) holds with $\mathbf{c}_2^*=\mathbf{c}_2$, (ii) does not. Thus for
some $L\in\mathbf{B}$, $D_{\mathbf{B}\setminus\{L\} = \mathbf{1}, L =
0,%
\mathbf{C}_2=\mathbf{c}_2}(\omega)=1$. By construction of ${\mathfrak
{B}}%
^{\ddag}$, for some $j$, $\mathbf{B}^\ddag_j = \mathbf{B}^{[\mathbf{c}%
]}\setminus\{ L\}$, so $(\bigwedge(\mathbf{B}%
^\ddag_j))_{\mathbf{c}} = 1$. Since $1=D_{\mathbf{c}}(\omega)=D_{\mathbf
{B}%
\setminus\{L\} = \mathbf{1}, L = 1, \mathbf{C}_2=\mathbf{c}_2}(\omega
)$, we
have $A_j^\ddag(\omega)=1$, as required.
\end{longlist}

($\Leftarrow$) Suppose for a contradiction, that for some $\omega^*$
and $%
\mathbf{c}_2^*$, (i) and (ii) hold, but $\mathbf{B}$ is not irreducible.
Then there exists a representation $(\mathbf{A},{\mathfrak{B}})$ such that
for all $\mathbf{B}_i \in{\mathfrak{B}}$, $\mathbf{B}\nsubseteq
\mathbf{B%
}_i$. By (i), $D_{\mathbf{B}=\mathbf{1}, \mathbf{c}_2^*}(\omega^*)=1$. Thus
for some pair $(A_j,\mathbf{B}_j)$, $A_j(\omega^*)=1$ and $\mathbf{B}_j
\subseteq\mathbf{B} \cup\mathbf{B}^{[\mathbf{c}_2^*]}$. Since $\mathbf{B}
\nsubseteq\mathbf{B}_j$ there exists some $L\in\mathbf{B}\setminus
\mathbf{B}_j$, but then since $A_j(\omega^*)=1$ and $(\bigwedge(\mathbf{B}_j))_{\mathbf{B}\setminus\{L\}=\mathbf{1}, L=0,
\mathbf{C}_2=\mathbf{c}_2^*} = 1$, we have $D_{\mathbf{B}\setminus\{L\}
=%
\mathbf{1}, L=0, \mathbf{C}_2=\mathbf{c}_2^*}(\omega^*)=1$, which is a
contradiction.
\end{pf}

%
%co3.3 #&#
\begin{corollary}
\label{corexpandpop} If $\mathbf{B}$ is irreducible for ${\mathcal{D}}(
\mathbf{C},\Omega)$, then for any $\Omega^* \supset\Omega$, $\mathbf
{B}$ is
irreducible for ${\mathcal{D}}(\mathbf{C},\Omega^*)$.
\end{corollary}

\begin{pf}
By Theorem~\ref{thmcharsci}, since if $\Omega$
satisfies (i) and (ii), then so does~$\Omega^*$.
\end{pf}

%s3.1 #&#
%s3.1 ###
\subsection{\texorpdfstring{$\mathbf{B}$ irreducible for ${\mathcal{D}}(\mathbf{C},\Omega)$
with $|\mathbf{B}|=|\mathbf{C}|$}
{B irreducible for D(C,Omega) with |B|=|C|}}\label{secbequalsc}

In the special case where $|\mathbf{B}|=|\mathbf{C}|$, the concepts of
minimal sufficient cause for some $\omega^*$ and irreducibility coincide.

%
%pr3.4 #&#
\begin{proposition}
\label{propmscandirred} If $\mathbf{B} \in\dot{\mathbb{P}}({\mathbb
{L}}%
(\mathbf{C}))$ and $|\mathbf{B}| = |\mathbf{C}|$, then $\mathbf{B}$ is a
minimal sufficient cause for some $\omega^* \in\Omega$
relative to $\mathbf{C}$
if and only if $\mathbf{B}$ is
irreducible for ${\mathcal{D}}(\mathbf{C},\Omega)$.
\end{proposition}

\begin{pf} If $|\mathbf{B}|=|\mathbf{C}|$, then condition (i)
in Theorem~\ref{thmcharsci} (taking $\mathbf{C}_2=\varnothing$) holds
if and only if $%
\mathbf{B}$ is a sufficient cause for $D$ for $\omega^*$, and similarly
condition (ii) holds if and only if $\mathbf{B}$ is a minimal
sufficient cause for $D$
(for $\omega^*$).
\end{pf}

Thus we have the following:

%
%co3.5 #&#
\begin{corollary}
If $\mathbf{B} \in\dot{\mathbb{P}}({\mathbb{L}}(\mathbf{C}))$,
$|\mathbf{B%
}|=|\mathbf{C}|$ and $\mathbf{B}$ is a minimal sufficient cause for $D$ for
some $\omega^*\in\Omega$, then $\mathbf{B}\in{\mathfrak{B}}$ for every
representation $(\mathbf{A},{\mathfrak{B}})$ for ${\mathcal{D}}(\mathbf
{C}%
,\Omega)$.
\end{corollary}

\begin{pf} Immediate from Proposition~\ref
{propmscandirred}.
\end{pf}

%s3.2 #&#
%s3.2 ###
\subsection{\texorpdfstring{$\mathbf{B}$ irreducible for ${\mathcal{D}}(\mathbf{C},\Omega)$
with $|\mathbf{B}|<|\mathbf{C}|$}
{B irreducible for D(C,Omega) with |B|<|C|}}

\label{secblessthanc} When $|\mathbf{B}|<|\mathbf{C}|$, the conditions for
irreducibility and for being a minimal sufficient cause are logically
distinct. Condition (i) in Theorem~\ref{thmcharsci} requires
$D_{\mathbf{B}%
=\mathbf{1}, \mathbf{C}_2=\mathbf{c}_2^*}(\omega^*)=1$ for one
\mbox{assignment} $%
\mathbf{c}_2^*$ (and some $\omega^*$), while if $\mathbf{B}$ is a sufficient
cause (for $\omega^*$), then this condition is required to hold for all
assignments $\mathbf{c}_2^*$; in contrast condition (ii) in Theorem~\ref{thmcharsci} requires that there exists a single $\mathbf{c}_2^*$ (and
some $\omega^*$) such that for all $L \in\mathbf{B}$, $D_{\mathbf{B}%
\setminus\{L\}=\mathbf{1}, L=0, \mathbf{C}_2=\mathbf{c}_2^*}(\omega^*)=0$,
while for $\mathbf{B}$ to be a \textit{minimal} sufficient cause for
$\omega^*$
merely requires that for all $L \in\mathbf{B}$, there exists an
assignment $%
\mathbf{c}_2^L$ such that $D_{\mathbf{B}\setminus\{L\}=\mathbf{1}, L=0,
\mathbf{C}_2=\mathbf{c}_2^{L}}(\omega^*)=0$.

%
%ex3 #&#
\setcounter{excont}{0}
\begin{excont}[(Revisited)]%{{\hyperref[ex:2omegas]
Let $\mathbf{C}=\{X_1,X_2,X_3\}$, $\Omega= \{2\}$. Relative to $%
\mathbf{C}$, $\{X_1,X_2\}$ is a minimal sufficient cause for $\omega=2$
since $D_{111}(2)=D_{110}(2)=1$, and $D_{011}(2)=D_{100}(2)=0$. However
$%
\{X_1,X_2\}$ is not irreducible for ${\mathcal{D}}(\mathbf{C},\Omega)$
because we have $D_{101}(2) = D_{010}(2) = 1$, hence condition (ii) in
Theorem~\ref{thmcharsci} is not satisfied for either $X_3=0$, or $X_3=1$.
Conversely $\{X_1\}$ is irreducible for $\Omega= \{2\}$ since
$D_{111}(2)=1$%
, while $D_{011}(2)=0$, but $\{X_1\}$ is not a sufficient cause because
$%
D_{100}(2)=0$.
\end{excont}

Though irreducibility of $\mathbf{B}$ for ${\mathcal{D}}(\mathbf
{C},\Omega)$
neither implies, nor is implied by $\mathbf{B}$ being a minimal sufficient
cause for some $\omega\in\Omega$, it does imply that every sufficient
cause representation for ${\mathcal{D}}(\mathbf{C},\Omega)$ contains at
least one conjunction $\mathbf{B}_j$ of which~$\mathbf{B}$ is a (possibly
proper) subset. However, \textit{prima facie} this still leaves open the
possibility that, for example, every representation either includes
$\mathbf{%
B}\cup\{L\}$ or $\mathbf{B}\cup\{\overline L\}$, for some $L$, but no
representation includes both. However, this cannot occur:

%
%co3.6 #&#
\begin{corollary}
\label{corcontains} Let $\mathbf{C}=\mathbf{C}_1 \,\dot{\cup}\, \mathbf
{C}_2$, $%
\mathbf{B} \in\dot{\mathbb{P}}({\mathbb{L}}(\mathbf{C}_1))$, if
$\mathbf{B%
}$ is irreducible for ${\mathcal{D}}(\mathbf{C},\Omega)$ then there
exists a
set $\mathbf{B}^* \in\dot{\mathbb{P}}({\mathbb{L}}(\mathbf{C}))$, with $|\mathbf{B}^*|=|\mathbf{C}|$ such that in every
representation $(\mathbf{A},{\mathfrak{B}})$ for ${\mathcal{D}}(\mathbf
{C}%
,\Omega)$ there exists $\mathbf{B}_j \in{\mathfrak{B}}$, with $\mathbf{B}
\subseteq\mathbf{B}_j \subseteq\mathbf{B}^*$.
\end{corollary}

Thus irreducibility of $\mathbf{B}$ further implies that there is a set
$%
\mathbf{B}^*$ of size $|\mathbf{C}|$ such that in every representation there
is at least one conjunct containing $\mathbf{B}$ that is itself
contained in
$\mathbf{B}^*$. However, it should be noted that, in general, there may be
more than one conjunct $\mathbf{B}_j$ with $\mathbf{B} \subseteq\mathbf
{B}%
_j \subseteq\mathbf{B}^*$.

\begin{pf} Immediate from Theorem~\ref{thmcharsci},
taking $%
\mathbf{B}^*=\mathbf{B} \cup\mathbf{B}^{[\mathbf{c}_2^*]}$.
\end{pf}

Finally, we note that a conjunction that is both irreducible and a minimal
sufficient cause corresponds to an ``essential prime implicant'' in digital
circuit theory~\cite{marcovitz2001}. The Quine--McCluskey algorithm
\cite{quine1952,quine1955,mccluskey1956} finds the set of essential prime
implicants for a given Boolean function, which here corresponds to the
potential outcomes ${\mathcal{D}}(\mathbf{C},\omega)$ for an individual.

%s3.3 #&#
%s3.3 ###
\subsection{Enlarging the set of potential causes}

As noted in Section~\ref{secscms} a set $\mathbf{B}$ may be a minimal sufficient
cause for $\mathbf{C}$ but not a superset $\mathbf{C}^\prime$.
Irreducibility is also not preserved without further conditions. To state
these conditions that
preserve irreducibility we require the following:

%
%de3.7 #&#
\begin{defin}
\label{defintervar} $X^\prime$ is said to be \textit{not causally influenced}
by a set $\mathbf{C}$ if for all $\omega\in\Omega$, the potential outcomes
$X^\prime_{\mathbf{C} = \mathbf{c}}(\omega)$ are constant as $\mathbf{c}$
varies.
\end{defin}

We will also assume that if every $X^\prime\in\mathbf{C}^\prime$ is not
causally influenced by $\mathbf{C}$, then the following \textit{relativized consistency
axiom} holds:
%
%e7 #&#
%e3.2 ###
\begin{equation}
\label{eqrcons} D_{\mathbf{C}=\mathbf{c}, \mathbf{C}^\prime= \mathbf{C}^\prime(\omega
)}(\omega)= D_{\mathbf{C}=\mathbf{c}}(\omega),
\end{equation}
that is, that if variables in $\mathbf{C}^\prime$ are not causally
influenced by
the variables in $\mathbf{C}$, then the counterfactual value of $D$
intervening to set $\mathbf{C}$ to $\mathbf{c}$ is the same as the
counterfactual value of $D$ intervening to set $\mathbf{C}$ to $\mathbf{c}$
and the variables in $\mathbf{C}^\prime$ to the values they actually took
on.

We now have the following corollary to Theorem~\ref{thmcharsci}:

%
%co3.8 #&#
\begin{corollary}
\label{corexpand} Let $\mathbf{C}=\mathbf{C}_1 \,\dot{\cup}\, \mathbf
{C}_2$, $%
\mathbf{B} \in\dot{\mathbb{P}}({\mathbb{L}}(\mathbf{C}_1))$, $|\mathbf
{B}%
| = |\mathbf{C}_1|$. If $\mathbf{B}$ is irreducible for ${\mathcal{D}}(%
\mathbf{C},\Omega)$, $\mathbf{C}^\prime\cap\mathbf{C}=\varnothing$ and for
all $X^\prime\in\mathbf{C}^\prime$, $X^\prime$ is not causally influenced
by $\mathbf{C}$, then $\mathbf{B}$ is irreducible for ${\mathcal
{D}}(\mathbf{%
C}\cup\mathbf{C}^\prime,\Omega)$.
\end{corollary}

\begin{pf} By Theorem~\ref{thmcharsci} there exists $%
\omega^* \in\Omega$ and an assignment $\mathbf{c}_2^*$ to $\mathbf{C}_2$
such that (i) and (ii) hold. Let $\mathbf{c}^\prime= \mathbf{C}^\prime
(\omega^*)$. Since variables in $\mathbf{C}^\prime$ are not causally
influenced by $\mathbf{C}$, for all assignments $\mathbf{b}$,
\[
D_{\mathbf{B}=\mathbf{b}, \mathbf{C}_2=\mathbf{c}_2^*, \mathbf
{C}^\prime=
\mathbf{c}^\prime}\bigl(\omega^*\bigr) = D_{\mathbf{B}=\mathbf{b}, \mathbf{C}_2=%
\mathbf{c}_2^*, \mathbf{C}^\prime= \mathbf{C}^\prime(\omega^*)}\bigl(\omega^*\bigr) =
D_{\mathbf{B}=\mathbf{b}, \mathbf{C}_2=\mathbf{c}_2^*}\bigl(\omega^*\bigr);
\]
the second equality here follows from (\ref{eqrcons}). It follows that
$%
\omega^*$ and $(\mathbf{c}_2^*,\mathbf{c}^\prime)$ obey (i) and (ii) in
Theorem~\ref{thmcharsci} with respect to $\mathbf{C} \cup\mathbf{C}%
^\prime$.
\end{pf}

The assumption that every variable in $\mathbf{C}^\prime$ is not causally
influenced by $\mathbf{C}$, is required because otherwise we may have $%
\mathbf{C}^\prime(\omega^*) \neq(\mathbf{C}^\prime)_{\mathbf
{B}=\mathbf{b}%
^*} (\omega^*)$ for some assignment $\mathbf{b}^*$ to $\mathbf{B}$.
%Thus even though $(\omega^*, \mathbf{c}_2^*)$ satisfy (i) and (ii) in
%Theorem~\ref{thmcharsci},
%establishing that $\mathbf{B}$ is irreducible for ${\mathcal D}(
%$(\omega^*, \mathbf{c}_2^*, \mathbf{C}^\prime(\omega^*))$ may not
%satisfy these conditions.
For example, let $\mathbf{C} = \{ X_{1}, X_{2}, X_{3}\}$, and suppose that
\begin{eqnarray*}
D_{X_{1}=x_{1},X_{2}=x_{2},X_{3}=x_{3}}(\omega)&=& x_{3},
\\
(X_3)_{X_{1}=x_{1},X_{2}=x_{2}}(\omega)&=&x_{1}\wedge
x_{2}
\end{eqnarray*}
for all $\omega\in\Omega$. In this case $\{X_1,X_2\}$ is irreducible
for ${%
\mathcal{D}}(\{X_1,X_2\},\Omega)$, but not for ${\mathcal{D}}%
(\{X_1,X_2,X_3\},\Omega)$.
%Corollary~\ref{corexpand} fails because $X_3$ is causally influenced
%by $X_1$ and $X_2$.
We saw earlier that if $\mathbf{B}$ is a minimal sufficient cause for $%
\mathbf{C}$, then this does not imply that $\mathbf{B}$ is a minimal
sufficient cause with respect to subsets of $\mathbf{C}$. Here we see that
if $\mathbf{B}$ is irreducible with respect for ${\mathcal{D}}(\mathbf
{C}%
,\Omega)$, then this does not imply irreducibility for supersets
$\mathbf{C}%
^* \supset\mathbf{C}$, unless every variable in $\mathbf{C}^* \setminus
\mathbf{C}$ is not causally influenced by a variable in $\mathbf{C}$.

%s4 #&#
%s4 ###
\section{Tests for irreducibility}
\label{sectestsnway}

In this section we derive empirical conditions which imply that a given
conjunction $\mathbf{B}$ is irreducible for ${\mathcal{D}}(\mathbf{C}%
,\Omega) $. Our first approach is via condition (\ref{eqcond}).

%s4.1 #&#
%s4.1 ###
\subsection{Adjusting for measured confounders}

To detect that (\ref{eqcond}) holds requires us to identify the mean
of potential outcomes in certain subpopulations.
This is only possible if we have no unmeasured confounders \cite
{Rosenbaum1983,Robins1986}:

%
%de4.1 #&#
\begin{defin}
A set of covariates $\mathbf{W}$ \textit{suffices to adjust for
confounding of
(the effect of) $\mathbf{C}$ on $D$} if
%
%e8 #&#
%e4.1 ###
\begin{equation}
\label{eqnoconfound} D_{\mathbf{C}=\mathbf{c}}\perp\!\!\!\perp
\mathbf{C} \mid \mathbf{W} = \mathbf{w}
\end{equation}
for all $\mathbf{c}$, $\mathbf{w}$.
\end{defin}

%
%pr4.2 #&#
\begin{proposition}
\label{propadjconfound} If a set $\mathbf{W}$ suffices to adjust for
confounding of $\mathbf{C}$ on $D$ and $P(\mathbf{%
C} = \mathbf{c}, \mathbf{W}=\mathbf{w})>0$, then
\[
E[ D_{\mathbf{C}=\mathbf{c}} \mid\mathbf{W}=\mathbf{w}] = E[ D \mid \mathbf{%
C} = \mathbf{c}, \mathbf{W}=\mathbf{w}].
\]
\end{proposition}

The proof of this is standard and hence omitted.
%E[ D_{\mathbf{C}=\mathbf{c}} \mid\mathbf{W}=\mathbf{w}] &=& E[ D_{
%=\mathbf{c}} \mid\mathbf{C} = \mathbf{c}, \mathbf{W}=\mathbf{w}] \\
%&=& E[ D \mid\mathbf{C} = \mathbf{c}, \mathbf{W}=\mathbf{w}],
%where the first equality is by (\ref{eqnoconfound}), and the second
%follows
%by consistency.\hfill$\Box$

Note that if $\mathbf{W}$ is sufficient to adjust for confounding of $%
\mathbf{C}$ on $D$, then $\mathbf{W}$ is also sufficient to adjust for
confounding of $\mathbf{B}$ on $D$, where $\mathbf{B} \in\dot{\mathbb
{P}}(%
{\mathbb{L}}(\mathbf{C}))$, $|\mathbf{B}| = |\mathbf{C}|$.

%s4.2 #&#
%s4.2 ###
\subsection{Tests for irreducibility without monotonicity}
\label{secwithout}

%
%th4.3 #&#
\begin{theorem}
\label{thmsciconfound} Let $\mathbf{C}=\mathbf{C}_1\,\dot{\cup}\, \mathbf
{C}_2$%
, $\mathbf{B} \in\dot{\mathbb{P}}({\mathbb{L}}(\mathbf{C}_1))$,
$|\mathbf{%
B}| = |\mathbf{C}_1|$. If $\mathbf{W}$ is sufficient to adjust for
confounding of $\mathbf{C}$ on $D$, and for some $\mathbf{c}_2$, $\mathbf{w}$%
,
%
%e9 #&#
%e4.2 ###
\begin{eqnarray}\label{eqempcond}
0&<&E[D \mid\mathbf{B} = \mathbf{1}, \mathbf{C}_2 =
\mathbf{c}_2,%
\mathbf{W} = \mathbf{w}] ,
\nonumber
\\[-8pt]
\\[-8pt]
\nonumber
&&{} - \sum_{L \in\mathbf{B}} E\bigl[D \mid\mathbf{B}
\setminus\{L\} = \mathbf{1}, L = 0, \mathbf{C}_2 =
\mathbf{c}_2, \mathbf{W} = \mathbf{w%
}\bigr],
\end{eqnarray}
then $\mathbf{B}$ is irreducible for ${\mathcal{D}}(\mathbf{C} ,\Omega)$.
\end{theorem}

\begin{pf} We prove the contrapositive. Suppose that
$\mathbf{B%
}$ is not irreducible for ${\mathcal{D}}(\mathbf{C} ,\Omega)$. Then by
Theorem~\ref{thmcharsci}, for all $\omega\in\Omega$, and all
$\mathbf{c}%
_2$,
\[
D_{\mathbf{B} = \mathbf{1},\mathbf{C}_2=\mathbf{c}_2}(\omega) - \sum_{L
\in
\mathbf{B}}
D_{\mathbf{B}\setminus\{L\} = \mathbf{1}, L = 0,\mathbf
{C}_2=%
\mathbf{c}_2}(\omega) \leq 0.
\]
Hence for any $\mathbf{w}$,
\[
E \biggl[ D_{\mathbf{B} = \mathbf{1},\mathbf{C}_2=\mathbf{c}_2} - \sum_{L \in\mathbf{B}}
D_{\mathbf{B}\setminus\{L\} = \mathbf{1}, L = 0,
\mathbf{C}_2=\mathbf{c}_2} \Big| \mathbf{W} = \mathbf{w} \biggr] \leq 0.
\]
Applying Proposition~\ref{propadjconfound} to each term implies the
negation of (\ref{eqempcond}).
\end{pf}

The condition provided in Theorem~\ref{thmsciconfound} can be empirically
tested with $t$-test type statistics if $\mathbf{W}$ consists of a small
number of categorical or binary variables or using regression or inverse
probability of treatment weighting methods
\cite{Robins1999,vanderweelemarginal2010,vansteelandtmarginal2012,Vansteelandt2003}, otherwise.

%The condition provided in Theorem~\ref{thmsciconfound} constitutes an
%empirical
%test for a set $\mathbf{B}$ to be irreducible for some ${\mathcal D}(

It follows from Corollary~\ref{corexpand} that condition (\ref
{eqempcond}%
) further establishes that $\mathbf{B}$ is irreducible for ${\mathcal
{D}}(%
\mathbf{C}\cup\mathbf{C}^\prime,\Omega)$ so long as every variable in
$%
\mathbf{C}^\prime$ is not causally influenced by variables in $\mathbf{C}$.

It may be shown that condition (\ref{eqempcond}) is the sole restriction
on the law of $(D,{\mathbf C},{\mathbf W})$
implied by the negation of irreducibility.

%s4.3 #&#
%s4.3 ###
\subsection{Graphs}

In the next section we develop more powerful tests under monotonicity
assumptions. However, to state these conditions we first introduce some
concepts from graph theory:

%
%de4.4 #&#
\begin{defin}
A \textit{graph} ${\mathfrak{G}}$ defined on a set $\mathbf{B}$ is a
collection of pairs of elements in $\mathbf{B}$, ${\mathfrak{G}}\equiv
\{
\mathbf{E} \mid\mathbf{E} = \{B_1,B_2\} \subseteq\mathbf{B}, B_1\neq
B_2\}
$.
\end{defin}

This is the usual definition of a graph, except that the vertex
set here is a set of literals. We will refer to sets in ${\mathfrak
{G}}$ as
\textit{edges}, which we will represent pictorially as $B_1 \hspace*{3pt}\rule[2.5pt]{12pt}{0.3pt}\hspace*{3pt} B_2$.

%
%de4.5 #&#
\begin{defin}
Two elements $L, L^* \in\mathbf{B}$ are said to be \textit{connected}
in ${%
\mathfrak{G}}$ if there exists a sequence $L = L_1,\ldots,L_p = L^*$
of distinct elements in $\mathbf{B}$ such that $\{L_i,L_{i+1}\}\in{%
\mathfrak{G}}$ for $i=1,\ldots,p-1$.
\end{defin}

The sequence of edges joining $L$ and $L^*$ is said to form a
\textit{path} in ${\mathfrak{G}}$.

%
%de4.6 #&#
\begin{defin}
A graph ${\mathfrak{G}}$ on $\mathbf{B}$ is said to form a \textit
{tree} if $|{%
\mathfrak{G}}|=|\mathbf{B}|-1$, and any pair of distinct elements in $%
\mathbf{B}$ are connected in ${\mathfrak{G}}$.
\end{defin}

In a tree there is a unique path between any two elements.

%
%pr4.7 #&#
\begin{proposition}
\label{propbijection} Let ${\mathfrak{T}}$ be a tree on $\mathbf{B}$. For
each element $R\in\mathbf{B}$ there is a natural bijection%
\[
\phi^{\mathfrak{T}}_R\dvtx \mathbf{B} \setminus\{R\}
\leftrightarrow \mathfrak{T}
\]
given by $\phi^{\mathfrak{T}}_R(L) = \mathbf{E}=\{L^{\prime},L\}$ where
$%
\mathbf{E}\in\mathfrak{T}$ is the last edge on a path from $R$ to $L$.
\end{proposition}

It is not hard to show that for a graph $\mathfrak{G}$, if the bijections
described in Proposition~\ref{propbijection} exist, then $\mathfrak
{G}$ is a
tree.

%
%th4.8 #&#
\begin{theorem}[(Cayley~\cite{cayley1889})]
\label{thmcayley} On a set $\mathbf{B}$ there are $|\mathbf
{B}|^{|\mathbf{%
B|-2}}$ different trees.
\end{theorem}

%s4.4 #&#
%s4.4 ###
\subsection{Monotonicity}

Sometimes it may be known that a certain cause has an effect on an outcome
that is either always positive or always negative.

%
%de4.9 #&#
\begin{defin}
$B_i$ has a \textit{positive monotonic effect on $D$ relative to a set}
$%
\mathbf{B}$ (with $B_i \in\mathbf{B}$) in a population $\Omega$ if for
all $\omega\in\Omega$ and all values $\mathbf{b}_{-i}$ for the variables
in $\mathbf{B}\setminus\{B_i\}$, $D_{\mathbf{B}\setminus\{B_i\} =
\mathbf{b%
}_{-i}, B_i=1}(\omega)\geq D_{\mathbf{B}\setminus\{B_i\} = \mathbf{b}%
_{-i}, B_i=0}(\omega)$.
\end{defin}

Similarly we say that $L$ has a \textit{negative monotonic effect}
relative to $\mathbf{B}\cup\{L\}$ if $\overline{L}$ has a {positive
monotonic effect} relative to $\mathbf{B}\cup\{\overline{L}\}$. Note that
the case in which $D_{\mathbf{B}\setminus\{B_i\} = \mathbf{b}_{-i},
B_i=1}(\omega)= D_{\mathbf{B}\setminus\{B_i\} = \mathbf{b}_{-i},
B_i=0}(\omega)$ for all $\omega$, and hence $B_i$ has no effect on $D$
relative to $\mathbf{B}$, is included as a degenerate case.

The definition of a positive monotonic effect requires that an intervention
does not decrease $D$ for every individual, not simply on average,
regardless of
the other interventions taken.
%i.e. $E[D_{\mathbf{B}\setminus\{B_i\} = \mathbf{b%
%}_{-i}, B_i=1}]\geq E[D_{\mathbf{B}\setminus\{B_i\} = \mathbf{b}%
%_{-i}, B_i=0}]$.
This is thus a strong
assumption; %equirements for the attribution of a monotonic effect are
%thus
%considerable. However whenever a particular intervention is always
%beneficial or neutral for all individuals, one will be able to
%attribute a
%positive monotonic effect; whenever the intervention is always harmful
%or
%neutral for all individuals, one will be able to attribute a negative
%monotonic effect. \
see~\cite{VanderWeele2008c} for further discussion.

Monotonic Boolean functions have been studied in other contexts:

%
%pr4.10 #&#
\begin{proposition} If for all $C_i \in{\mathbf C}$, $C_i$ has a
(positive or negative) monotonic effect on $D$ relative to $\mathbf C$,
and $k= |{\mathbf C}|$, then the number of distinct sets of potential outcomes
in ${\mathcal{D}}(\mathbf{C}, \Omega)$
is given by the $k$th Dedekind number (Dedekind~\cite{dedekind1897},
Wiedemann~\cite{wiedemann1991}).
\end{proposition}

%s4.5 #&#
%s4.5 ###
\subsection{Tests for irreducibility with monotonicity}

Knowledge of the monotonicity of certain potential causes allows for the
construction of more powerful statistical tests for irreducibility than
those given by Theorem~\ref{thmsciconfound}.

%
%th4.11 #&#
\begin{theorem}
\label{thmmonotone-cond} Let $\mathbf{C}=\mathbf{C}_1\,\dot{\cup}\, \mathbf
{C}%
_2 $, $\mathbf{B} = (\mathbf{B}_+ \,\dot{\cup}\, \mathbf{B}^\prime) \in\dot
{%
\mathbb{P}}({\mathbb{L}}(\mathbf{C}_1))$, $|\mathbf{B}| = |\mathbf{C}_1|$
and suppose that every $L \in\mathbf{B}_+$ has a positive monotonic effect
on $D$ relative to $\mathbf{C}$. If for some tree ${\mathfrak{T}}$ on $%
\mathbf{B}_+$, $\omega^* \in\Omega$ and some $\mathbf{c}_2$, we have
%
%e10 #&#
%e4.3 ###
\begin{eqnarray}
\label{eqcondmonotone} 0&<& D_{\mathbf{B} = \mathbf{1},\mathbf{C}_2=\mathbf{c}_2}\bigl(\omega^*\bigr),
\nonumber
\\[-8pt]
\\[-8pt]
\nonumber
&& {} - \sum_{L \in\mathbf{B}} D_{\mathbf{B}\setminus\{L\} =
\mathbf{1%
}, L = 0,\mathbf{C}_2=\mathbf{c}_2}\bigl(
\omega^*\bigr) + \sum_{\mathbf{E} \in{%
\mathfrak{T}}} D_{\mathbf{B}\setminus\mathbf{E} = \mathbf{1}, \mathbf
{E} =
\mathbf{0},\mathbf{C}_2=\mathbf{c}_2}\bigl(
\omega^*\bigr),
\end{eqnarray}
then $\mathbf{B}$ is irreducible for ${\mathcal{D}}(\mathbf{C} ,\Omega)$.
\end{theorem}

If we know that $X$ has a \textit{negative} monotonic effect on $D$,
then we
may use this theorem to construct more powerful tests of the irreducibility
of sets containing $\overline{X}$ with respect to ${\mathcal{D}}(\mathbf{C}
,\Omega)$.
Under the assumption that every $L\in\mathbf{C}$ has a monotonic
effect on $D$, we
have shown via direct calculation using
\texttt{cddlib}~\cite{fukuda2005}
that for $|\mathbf{C}|\leq4$, the conditions in (\ref
{eqcondmonotone}) are the sole restrictions
on the law of $(D,{\mathbf C},{\mathbf W})$
implied by the negation of irreducibility. We conjecture that this
holds in general.
%However, knowledge that $X$ has a negative effect cannot be used
%to build tests for the irreducibility of sets containing $X$. There are
%other tests that could be formed by adding terms to (\ref{eqcond})
%however,
%these lead to weaker tests than (\ref{eqcondmonotone})

\begin{pf}
%Since (\ref{eqcondmonotone}) adds nonnegative
%terms to the RHS of (\ref{eqcond}), it is clear that (\ref{eqcond})
%implies (\ref{eqcondmonotone}). Thus
By Theorem~\ref{thmcharsci} it is
sufficient to prove that under the monotonicity assumption on $\mathbf
{B}_+$%
, (\ref{eqcondmonotone}) implies (\ref{eqcond}). Suppose that (\ref{eqcond}) does not hold, so that for all values $\mathbf{c}_2$, and
all $%
\omega^* \in\Omega$,
\[
D_{\mathbf{B} = \mathbf{1},\mathbf{C}_2=\mathbf{c}_2}\bigl(\omega^*\bigr) - \sum_{L
\in\mathbf{B}}
D_{\mathbf{B}\setminus\{L\} = \mathbf{1}, L = 0,\mathbf
{C}_2=%
\mathbf{c}_2}\bigl(\omega^*\bigr) \leq 0.
\]
Then for each $\omega^* \in\Omega$, there exists $R \in\mathbf{B}_+$ such
that
\[
D_{\mathbf{B} = \mathbf{1},\mathbf{C}_2=\mathbf{c}_2}\bigl(\omega^*\bigr) - \sum_{L
\in\mathbf{B}^\prime\cup\{R\}}
D_{\mathbf{B}\setminus\{L\} = \mathbf{1},
L = 0,\mathbf{C}_2=\mathbf{c}_2}\bigl(\omega^*\bigr) \leq0.
\]
For a given tree $\mathfrak{T}$, the remaining terms on the right-hand
side of (\ref{eqcondmonotone}) are
\begin{eqnarray*}
&&- \sum_{L \in\mathbf{B}_+\setminus\{R\}} D_{\mathbf
{B}\setminus\{L\} =
\mathbf{1}, L = 0,\mathbf{C}_2=\mathbf{c}_2}\bigl(
\omega^*\bigr) + \sum_{\mathbf
{E} \in{\mathfrak
T}} D_{\mathbf{B}\setminus\mathbf{E} = \mathbf{1}, \mathbf{E} = \mathbf
{0},\mathbf{C}_2=\mathbf{c}_2}\bigl(
\omega^*\bigr)\\
&&\qquad= \sum_{L \in\mathbf{B}_+\setminus\{R\}} \bigl( D_{\mathbf
{B}\setminus\phi_{R}^{%
\mathfrak{T}}(L) = \mathbf{1}, \phi_{R}^{\mathfrak{T}}(L) = \mathbf{0},%
\mathbf{C}_2=\mathbf{c}_2}\bigl(
\omega^*\bigr) - D_{\mathbf{B}\setminus\{L\} =
\mathbf{%
1}, L = 0,\mathbf{C}_2=\mathbf{c}_2}\bigl(\omega^*\bigr) \bigr),
\end{eqnarray*}
by Proposition~\ref{propbijection}. Finally since $\phi_{R}^{\mathfrak
{T}%
}(L) = \{L,L^\prime\} \subseteq\mathbf{B}_+$, $L^\prime$ has a positive
monotonic effect
%f1 #&#
%f1 ###
\begin{figure}

\includegraphics{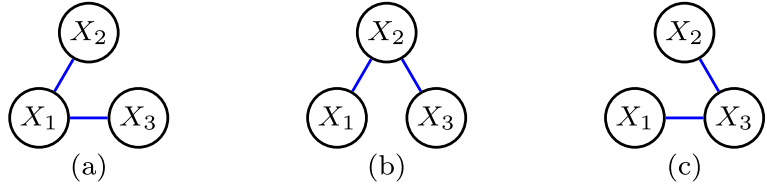}

\caption{The three trees on $\{X_1,X_2,X_3\}$.}\label{figone}
\end{figure}
on $D$ relative to $\mathbf{C}$, thus no term in the sum
is positive. Consequently for all $\omega^* \in\Omega$, (\ref{eqcondmonotone})
does not hold for any tree $\mathfrak{T}$.
\end{pf}

%A quick query for you:
%
%At the end of the proof of Theorem~4.10, p.16, line -2,
%am I correct in thinking that we should say:
%"thus each term in the nonpositive",
%or "no term in the sum is positive",
%
%rather than "each term is negative", as we have now.
%
%The last seems wrong to me, since positive monotonicity
%allows both terms to be zero, or both 1.
%
%?Or am I forgetting something?
%
%While I'm at it, I have another query.
%Do you agree that we should add at line -6 on p.16, to:
%
%"The remaining terms on the RHS of (4.3) are:"
%
%"For a given tree T, the remaining terms on the RHS of (4.3) are:"
%
%and then on the last line of the proof:
%
%"Consequently for all omega* in Omega, (4.3) does not hold for any
%tree T."
%
%??

%effects on }$D$\textit{\ relative to }$X_{1},,\ldots,,X_{m}$\textit{\ and
%that
%none of }$X_{m+1},,\ldots,,X_{s}$\textit{\ are intermediate variables
%between }$%
%X_{1},,\ldots,,X_{m}$\textit{\ and }$D$\textit{. \ Let }$\mathcal{U}%
%=\{(x_{1},,\ldots,,x_{m})\in\{0,1\}^{m}:\dsum_{i=1}^{m}x_{i}=m-1\}$%
%D_{X_{1}=1,\ldots,X_{m}=1}(\omega)-\dsum_{(x_{1},,\ldots,,x_{m})\in\mathcal{%
%U}}D_{x_{1},\ldots,x_{m}}(\omega)+\dsum_{(x_{1},,\ldots,,x_{m})\in\mathcal{S%
%}}D_{x_{1},\ldots,x_{m}}(\omega)>0$\textit{\ for some subordinate set }$S$%
%X_{1},,\ldots,,X_{m}$\textit{\ have a sufficient cause interaction.}

%
%ex2 #&#
\begin{example}
In the case $\mathbf{B} = \{X_1,X_2\} = \mathbf{B}_+=\mathbf{C}$,
there is only one tree on $\mathbf{B}_+$, consisting of a single
edge $%
X_1\hspace*{3pt}\rule[2.5pt]{12pt}{0.3pt}\hspace*{3pt} X_2$. Thus if $X_1$ and $X_2$ have a
positive monotonic effect on $D$ (relative to $\mathbf{C}$) then
Theorem~\ref{thmmonotone-cond} implies that if the following inequality holds for
some $%
\omega\in\Omega$,
\[
D_{11}(\omega)- \bigl(D_{10}(\omega)+D_{01}(
\omega) \bigr)+D_{00}(\omega )>0,
\]
then $\{X_1,X_2\}$ is irreducible for ${\mathcal{D}}(\mathbf{C}
,\Omega)$.
\end{example}%

%% \vspace{-0.5cm}
%{\scriptsize\ \psset{linewidth=0.8pt} \newlength{\MyLength} %
%% \newcommand{\mynewNode}[2]{\circlenode[linestyle=none]{#1}{
%0){%
%% \ncline[linecolor=blue]{->}{0}{1}
%% \end{pdfpic}

%
%ex3 #&#
\begin{example}
If $\mathbf{B} = \{X_1,X_2,X_3\} = \mathbf{B}_+=\mathbf{C}$, then
there are three trees on $\mathbf{B}_+$; see Figure~\ref{figone}. These
correspond to the following conditions:
\begin{eqnarray*}
&&\mathrm{(a)}\quad D_{111}(\omega)- \bigl( D_{110}(\omega)+D_{101}(
\omega )+D_{011}(\omega) \bigr)+ \bigl(D_{010}(
\omega)+D_{001}(\omega) \bigr) > 0,
\\
&&\mathrm{(b)}\quad D_{111}(\omega)- \bigl( D_{110}(\omega)+D_{101}(
\omega )+D_{011}(\omega) \bigr)+ \bigl(D_{100}(
\omega)+D_{001}(\omega) \bigr) > 0,
\\
&&\mathrm{(c)}\quad D_{111}(\omega)- \bigl(D_{110}(\omega)+D_{101}(
\omega )+D_{011}(\omega) \bigr)+ \bigl(D_{100}(
\omega)+D_{010}(\omega) \bigr) > 0.
\end{eqnarray*}
Thus $\mathbf{B}$ is irreducible for ${\mathcal{D}}(\mathbf{C} ,\Omega
)$ if
at least one holds for some $\omega\in\Omega$.
\end{example}

%
%co4.12 #&#
\begin{corollary}
\label{cormonotone-cond-empirical} Let $\mathbf{C}=\mathbf{C}_1\,\dot
{\cup}\,
\mathbf{C}_2$, $\mathbf{B} = (\mathbf{B}_+ \,\dot{\cup}\, \mathbf{B}^\prime)
\in\dot{\mathbb{P}}({\mathbb{L}}(\mathbf{C}_1))$, $|\mathbf{B}| = |%
\mathbf{C}_1|$. Suppose that every $L \in\mathbf{B}_+$ has a positive
monotonic effect on $D$ relative to $\mathbf{C}$, and $\mathbf{W}$ is
sufficient to adjust for confounding of $\mathbf{C}$ on $D$. If for some
tree ${\mathfrak{T}}$ on $\mathbf{B}_+$, and some $\mathbf{c}_2$,
$\mathbf{w}
$ we have
%
%e11 #&#
%e4.4 ###
\begin{eqnarray}\label{eqcondmonotoneempirical}
0&<& E[D \mid\mathbf{B} = \mathbf{1}, \mathbf{C}_2 =
\mathbf{c}_2, \mathbf{W} = \mathbf{w}]
\nonumber
\\
&&{} - \sum_{L \in\mathbf{B}} E\bigl[D\mid\mathbf{B}
\setminus\{L\} = %
\mathbf{1}, L = 0, \mathbf{C}_2 =
\mathbf{c}_2, \mathbf{W} = \mathbf{w%
}\bigr]
\\
&&{}  + \sum_{\mathbf{E} \in{\mathfrak{T}}} E[D \mid
\mathbf{B}%
\setminus\mathbf{E} = \mathbf{1}, \mathbf{E} = \mathbf{0},
\mathbf{C%
}_2 = \mathbf{c}_2, \mathbf{W} =
\mathbf{w}],
\nonumber
\end{eqnarray}
then $\mathbf{B}$ is irreducible for ${\mathcal{D}}(\mathbf{C} ,\Omega)$.
\end{corollary}

\begin{pf} Directly analogous to the proof of Theorem~\ref
{thmsciconfound}.
\end{pf}

The special case of the previous Corollary where $|\mathbf
{B_+}|=|\mathbf{C}%
|=2$, and $\mathbf{W}=\varnothing$, appears in Rothman and Greenland \cite
{Rothman1998}; see also Koopman~\cite{Koopman1981}.
Theorem~\ref{thmcayley}
implies that if every literal in $\mathbf{B}$ has a positive monotonic
effect on $D$, then we will have $|\mathbf{B}|^{|\mathbf{B}|-2}$ conditions
to test, each of which is sufficient to establish the irreducibility of
$%
\mathbf{B}$ for ${\mathcal{D}}(\mathbf{C} ,\Omega)$. As before, the
conditions (\ref{eqcondmonotoneempirical}) may be tested via $t$-test
type statistics or using various statistical models.

As with the results in Section~\ref{secwithout}, we may apply Corollary~\ref
{corexpand} to establish that~$\mathbf{B}$ is irreducible for
${\mathcal{D}}%
(\mathbf{C}\cup\mathbf{C}^\prime,\Omega)$ if every variable in
$\mathbf{C}%
^\prime$ is not causally influenced by variables in $\mathbf{C}$.

%s4.6 #&#
%s4.6 ###
\subsection{Tests for a minimal sufficient cause under monotonicity}

As noted in Section~\ref{secbequalsc} if $|\mathbf{B}|=|\mathbf{C}|$,
then irreducible conjunctions are also minimal sufficient causes. Thus in
this special case, the tests of irreducibility given in Theorem~\ref{thmsciconfound}
and Corollary~\ref{cormonotone-cond-empirical} also
establish that $\mathbf{B}$ is a minimal sufficient cause relative to~$\mathbf{C}$.
When $|\mathbf{B}|<|\mathbf{C}|$ these tests do not in general
establish this. However, under positive monotonicity assumptions on
$\mathbf{%
C}_2$, such tests may be obtained by taking $\mathbf{c}_2=\mathbf{0}$:

%
%pr4.13 #&#
\begin{proposition}
\label{propmsccond} Let $\mathbf{C}=\mathbf{C}_1 \,\dot{\cup}\, \mathbf{C}_2$,
$\mathbf{B} \in\dot{\mathbb{P}}({\mathbb{L}}(\mathbf{C}_1))$, $|\mathbf
{B}%
| = |\mathbf{C}_1|$. Suppose every $L\in\mathbf{C}_2$ has a positive
monotonic effect on $D$ relative to $\mathbf{C}$. If \textup{(i)} $D_{%
\mathbf{B} = \mathbf{1}, \mathbf{C}_2=\mathbf{0}}(\omega^*)=1$ and
\textup{(ii)} for all $L \in\mathbf{B}$, $D_{\mathbf{B}\setminus\{L\} = \mathbf{1},
L = 0, \mathbf{C}_2=\mathbf{0}}(\omega^*)=0$, then $\mathbf{B}$ is a minimal
sufficient cause for $D$ relative to $\mathbf{C}$ for $\omega^*$.
\end{proposition}

\begin{pf} For any $\mathbf{c}_2$, $D_{\mathbf{B} =
\mathbf{1}%
, \mathbf{C}_2=\mathbf{c}_2}(\omega^*)\geq D_{\mathbf{B} = \mathbf{1},
\mathbf{C}_2=\mathbf{0}}(\omega^*)=1$ by the monotonicity assumption.
Hence $%
\mathbf{B}$ is a sufficient cause for $D$ relative to $\mathbf{C}$ for $
\omega^*$. Minimality follows directly from condition (ii).
\end{pf}

We have the following corollaries which provide conditions under which $
\mathbf{B}$ is a minimal sufficient cause for $D$ relative to $\mathbf{C}$
for some $\omega\in\Omega$, in addition to being irreducible for
$\mathcal{D%
}(\mathbf{C},\Omega)$:

%
%co4.14 #&#
\begin{corollary}
\label{cormsctest} Let $\mathbf{C}=\mathbf{C}_1 \,\dot{\cup}\, \mathbf
{C}_2$, $%
\mathbf{B} \in\dot{\mathbb{P}}({\mathbb{L}}(\mathbf{C}_1))$, $|\mathbf
{B}%
| = |\mathbf{C}_1|$. Suppose every $L\in\mathbf{C}_2$ has a positive
monotonic effect on $D$ relative to~$\mathbf{C}$, and $\mathbf{W}$ is
sufficient to adjust for confounding of $\mathbf{C}$ on $D$. If (\ref{eqempcond})
holds with $\mathbf{c}_2=\mathbf{0}$ for some $\mathbf{w}$,
then $\mathbf{B}$ is a minimal sufficient cause of $D$ relative to
$\mathbf{C%
}$ for some $\omega\in\Omega$.
\end{corollary}

\begin{pf} The proof follows immediately from
Proposition~\ref{propmsccond} and
Theorem~\ref{thmsciconfound}.
\end{pf}

%
%co4.15 #&#
\begin{corollary}
\label{cormsc-monotone-cond-empirical} Let $\mathbf{C}=\mathbf{C}_1\,\dot
{\cup%
}\, \mathbf{C}_2$, $\mathbf{B} = (\mathbf{B}_+ \,\dot{\cup}\, \mathbf
{B}^\prime)
\in\dot{\mathbb{P}}({\mathbb{L}}(\mathbf{C}_1))$, $|\mathbf{B}| = |%
\mathbf{C}_1|$. Suppose that every $L \in\mathbf{B}_+\cup\mathbf
{C}_2$ has
a positive monotonic effect on $D$ relative to~$\mathbf{C}$, and
$\mathbf{W}$
is sufficient to adjust for confounding of $\mathbf{C}$ on $D$. If (%
\ref{eqcondmonotoneempirical}) holds with $\mathbf{c}_2=\mathbf{0}$ for
some $\mathbf{w}$ and some tree ${\mathfrak{T}}$ on $\mathbf{B}_+$,
then $%
\mathbf{B}$ is a minimal sufficient cause of $D$ relative to $\mathbf{C}$
for some $\omega\in\Omega$.
\end{corollary}

\begin{pf} The proof follows immediately from Proposition~\ref{propmsccond} and
Corollary~\ref{cormonotone-cond-empirical}.
\end{pf}

%Note that for the the previous
%Theorem and Corollary
%we don't actually need the assumption
%of a positive monotonic effect for
%elements in B_1 for *any* setting
%of B_1, B_2, we only need that each L in B
%is positive monotonic for B_2=1, and B_1\{L} with at most one zero

%s5 #&#
%s5 ###
\section{Singular interactions}
\label{secsingular}

In the genetics literature, in the context of two binary genetic
factors, $%
X_{1}$ and $X_{2}$, ``compositional'' \textit{epistasis}
\cite{Bateson1909,Cordell2002,Phillips2008} is said to be present if for
some $\omega^{\ast}$, $D_{11}(\omega^{\ast})=1$ but
$D_{10}(\omega^{\ast})=D_{01}(\omega^{\ast})=D_{00}(\omega^{\ast})=0$; in this case
the effect of one genetic factor is effectively masked when the other
genetic factor is absent. If $\{X_1,X_2\}$ is a minimal sufficient
cause of $%
D$ relative to $\{X_1,X_2\}$ for $\omega^*$ then although this implies $
D_{11}(\omega^{\ast})=1$ and $D_{10}(\omega^{\ast})=D_{01}(\omega^{\ast
})=0$, it does \textit{not} imply $D_{00}(\omega^{\ast})=0$. This motivates
the following:

%
%de5.1 #&#
\begin{defin}
A minimal sufficient cause $\mathbf{B}$ for $D$ relative to $\mathbf
{C}$ for
$\omega^*$ is said to be \textit{singular} if there is no $\mathbf
{B}^\prime
\in\dot{\mathbb{P}}(\mathbb{L}(\mathbf{C}))$, $\mathbf{B}^\prime\neq
\mathbf{B}$, forming a minimal sufficient cause for $D$ relative to
$\mathbf{%
C}$ for $\omega^*$. $\mathbf{B}$ is \textit{singular for }$\mathcal
{D}(\mathbf{%
C},\Omega)$ if $\mathbf{B}$ is singular relative to $\mathbf{C}$ for
some $%
\omega^* \in\Omega$. %
%Let $\mathbf{B}\in\dot{\mathbb P}(\mathbb{L}(%
%for $\mathcal{D}(\mathbf{C},\mathbf{\Omega})$ if there exists $\omega
%^{\ast}\in\mathbf{\Omega}$ such that $(1,\{\mathbf{B}\})$
%constitutes a
%sufficient cause representation for $\mathcal{D}(\mathbf{C},\{ \omega^{
\end{defin}

If $\mathbf{B}$ is singular for ${%
\mathcal{D}}(\mathbf{C},\Omega)$, then we will also refer to a
\textit{singular interaction} between the components of $\mathbf{B}$.
We now characterize singularity in terms of potential outcomes:

%
%th5.2 #&#
\begin{theorem}
\label{thmsingularalt} Let $\mathbf{C}=\mathbf{C}_{1}%
\cup \mathbf{C}_{2}$, $\mathbf{B}\in \dot{\mathbb{P}}(\mathbb
{L}(\mathbf{%
C}_{1}))$, $|\mathbf{B}|=|\mathbf{C}_{1}|$. Then $\mathbf{B}$ is
singular for
$\mathcal{D}(\mathbf{C},{\Omega})$ if and only if there exists $\omega^{\ast}\in{%
\Omega}$ such that
%
%e12 #&#
%e5.1 ###
\begin{equation}
\label{eqsingular-char} \mbox{for all values }\mathbf{c}_{2}^{\ast},
\mathbf{b}\dvtx\qquad
D_{\mathbf{B}=\mathbf{b },\mathbf{C}_{2}=\mathbf{c}_{2}^{\ast
}}\bigl(\omega^{\ast}
\bigr)=1 \quad\Longleftrightarrow\quad \mathbf{b}=\mathbf{1}.
\end{equation}
%
%{\rm(ii)} for all $\mathbf{B}^{\ast}\in
% \dot{\mathbb P}(\mathbb{L}(\mathbf{C}_{1}))$, $|\mathbf{B}^{\ast}|=|
%_{1}|$ with $\mathbf{B}^{\ast}\neq\mathbf{B}$, $D_{\mathbf{B}^{\ast
%}=%
\end{theorem}

Note that (\ref{eqsingular-char}) is equivalent to
%
%e13 #&#
%e5.2 ###
\begin{equation}
\label{eqsingular-char-alt} D_{\mathbf{C}=\mathbf{c}}\bigl(\omega^{\ast}\bigr)=
\Bigl(\bigwedge (\mathbf {B})\Bigr)_{\mathbf{c}%
} \qquad\mbox{for
all }\mathbf{c}.
\end{equation}
Thus if $\mathbf{B}$ is singular for $\mathcal{D}(\mathbf{C},\Omega)$, then
there is some individual $\omega^*$ whose potential outcomes are given by
the single conjunction $\mathbf{B}$.

\begin{pf} By definition, $\mathbf{B}$ is a sufficient cause
for $D$ for $\omega^*$ if and only if ($\mathbf{b} = \mathbf{1}$ $
\Rightarrow$ $%
D_{\mathbf{B}=\mathbf{b}, \mathbf{C}_2=\mathbf{c}_2^*}(\omega^*) = 1$).
Thus it is sufficient to show that, assuming $\mathbf{B}$ is a minimal
sufficient cause for $D$ for $\omega^*$, there are no other minimal
sufficient causes of $D$ for $\omega^*$ if and only if ($D_{\mathbf
{B}=\mathbf{b},
\mathbf{C}_2=\mathbf{c}_2^*}(\omega^*) = 1$ $ \Rightarrow$ $\mathbf{b}%
= \mathbf{1}$).\vspace*{1pt}

Suppose $\mathbf{B}$ is the only minimal sufficient cause for $D$ for $%
\omega^*$, but that for some $\mathbf{b}^* \neq \mathbf{1}$,
$D_{\mathbf{%
B}=\mathbf{b}^*, \mathbf{C}_2=\mathbf{c}_2^*}(\omega^*)=1$. Let $\mathbf
{B}%
^\dag\equiv\mathbf{B}^{[\mathbf{B}=\mathbf{b}^*,\mathbf{C}_2=\mathbf
{c}%
_2^*]}$. $\mathbf{B}^\dag$ forms a sufficient cause for $D$ for $\omega^*$,
and $\mathbf{B} \nsubseteq\mathbf{B}^\dag$. Hence there is some
$\mathbf{B}%
^\prime\subseteq\mathbf{B}^\dag$ that is a minimal sufficient cause
for $D
$ for $\omega^*$, and $\mathbf{B}\neq\mathbf{B}^\dag$, a
contradiction.

Conversely suppose ($D_{\mathbf{B}=\mathbf{b}, \mathbf{C}_2=\mathbf{c}%
_2^*}(\omega^*) = 1$ $ \Rightarrow$ $\mathbf{b} = \mathbf{1}$) but
there exists another minimal sufficient cause $\mathbf{B}^\prime$ for $D$
for $\omega^*$, and $\mathbf{B}\neq\mathbf{B}^\dag$. Since $\mathbf{B}%
^\prime$ is minimal, $\mathbf{B} \nsubseteq\mathbf{B}^\prime$. Thus there
exists a $\tilde{\mathbf{c}}$ such that $(\mathbf{B})_{\tilde{\mathbf
{c}}%
}\neq\mathbf{1}$, but $(\mathbf{B}^\prime)_{\tilde{\mathbf{c}}}=\mathbf{1}$
and hence $D_{\mathbf{C}={\tilde{\mathbf{c}}}}(\omega^*)=1$, a
contradiction.
\end{pf}

%
%co5.3 #&#
\begin{corollary}
\label{corsingularimplies} For $\mathcal{D}(\mathbf{C},{\Omega})$,
if $%
\mathbf{B}$ is singular then $\mathbf{B}$ is irreducible.
\end{corollary}

\begin{pf} The proof follows immediately from (\ref
{eqsingular-char}) and the
definition of irreducibility.
\end{pf}

%Corollary~\ref{corexpand} showed that irreducibility was retained
% if the set of potential causes $\mathbf{C}$ were
%enlarged so long as the set of additional potential causes $\mathbf{C}%
%^{\prime}$ were not causally influenced by $\mathbf{C}$. \ The next
%corollary gives the analogous result for singular interactions. \ The
%proof
%is essentially identical that of Corollary~16 and therefore omitted.

%Let $\mathbf{C}=\mathbf{C}_{1}\dbigcup\mathbf{C}_{2}$%
%, $\mathbf{B}\in\dot{\mathbb P}(\mathbb{L}(\mathbf{C}_{1}))$, $|
%|=|\mathbf{C}_{1}|$. If $\mathbf{B}$ exhibits a singular interaction
%for $%
%^{\prime}=\varnothing$ and for all $X^{\prime}\in\mathbf{C}^{\prime
%}$, $%
%X^{\prime}$ is not causally influenced by $\mathbf{C}$, then $
%exhibits a singular interaction for $\mathcal{D}(\mathbf{C}\cup
%^{\prime},\mathbf{\Omega})$.

Theorem~\ref{thmsingularchar} relates singular interactions to properties
of the set of sufficient cause representations for $\mathcal{D}(\mathbf
{C},{%
\Omega})$. %First we prove the following:

%Let ${\mathfrak B} \subseteq\dot{\mathbb P}(\mathbb{L}(\mathbf{C}))
%$(\dnf({\mathfrak B}))_\mathbf{c} = (\wedge(\mathbf{B}))_\mathbf{c}$
%iff {\rm(i)} for all $%
%{\rm(ii)} if $\mathbf{B}_{i}\in\mathfrak{B}$ then $\mathbf{B}

%minimal subset of $\mathbf{C}$ such that
%${\mathfrak B} \subseteq\dot{\mathbb P}(\mathbb{L}(\mathbf{C}_1))$.
%Let
%$\mathbf{C}_2 = \mathbf{C}\setminus\mathbf{C}_1$. The proof is by
%induction on $k\equiv|\mathbf{C}_1|-|\mathbf{B}|$.
% If $k=0$ then the result holds trivially because (i) and (ii) imply $
%For the inductive case let $L\in\mathbf{C}_1$ be a literal such that
%both $L,\bar{L} \notin\mathbf{B}$.
%Let ${\mathfrak B}^\prime\equiv\{ \mathbf{B}_i\setminus\{L,\bar{L}\}
%Since, by hypothesis, $\mathfrak B$ obeys (i) and (ii), it follows
%that ${\mathfrak B}^\prime$ also does.
%Hence for all $\bf c$, $(\dnf({\mathfrak B}^\prime))_\mathbf{c} = (
%$\mathbf{B}_i = \mathbf{B}^\prime\cup\{X\}$, where $X \in\{L,\bar{L}
%then $\mathbf{B}_i $ is not minimal in ${\mathfrak B}$, so may be
%removed without changing
% $(\dnf({\mathfrak B}^\prime))_\mathbf{c}$. If $\mathbf{B}^\prime
% there exists, $\tilde\mathbf{B}_i = \mathbf{B}^\prime\cup\{ \bar{X}

%
%th5.4 #&#
\begin{theorem}
\label{thmsingularchar} Let $\mathbf{B}\in\dot{\mathbb{P}}(\mathbb
{L}(%
\mathbf{C}))$. $\mathbf{B}$ is singular for $\mathcal{D}(\mathbf
{C},{\Omega}%
)$ if and only if there exists $\omega^{\ast}\in{\Omega}$ such that
in every
representation $(\mathbf{A},\mathfrak{B})$ for $\mathcal{D}(\mathbf{C},{
\Omega})$, \textup{(i)} for all $\mathbf{B}^{\ast}\in\dot{\mathbb{P}}(%
\mathbb{L}(\mathbf{C}))$, with $|\mathbf{B}^{\ast}|=|\mathbf{C}|$ and $
\mathbf{B}\subseteq\mathbf{B}^{\ast}$ there exists $\mathbf{B}_{i}\in
\mathfrak{B}$ with $\mathbf{B}_{i}\subseteq\mathbf{B}^{\ast}$ and $%
A_{i}(\omega^{\ast})=1$; \textup{(ii)} for all $\mathbf{B}_{i}\in
\mathfrak{B}$ such that $\mathbf{B}\nsubseteq\mathbf{B}_{i}$,
\mbox{$A_{i}(\omega^{\ast})=0$}.
\end{theorem}

\begin{pf} Let $\mathbf{C}=\mathbf{C}_1\,\dot{\cup}\, \mathbf{C}_2$,
where $%
\mathbf{B}\in\dot{\mathbb{P}}(\mathbb{L}(\mathbf{C}_1))$, and $|\mathbf{B}|
= |\mathbf{C}_1|$.

$(\Rightarrow)$ Suppose $\mathbf{B}$ is singular for $\mathcal
{D}(\mathbf{C}%
,{\Omega})$, so that some $\omega^{\ast}\in{\Omega}$ satisfies~(\ref{eqsingular-char-alt}). Then for any representation $(\mathbf
{A},\mathfrak{B%
})$ for $\mathcal{D}(\mathbf{C},{\Omega})$ and any $\mathbf{B}^{\ast}$
such that $|\mathbf{B}^{\ast}|=|\mathbf{C}|$ and $\mathbf{B}\subseteq
\mathbf{B}^{\ast}$, we can select values $\mathbf{c}_{2}^{\ast}$ so
that $%
\mathbf{B}^{\ast}=\mathbf{B}^{[\mathbf{B}=\mathbf{1},\mathbf
{C}_{2}=\mathbf{%
c}_{2}^{\ast}]}$. Since $D_{\mathbf{B}=\mathbf{1},\mathbf
{C}_{2}=\mathbf{c}%
_{2}^{\ast}}(\omega^{\ast})=1$ there exists $A_{i}\in\mathbf{A}$, $%
\mathbf{B}_{i}\in\mathfrak{B}$ with $A_{i}(\omega^{\ast})=1$ and
$(\bigwedge (%
\mathbf{B}_{i}))_{\mathbf{B}=\mathbf{1},\mathbf{C}_{2}=\mathbf
{c}_{2}^{\ast
}}=1$. Thus $\mathbf{B}_{i}\subseteq\mathbf{B}^{\ast}$, so (i) holds. For
all $\mathbf{B}_{i}\in\mathfrak{B}$ such that $\mathbf{B}\nsubseteq
\mathbf{%
B}_{i}$, we can choose $\widetilde{\mathbf{B}}\in\dot{\mathbb{P}}(\mathbb
{L}(%
\mathbf{C}_{1}))$, $|\widetilde{\mathbf{B}}|=|\mathbf{C}_{1}|$ with $\widetilde{
\mathbf{B}}\neq\mathbf{B}$ and values $\tilde{\mathbf{c}}_{2}$ so that
$%
\mathbf{B}_{i}\subseteq\mathbf{B}^{[\widetilde{\mathbf{B}}=\mathbf
{1},\mathbf{C}%
_{2}=\tilde{\mathbf{c}}_{2}]}$. Since $D_{\widetilde{\mathbf{B}}=\mathbf{1},
\mathbf{C}_{2}=\tilde{\mathbf{c}}_{2}}(\omega^{\ast})=0$ we have $%
A_{i}(\omega^{\ast})=0$ since $(\bigwedge  (\mathbf
{B}_{i}))_{\widetilde{\mathbf{B}%
}=\mathbf{1},\mathbf{C}_{2}=\tilde{\mathbf{c}}_{2}}=1$, so (ii) holds as
required.

$(\Leftarrow)$ Suppose there exists $\omega^{\ast}\in{\Omega} $ such
that every representation $(\mathbf{A},\mathfrak{B})$ satisfies (i) and
(ii). %for $\mathcal{D}(\mathbf{C},\mathbf{\Omega})$, for all $%
%with $\mathbf{B}_{i}\subseteq\mathbf{B}^{\ast}$ and $A_{i}(\omega^{
%})=1$ and for all $\mathbf{B}_{i}\in\mathfrak{B}$ such that $
We will show that (\ref{eqsingular-char}) holds.
%for all values $\mathbf{c}%
%_{2}^{\ast}$ for $\mathbf{C}_{2}$: (i) $D_{\mathbf{B}=\mathbf{1},
%_{2}=\mathbf{c}_{2}^{\ast}}(\omega^{\ast})=1$; (ii) for all $
%^{\ast}\in\dot{\mathbb P}(\mathbb{L}(\mathbf{C}_{1}))$, $|\mathbf{B}%
%^{\ast}|=|\mathbf{C}_{1}|$ with $\mathbf{B}^{\ast}\neq\mathbf{B}$,
%$D_{%
%^{\ast})=0$. \
For any values $\mathbf{c}_{2}^{\ast}$ let $\mathbf{B}^{\ast}\equiv
\mathbf{B}^{[\mathbf{B}=\mathbf{1},\mathbf{C}_{2}=\mathbf{c}_{2}^{\ast}]}$,
so $|\mathbf{B}^{\ast}| = |\mathbf{C}|$ and $\mathbf{B}\subseteq
\mathbf{B%
}^{\ast}$. Thus by (i) there exists $\mathbf{B}_{i}\in\mathfrak{B}$
with $%
\mathbf{B}_{i}\subseteq\mathbf{B}^{\ast}$ and $A_{i}(\omega^{\ast})=1$.
Hence $D_{\mathbf{B}=\mathbf{1},\mathbf{C}_{2}=\mathbf{c}_{2}^{\ast
}}(\omega^{\ast})=1$ since $A_{i}(\omega^{\ast})=1$ and $(\bigwedge  (%
\mathbf{B}_{i}))_{\mathbf{B}=\mathbf{1},\mathbf{C}_{2}=\mathbf
{c}_{2}^{\ast
}}=1$. Conversely for any $\mathbf{b}^\prime \neq \mathbf{1}$, let $%
\mathbf{B}^\prime\equiv\mathbf{B}^{[\mathbf{B}=\mathbf{b}^\prime]}$,
so $|%
\mathbf{B}^{\prime}|=|\mathbf{C}_{1}|$ with $\mathbf{B}^{\prime}\neq
\mathbf{B}$. Thus for all $\mathbf{B}_{i}\in\mathfrak{B}$ such that $%
(\bigwedge  (\mathbf{B}_{i}))_{\mathbf{B}^{\prime}=\mathbf
{1},\mathbf{C}_{2}=%
\mathbf{c}_{2}^{\ast}}=1$, $\mathbf{B}\nsubseteq\mathbf{B}_{i}$ and thus
by (ii) $A_{i}(\omega^{\ast})=0$. Hence $D_{\mathbf{B}^{\prime
}=\mathbf{1}%
,\mathbf{C}_{2}=\mathbf{c}_{2}^{\ast}}(\omega^{\ast})=0$.
\end{pf}

We now consider results relevant for testing for singular interactions with
or without monotonicity assumptions.

%
%th5.5 #&#
\begin{theorem}
\label{thmsingular-test} Let $\mathbf{B}=\mathbf{B}_{+}\,\dot{\cup
}\,\mathbf{B}%
^{\prime}\in\dot{\mathbb{P}}(\mathbb{L}(\mathbf{C}))$, $|\mathbf{B}|=|
\mathbf{C}|$ and suppose that every $L\in\mathbf{B}_{+}$ has a positive
monotonic effect on $D$ relative\vadjust{\goodbreak} to $\mathbf{C}$.  If for some tree $%
\mathfrak{T}$ on $\mathbf{B}_{+}$ and some $\omega^{\ast}\in{\Omega
}$, we
have%
%
%e14 #&#
%e5.3 ###
\begin{eqnarray}\label{eqfourterms}
&&D_{\mathbf{B}=\mathbf{1}}\bigl(\omega^{\ast}\bigr) -\sum
_{L\in\mathbf
{B}%
_{+}}D_{\mathbf{B}\setminus\{L\}=\mathbf{1},L=0}\bigl(\omega^{\ast}\bigr)
\nonumber
\\[-8pt]
\\[-8pt]
\nonumber
&&\qquad{} - \sum_{\widetilde{\mathbf{B}}\dvtx \varnothing\neq\widetilde{%
\mathbf{B}}\subseteq\mathbf{B}^{\prime}} D_{\mathbf{B}\setminus
\widetilde{\mathbf{B}}=\mathbf{1},\widetilde{\mathbf{B}}=\mathbf
{0}}\bigl(
\omega^{\ast}\bigr) + \sum_{\mathbf{E}\in\mathfrak{T}}D_{\mathbf{B}\setminus
\mathbf{E}=\mathbf{1},\mathbf{E}=\mathbf{0}}
\bigl(\omega^{\ast}\bigr) > 0,
\end{eqnarray}
then $\mathbf{B}$ is singular for $\mathcal{D}(\mathbf{C},{\Omega})$.
\end{theorem}

\begin{pf} By Theorem~\ref{thmsingularchar}, $\mathbf
{B}$ is
singular for $\mathcal{D}(\mathbf{C},{\Omega})$ if and only if
%
%e15 #&#
%e5.4 ###
\begin{equation}
\label{eqsingular-condition} \mbox{for some }\omega^{\ast}\in{\Omega},\qquad
D_{\mathbf{B}=\mathbf{1}%
}\bigl(\omega^{\ast}\bigr) -\sum
_{\widetilde{\mathbf{B}}\subseteq\mathbf{B}}D_{
\mathbf{B}\setminus\widetilde{\mathbf{B}}=\mathbf{1},\widetilde{\mathbf
{B}}=%
\mathbf{0}}\bigl(\omega^{\ast}\bigr)>0.
\end{equation}
Suppose for a contradiction that (\ref{eqsingular-condition}) does not hold
but (\ref{eqfourterms}) holds for some $\omega^* \in\Omega$. Since $%
\mathbf{B}_+$ has a positive monotonic effect on $D$ relative to
$\mathbf{C}$%
, if $\widetilde{\mathbf{B}}\subseteq\mathbf{B}$ is such that
$\widetilde{%
\mathbf{B}}\cap\mathbf{B}_+\neq\varnothing$, then $D_{\mathbf
{B}\setminus
\widetilde{\mathbf{B}}=\mathbf{1},\widetilde{\mathbf{B}}=\mathbf{0}%
}(\omega^*)=1$ implies $D_{\mathbf{B}\setminus\{L\}=\mathbf{1}%
,L=0}(\omega^*)=1$ for some $L\in\widetilde{\mathbf{B}}$. Hence for
all $%
\omega\in\Omega$,
%
%e16 #&#
%e5.5 ###
\begin{equation}
\label{eqproofstep}\qquad D_{\mathbf{B}=\mathbf{1}}(\omega) -\sum
_{L\in\mathbf{B}}D_{\mathbf{B}%
\setminus\{L\}=\mathbf{1},L=0}(\omega) - \sum
_{ \widetilde{%
\mathbf{B}}\subseteq\mathbf{B}^\prime, |\widetilde{\mathbf{B}}| \geq2 } D_{\mathbf{B}\setminus\widetilde{\mathbf{B}}=\mathbf{1},\widetilde{%
\mathbf{B}}=\mathbf{0}}(\omega)\leq0.
\end{equation}
By applying the same argument to the first two terms on the left-hand
side of (\ref{eqproofstep}) %$D_{\mathbf{B}=\mathbf{1}}(\omega^{\ast
%})-\sum_{L\in\mathbf{B}}D_{\mathbf{B}\setminus\{L\}=\mathbf{1}%
%,L=0}(\omega^{\ast})$
as was applied in the proof of Theorem~\ref{thmmonotone-cond}, we have
that (%
\ref{eqfourterms}) does not hold for all $\omega\in\Omega$, which is a
contradiction.
\end{pf}
%&&D_{\mathbf{B}=\mathbf{1}}(\omega)-\sum_{L\in\mathbf{B}%
%}D_{\mathbf{B}\setminus\{L\}=\mathbf{1},L=0}(\omega)-\sum_{
%,\{L_{1},,\ldots,,L_{k}\}=\mathbf{0}}(\omega) \\
%&&\text{ \ \ \ \ \ }+\sum_{\mathbf{E}\in\mathfrak{T}}D_{\mathbf{B}%
%0
%and this completes the proof.

The following corollary to Theorem~\ref{thmsingular-test} generalizes the
discussion in~\cite{VanderWeele2010a,VanderWeele2010b} to an arbitrary
number of dichotomous factors:

%
%co5.6 #&#
\begin{corollary}
\label{corsingular-test} Let $\mathbf{B}=\mathbf{B}_{+}\,\dot{\cup
}\,\mathbf{B}%
^{\prime}\in \dot{\mathbb{P}}(\mathbb{L}(\mathbf{C}))$, $|\mathbf
{B}|=|%
\mathbf{C}|$. Suppose that every $L\in\mathbf{B}_{+}$ has a positive
monotonic effect on $D$ relative to $\mathbf{B}$, and $\mathbf{W}$ is
sufficient to adjust for confounding of $\mathbf{C}$ on $D$.  If for some
tree $\mathfrak{T}$ on $\mathbf{B}_{+}$, and some $\mathbf{w}$, we have%
%
%e17 #&#
%e5.6 ###
\begin{eqnarray}\label{eqsingular-test}
0 &<&E[D\mid\mathbf{B} = \mathbf{1},\mathbf{W} = \mathbf{w}%
]-\sum
_{L\in\mathbf{B}_{+}}E\bigl[D\mid\mathbf{B}\setminus\{L\} =
\mathbf{1},L=0,\mathbf{W} = \mathbf{w}\bigr]
\nonumber\\
&&{}- \sum_{\widetilde{\mathbf{B}}\dvtx \varnothing\neq\widetilde{\mathbf{B}}%
\subseteq\mathbf{B}^{\prime}}E[D\mid\mathbf{B}\setminus
\widetilde{%
\mathbf{B}}=\mathbf{1},\widetilde{\mathbf{B}}=\mathbf{0},
\mathbf {W}=\mathbf{w%
}]
\\
&&{}+\sum_{\mathbf{E}\in\mathfrak{T}}E[D\mid\mathbf{B}\setminus\mathbf
{E}=%
\mathbf{1},\mathbf{E}=\mathbf{0},\mathbf{W}=\mathbf{w}],
\nonumber
\end{eqnarray}
then $\mathbf{B}$ is singular for $\mathcal{D}(\mathbf{C},{\Omega})$.
\end{corollary}

\begin{pf} By applying Proposition~\ref
{propadjconfound} to
each term in (\ref{eqfourterms}), the proof is complete.
\end{pf}

Condition (\ref{eqsingular-test}) leads directly to a statistical test of
compositional epistasis. This is notable since some claims in the genetics
literature appear to suggest that such tests did not exist~\cite{Cordell2002}.

As stated in the next corollary, from Theorem~\ref{thmsingular-test},
if all
or all but one of the elements of $\mathbf{B}$ have positive monotonic
effects on $D$, then singularity and irreducibility coincide:

%
%co5.7 #&#
\begin{corollary}
\label{corsingularirreducible} Suppose $|\mathbf{B}|=|\mathbf{C}|$ and
that for all or all but one of $B_{i}\in\mathbf{B}$, $B_{i}$ has a positive
monotonic effect on $D$ relative to $\mathbf{B}$, then $\mathbf{B}$ is
singular for $\mathcal{D}(\mathbf{C},{\Omega})$ if and only if $\mathbf
{B}$ is
irreducible for $\mathcal{D}(\mathbf{C},{\Omega})$.
\end{corollary}

An important consequence of this corollary is that when there is
at most one variable that does not have a positive monotonic effect,
condition (\ref{eqcondmonotoneempirical}) establishes that $\mathbf{B}$
is singular in addition to being irreducible for $\mathcal{D}(\mathbf
{C},{%
\Omega})$.

\begin{pf} Let $\mathbf{B}^\prime$ denote the one or zero
elements of $\mathbf{B}$ that do not have a monotonic effect on $D$ relative
to $\mathbf{C}$. If $\mathbf{B}$ is irreducible for $\mathcal{D}(\mathbf
{C},{%
\Omega})$, then by the argument in the proof of Theorem~\ref
{thmmonotone-cond},
\[
D_{\mathbf{B}=\mathbf{1}}\bigl(\omega^{\ast}\bigr)-\sum
_{L\in\mathbf{B}}D_{%
\mathbf{B}\setminus\{L\}=\mathbf{1},L=0}\bigl(\omega^{\ast}\bigr)+\sum
_{%
\mathbf{E}\in\mathfrak{T}}D_{\mathbf{B}\setminus\mathbf{E}=\mathbf
{1},%
\mathbf{E}=\mathbf{0}}\bigl(\omega^{\ast}
\bigr)>0.
\]
Since the third term on the left-hand side of (\ref{eqfourterms})
vanishes when $|%
\mathbf{B}^\prime|\leq1$, it follows that $\mathbf{B}$ is singular for
$%
\mathcal{D}(\mathbf{C},{\Omega})$. The converse is given in Corollary~\ref{corsingularimplies}.
\end{pf}

%
%co5.8 #&#
\begin{corollary}
\label{corsingularirreducible2} Suppose $|\mathbf{B}|=|\mathbf{C}|$ and
that for all or all but one of $B_{i}\in\mathbf{B}\,\dot{\cup}\, \mathbf
{C}^\prime$, $%
B_{i}$ has a positive monotonic effect on $D$ relative to $\mathbf
{B}\cup
\mathbf{C}^\prime$, for all $X^\prime\in\mathbf{C}^\prime$, $X^\prime
$ is
not causally influenced by $\mathbf{C}$ and $\mathbf{B}$ is singular
for $%
\mathcal{D}(\mathbf{C},{\Omega})$, then $\mathbf{B}$ is singular for $%
\mathcal{D}(\mathbf{C}\cup\mathbf{C}^\prime,{\Omega})$.
\end{corollary}

\begin{pf} By Corollary~\ref{corsingularirreducible},
$\mathbf B$ is irreducible relative to $\mathcal{D}(\mathbf{C},{\Omega})$.
Hence by Corollary~\ref{corexpand} $\mathbf B$ is irreducible relative to
$\mathcal{D}(\mathbf{C}\cup\mathbf{C}^\prime,{\Omega})$.
The conclusion then follows from a further application of
Corollary~\ref{corsingularirreducible}.
\end{pf}

%s5.1 #&#
%s5.1 ###
\subsection{Relation to Pearl's probability of causation}
\label{secprobcause}

%
%There are certain relations between the conditions for sufficient cause
%interactions and Pearl's ``probability of causation" that are worth
%noting,
%specifically concerning Pearl's discussion of necessary and sufficient
%causes.
Pearl~\cite{Pearl2000}, Chapter~9, defined the \textit{probability of necessity
and sufficiency {(PNS)} of cause $X$ for outcome $D$} to be $%
P(D_{X=1}(\omega)=1,D_{X=0}(\omega)=0)$. In other words $\operatorname{PNS}(D,X)$ is the
probability that $D$ would occur if $X$ occurred and would not have
done so
had $X$ not occurred. We generalize this to the setting in which there are
multiple causes $\mathbf{B}$:

%
%de5.9 #&#
\begin{defin}
\label{defpns} For $\mathbf{B} \subseteq\dot{\mathbb{P}}({\mathbb{L}}(
\mathbf{C}))$, the \textit{probability of necessity and sufficiency of $
\mathbf{B}$ causing $D$} is
\[
\operatorname{PNS}(D,\mathbf{B}) \equiv P ( D_{\mathbf{B}=\mathbf{1}%
} =1 \mbox{ and for all }
\mathbf{b}\neq\mathbf{1}, D_{\mathbf{B}=%
\mathbf{b}} =0 ).
\]
\end{defin}

Thus $\operatorname{PNS}(D,\mathbf{B})$ is the probability that $D$ would occur
if every literal $L\in\mathbf{B}$ occurred and would not have done so had
at least one literal in $\mathbf{B}$ not occurred.

%
%pr5.10 #&#
\begin{proposition}
If $|\mathbf{B}|= |\mathbf{C}|$, then $\operatorname{PNS}(D,\mathbf{B})>0$ if and only
if $%
\mathbf{B}$ is singular for $\mathcal{D}(\mathbf{C},{\Omega})$.
\end{proposition}

\begin{pf}
The proof follows directly from Theorem~\ref{thmsingularchar} and Definition~\ref{defpns}.
\end{pf}

This connection also provides an interpretation for condition (\ref{eqsingular-test}).
For expositional convenience in the following
proposition, we assume that $\mathbf{B}$ and $D$ are unconfounded and
do not make monotonicity assumptions;
it would be straightforward to do so.

%
%pr5.11 #&#
\begin{proposition}
\label{proppns-connection} Under the conditions of Corollary~\ref{corsingular-test},
with $\mathbf{W}=\varnothing=\mathbf{B}_+$, $\operatorname{PNS}(%
\mathbf{B},D)$ is bounded below by the right-hand side of (\ref{eqsingular-test}%
).
\end{proposition}

\begin{pf} %Thus, the conditions for
%detecting a sufficient cause interaction in Theorem~\ref{thmcharsci}(i) are satisfied if
%and only if the probability of necessity and sufficiency of
%$X_{1},\ldots,X_{m}$
%causing $D$ is greater than zero. Furthermore, the expression in
%Theorem~\ref{thmsciconfound}
%provides a lower bound for $PNS(D,X_{1},\ldots,X_{m})=0$ since
%
\begin{eqnarray*}
\operatorname{PNS}(D,\mathbf{B}) &=& P( D_{\mathbf{B}=\mathbf{1}} =1 %
\mbox{ and for all }\mathbf{b}\neq\mathbf{1}, D_{\mathbf
{B}=\mathbf{b}%
} =0)
\\
&\geq& P( D_{\mathbf{B}=\mathbf{1}} =1) + P( %
\mbox{for all }\mathbf{b}\neq
\mathbf{1}, D_{\mathbf{B}=\mathbf{b}} =0) - 1
\\
&=& P( D_{\mathbf{B}=\mathbf{1}} =1) - P( %
\mbox{for some }\mathbf{b}
\neq\mathbf{1}, D_{\mathbf{B}=\mathbf{b}} =1)
\\
&\geq& P( D_{\mathbf{B}=\mathbf{1}} =1) - \sum_{\mathbf{b}%
\neq\mathbf{1}} P(
D_{\mathbf{B}=\mathbf{b}} =1)
\\
&=& E [ D = 1 \mid{\mathbf{B}=\mathbf{1}} ] - \sum_{\mathbf
{b}\neq
\mathbf{1}}
E [ D \mid{\mathbf{B}=\mathbf{b}} ]
\end{eqnarray*}
which is the right-hand side of (\ref{eqsingular-test}) with $\mathbf
{W}%
=\varnothing= \mathbf{B}_+$.
\end{pf}

This generalizes some of the lower bounds on $\operatorname{PNS}(D,X)$ given
by Pearl~\cite{Pearl2000}, Section~9.2.

%s6 #&#
%s6 ###
\section{Relation to statistical models with linear links}
\label{seclinearlink}

In related work~\cite{VanderWeele2008a} it is noted that the presence of
interaction terms in statistical models do not, in general, correspond to
sufficient conditions for irreducibility. Consider, for example, a saturated
Bernoulli regression model for $D$ with identity link and binary
regressors $%
\mathbf{C}=\{X_{1},\ldots,X_{p}\}$,
%
%e18 #&#
%e6.1 ###
\begin{equation}
E[D\mid\mathbf{C} = \mathbf{c}]=\sum_{\widetilde{\mathbf{B}}\subseteq
\mathbf{C}}
\beta_{\widetilde{\mathbf{B}}}\Bigl(\bigwedge  (\widetilde{\mathbf{B}})
\Bigr)_{%
\mathbf{c}}. \label{eqbernoulli}
\end{equation}
Note that with $\mathbf{c}=(x_{1},\ldots,x_{p})$, then $(\bigwedge  (\widetilde{%
\mathbf{B}}))_{\mathbf{c}}=\prod_{X_{i}\in\widetilde{\mathbf{B}}}x_{i}$,
the usual product interaction term. The conditions, given earlier,
for detecting the presence of irreducibility and singularity lead to
linear restrictions on the regression coefficients~$\beta_{\widetilde{%
\mathbf{B}}}$.
%
%e19 #&#
%e6.2 ###
\begin{equation}
\sum_{\widetilde{\mathbf{B}}\subseteq\mathbf{C}}m_{\widetilde{\mathbf
{B}}%
}\beta_{\widetilde{\mathbf{B}}}>0.
\label{eqlinearconstraint}
\end{equation}
Note that (\ref{eqlinearconstraint}) includes an intercept $\beta_{\varnothing}$. First we define
\[
\deg_{\mathfrak{T}}(L)\equiv\bigl|\{\mathbf{E}\mid\mathbf{E}\in{\mathfrak
{T}}%
,L\in\mathbf{E}\}\bigr|,
\]
the \textit{degree} of $L$ in a tree $\mathfrak{T}$, to be the number
of edges
in $\mathfrak{T}$ that contain $L$.

%
%pr6.1 #&#
\begin{proposition}
\label{proplinearconstraintirred} Under the conditions of Theorem~\ref{thmmonotone-cond}, with $\mathbf{B} = \mathbf{C}$, condition (%
\ref{eqcondmonotone}) is equivalent to restriction (\ref{eqlinearconstraint}) with $m_{\widetilde{\mathbf{B}}}=
m_{\widetilde{\mathbf{B}}}^{\mathrm{irred}}$ where
%
%e20 #&#
%e6.3 ###
\begin{eqnarray}
\label{eqmdef} m_{\widetilde{\mathbf{B}}}^{\mathrm{irred}}& \equiv& 1 - |\mathbf{B} \setminus
\widetilde{\mathbf{B}}| + |\mathfrak{T}|
\nonumber
\\[-8pt]
\\[-8pt]
\nonumber
&&{}- \sum_{L \in
\widetilde{\mathbf{B}} \cap\mathbf{B}_+}
\deg_{\mathfrak{T}}(L) +\bigl | \{ \mathbf{E} \mid\mathbf{E} \in{\mathfrak{T}},
\mathbf{E}\subseteq \widetilde{\mathbf{B}} \cap\mathbf{B}_+\}\bigr|.
\end{eqnarray}
\end{proposition}

Note that since $\mathfrak{T}$ is a tree on $\mathbf{B}_+$, the last
term in
(\ref{eqmdef}) has a natural graphical interpretation\vspace*{1pt} as the number of
edges in the ``induced subgraph'' of $\mathfrak{T}$ on the subset~$\widetilde{%
\mathbf{B}}$. Definition (\ref{eqmdef}) also subsumes condition (\ref{eqempcond})
given in Theorem~\ref{thmsciconfound} (without
monotonicity), in which case the last three terms in (\ref{eqmdef}) vanish.
Though Proposition~\ref{proplinearconstraintirred} assumes that
$\mathbf{C}%
_2=\varnothing$, the condition given by (\ref{eqlinearconstraint}) and
(\ref{eqmdef}) continues to apply in the case where $\mathbf{c}_2=\mathbf
{0}$, as
in Corollaries~\ref{cormsctest} and~\ref{cormsc-monotone-cond-empirical};
obvious extensions apply to the case where $\mathbf{c}_2\neq\mathbf
{0}$.\looseness=1

\begin{pf} This follows by counting the number (and sign) of
expectations in~(\ref{eqcondmonotone}) for which $\widetilde{\mathbf{B}}$
is a subset of the variables assigned the value $1$ in the conditioning
event. The first two terms in (\ref{eqcondmonotone}) correspond,
respectively, to the first two terms in (\ref{eqmdef}). The remaining three
terms in (\ref{eqmdef}) correspond to the last sum in~(\ref
{eqcondmonotone}): $|\mathfrak{T}%
| $, the number of edges in $\mathfrak{T}$, is the total number of
terms in
the sum. The sum over degrees subtracts the number of terms in which
some\vspace*{1pt} $L
\in\widetilde{\mathbf{B}}$ is assigned zero. Since this double counts terms
corresponding to edges contained in $\widetilde{\mathbf{B}}$, the last term
corrects for this.
\end{pf}

%
%pr6.2 #&#
\begin{proposition}
\label{proplinearconstraintsingular} Under the conditions of
Theorem~\ref{thmsingular-test}, with $\mathbf{B} = \mathbf{C}$, condition
(\ref{eqfourterms}) is equivalent to restriction (\ref{eqlinearconstraint}) with $m_{\widetilde{\mathbf{B}}}=
m_{\widetilde{\mathbf{B}}}^{\mathrm{sing}}$ where
%
%e21 #&#
%e6.4 ###
\begin{equation}
\label{eqmdefsingular} m_{\widetilde{\mathbf{B}}}^{\mathrm{sing}} \equiv
m_{\widetilde
{%
\mathbf{B}}}^{\mathrm{irred}} + \bigl(\bigl|\mathbf{B}^\prime\setminus
\widetilde{\mathbf{B}}\bigr|\bigr) - \bigl(2^{|\mathbf{B}^\prime\setminus\widetilde{%
\mathbf{B}}|} - 1\bigr).
\end{equation}
\end{proposition}

\begin{pf} Expression (\ref{eqmdefsingular}) follows from
another counting argument similar to the proof of Proposition~\ref{proplinearconstraintirred},
together with the observation that
conditions (\ref{eqcondmonotone}) and (\ref{eqfourterms}) only differ in that the
terms in the sum over $L$ in (\ref{eqcondmonotone}) for $L\in\mathbf
{B}%
^\prime$ are replaced by a sum over all subsets of $\mathbf{B}^\prime
$.
\end{pf}

%
%ex4 #&#
\begin{example}[(Two-way interactions)]
Consider the saturated Bernoulli regression with identity link
with $\mathbf{C}=\{X_{1},X_{2}\}$.
\[
E[D\mid X_{1}=x_{1},X_{2}=x_{2}]=
\beta_{\varnothing}+\beta_{1}x_{1}+\beta_{2}x_{2}+
\beta_{12}x_{1}x_{2}.
\]
Suppose that $X_{1}$ and $X_{2}$ are unconfounded with respect to $D$,
so (\ref{eqnoconfound}) holds with $\mathbf{W}=\varnothing$. Proposition~\ref{proplinearconstraintirred}
implies that $\{X_{1},X_{2}\}$ is irreducible
relative to $\mathbf{C}$ if $\beta_{12}>\beta_{\varnothing}$;
Proposition~\ref{proplinearconstraintsingular} implies that $\{X_{1},X_{2}\}$ is
singular relative to $\mathbf{C}$ if $\beta_{12}>2\beta_{\varnothing
}$. If
one of $X_{1}$ or $X_{2}$ have positive monotonic effects on $D$
relative to
$\mathbf{C}$, then Proposition~\ref{proplinearconstraintirred} and
Corollary~\ref{corsingularirreducible} imply that $\{X_{1},X_{2}\}$ is both
irreducible and
singular relative to $\mathbf{C}$ if $\beta_{12}>\beta_{\varnothing}$.
If $%
X_{1}$ and $X_{2}$ have positive monotonic effects on $D$ relative to $%
\mathbf{C}$, then Proposition~\ref{proplinearconstraintirred} and
Corollary~\ref{corsingularirreducible} imply that $\{X_{1},X_{2}\}$ is both
irreducible and
singular relative to $\mathbf{C}$ if $\beta_{12}>0$.

Thus only under the assumption of positive monotonic effects for both $%
X_{1}$ and $X_{2}$ does the sufficient condition for the irreducibility and
singularity of $\{X_{1},X_{2}\}$ coincide with the classical two-way
interaction term $\beta_{12}$ being positive. Note that under the
assumption of negative monotonic effects of
$X_{1}$ and $X_{2}$ on $D$, $\beta_{12} < 0$ is equivalent to
irreducibility and singularity for $\bar{D} \equiv(1-D)$;
see~\cite{vanderweeleremarks2011} for this and other remarks on
recoding of
exposures or outcomes.

It also follows from Proposition~\ref{propmscandirred} that if $%
\{X_{1},X_{2}\}$ is irreducible relative to~$\mathbf{C}$, then there exists
some $\omega\in\Omega$ for whom $\{X_{1},X_{2}\}$ is a minimal sufficient
cause relative to $\mathbf{C}$ (since $|\mathbf{B}|=|\mathbf{C}|$).
\end{example}

%
%ex5 #&#
\begin{example}[(Three-way interactions)]\label{exthree-way}
The saturated Bernoulli regression with three binary variables and a
identity link can be written as%
\begin{eqnarray*}
&&E[D=1\mid X_{1}=x_{1},X_{2}=x_{2},X_{3}=x_{3}]
\\
&&\qquad=\beta_{\varnothing}+\beta_{1}x_{1}+
\beta_{2}x_{2}+\beta_{3}x_{3} +
\beta_{12}x_{1}x_{2}+\beta_{13}x_{1}x_{3}\\
&&\qquad\quad{}+
\beta_{23}x_{2}x_{3}+\beta_{123}x_{1}x_{2}x_{2}.
\end{eqnarray*}
Suppose that $\mathbf{C} = \{ X_{1}, X_{2}, X_{3}\}$ is unconfounded
for $D$. Proposition~\ref{proplinearconstraintirred} implies that $\{
X_1,X_2,X_3\}$
is irreducible relative to $\mathbf{C}$ if
%
%e22 #&#
%e6.5 ###
\begin{equation}
\label{eqfirst3waybeta} \beta_{123}> 2\beta_{\varnothing} +
\beta_{1} + \beta_{2} + \beta_{3}.\vadjust{\goodbreak}
\end{equation}
It follows from Proposition~\ref{propmscandirred} that if $%
\{X_1,X_2,X_3\} $ is irreducible relative to $\mathbf{C}$, then there exists
some $\omega\in\Omega$ for whom $\{X_1,X_2,X_3\}$ is a minimal sufficient
cause relative to $\mathbf{C}$ (since $|\mathbf{B}|=|\mathbf{C}|$).

Proposition~\ref{proplinearconstraintsingular} implies $%
\{X_1,X_2,X_3\}$ is singular relative to $\mathbf{C}$ if
\[
\beta_{123} > 6 \beta_{\varnothing} + 2 \beta_{1} + 2
\beta_{2} + 2\beta_{3}.
\]

However, if $X_{1}$, $X_{2}$ and $X_{3}$ have positive monotonic
effects on $D$ (relative to~$\mathbf{C}$), then Proposition~\ref{proplinearconstraintirred}
implies $\{X_1,X_2,X_3\}$ is irreducible
relative to $\mathbf{C}$ if any of the following hold:
%
%e23 #&#
%e6.6 ###
\begin{equation}
\label{eqsecond3waybeta} \beta_{123} > \beta_{1},\qquad
\beta_{123} > \beta_{2},\qquad \beta_{123} >
\beta_{3};
\end{equation}
equivalently, $\beta_{123} > \min\{\beta_1,\beta_2,\beta_3\}$. By
Corollary~\ref{corsingularirreducible} this also establishes that $\{
X_1,X_2,X_3\}$
is singular relative to $\mathbf{C}$.

If only $X_{1}$ and $X_{2}$ have positive monotonic effects on $D$
relative to $\mathbf{C}$, then Proposition~\ref{proplinearconstraintirred}
implies that $\{X_1,X_2,X_3\}$ is irreducible relative to $\mathbf{C}$ if
%
%e24 #&#
%e6.7 ###
\begin{equation}
\label{eqthird3waybeta} \beta_{123} > \beta_{\varnothing}+
\beta_{1}+\beta_{2}.
\end{equation}
By Corollary~\ref{corsingularirreducible}, condition (\ref
{eqthird3waybeta}%
) also implies that $\{X_1,X_2,X_3\}$ is singular relative to $\mathbf{C}$
(since only $X_3$ does not have a positive monotonic effect on $D$). As we
would expect condition (\ref{eqthird3waybeta}) is weaker than (\ref{eqfirst3waybeta}),
but stronger than any of the conditions in (\ref{eqsecond3waybeta}). If only one potential cause has a monotonic
effect on $%
D$ relative to $\mathbf{C}$, then we can only use (\ref{eqfirst3waybeta}) to
establish irreducibility.

Thus for three-way interactions, $\beta_{123}>0$ does not correspond
to any of the sufficient conditions for irreducibility or singularity
of $%
\{X_1,X_2,X_3\}$ relative to $\mathbf{C}$, regardless of whether or not
monotonicity assumptions hold.
\end{example}

%A saturated Bernoulli regression
%with linear link with $s$ binary variables can be written as
%pr(D &=&1|X_{1}=x_{1},,\ldots,,X_{s}=x_{s})=\beta_{0}+\dsum_{i}\beta
%_{i}x_{i}+\dsum_{i_{1}\neq i_{2}}\beta
%_{i_{1}i_{2}}x_{i_{1}}x_{i_{2}} \\
%&&\text{ }+,\ldots,+\dsum_{i_{1},i_{2},,\ldots,,i_{s-1}}\beta
%_{i_{1},i_{2},,\ldots,,i_{s-1}}x_{i_{1}},\ldots,x_{i_{s-1}}+\beta
%_{12,\ldots,s}x_{1}x_{2},\ldots,x_{s}
%where the $k$th sum is over all combinations of $k$ distinct indices.
%can be verified by induction that Theorem~\ref{thmsciconfound}
%implies that an $s$-way
%sufficient cause interaction is present if%
%_{i}+(s-3)\dsum_{i_{1}\neq i_{2}}\beta_{i_{1}i_{2}}+,\ldots, \\
%&&+2\dsum_{i_{1},i_{2},,\ldots,,i_{s-3}}\beta
%_{i_{1},i_{2},,\ldots,,i_{s-3}}+\dsum_{i_{1},i_{2},,\ldots,,i_{s-2}}\beta
%_{i_{1},i_{2},,\ldots,,i_{s-2}}.
%There are no simple expressions in terms of regression coefficients
%for the
%empirical conditions implied by Theorems~\ref{thmsubordtest1}-
%cause interactions under the various monotonicity assumption; as can
%be seen
%from the list of subordinate sets in section~\ref{sectestsnway},
%there will be many different
%possible conditions on the coefficients that suffice to conclude the
%presence of a sufficient cause interaction under monotonicity.

%s7 #&#
%s7 ###
\section{Interpretation of sufficient cause models}
\label{secinterpretation}

%%% Use this command line to make the figure
%%% m4 sop.m4 | dpic -p > sop.tex
% \begin{tabular}[c]{p{5cm}p{5cm}}
%&

As mentioned in Section~\ref{secpotoutcomes} the observed data for an
individual
$(\mathbf{C}(\omega),D(\omega))$ represents a strict subset of the
potential outcomes ${\mathcal{D}}(\mathbf{C},\omega)$; this is the
``fundamental problem of causal inference.'' Further, as we have seen,
for a
given set of potential outcomes there can exist different determinative sets
of minimal sufficient causes ${\mathfrak{B}}$ for the same set of potential
outcomes; see (\ref{eqomega2a}) and (\ref{eqomega2b}). Thus we have the
following for an individual:
%
%e25 #&#
%e7.1 ###
\begin{equation}
\begin{array} {ccccc} \vdots& \searrow& \vdots& \searrow&
\\[3pt]
{\mathfrak{B}} & \rightarrow& {\mathcal{D}}(\mathbf{C},\omega) & \rightarrow&
\bigl(\mathbf{C}(\omega),D(\omega)\bigr).
\\
\vdots& \nearrow& \vdots& \nearrow&
\\
& \mbox{\scriptsize many-one} & & \mbox{\scriptsize many-one} & \end{array} %
\label{eqmany-one}
\end{equation}
%
%This double many-one relationship is also present at the population
%level:
% \begin{array}{ccccc}
%(\mathbf{A},{\mathfrak B}) & \rightarrow& {\mathcal D}(\mathbf{C},
%& \mbox{\tiny many-one} && \mbox{\tiny many-one}
It is typically impossible to know the set of potential outcomes for an
individual ${\mathcal{D}}(\mathbf{C},\omega)$, even when $\mathbf{C}=\{
X\}$%
, even from randomized experiments. However, possession of this knowledge
would permit one to predict how a given individual would respond under
any given pattern of exposures $\mathbf{C}=\mathbf{c}$.\vadjust{\goodbreak}

The results in this paper demonstrate that, given data from a randomized
experiment (or when sufficient variables are measured to adjust for
confounding), it is possible to infer the existence of an individual
for whom
all sets of minimal sufficient causes $\mathfrak{B}$ have certain features
in common. However, given the double many-one relationship (\ref
{eqmany-one}%
), and the fact that the set of potential outcomes ${\mathcal
{D}}(\mathbf{C}%
, \omega)$, if they were known, apparently address all potential
counterfactual queries, it is natural to ask what is to be gained by
considering such representations. We now motivate our results by presenting
several different interpretations of sufficient cause representations.

%&{\mathcal D}(\mathbf{C}, \omega) &\xleftrightarrow{\mbox{\tiny
%many-one}}
%(\mathbf{A},{\mathfrak B}) \xleftrightarrow{\mbox{\tiny many-one}}& {
%many-one}}

%s7.1 #&#
%s7.1 ###
\subsection{The descriptive interpretation}
Under this view, sets of minimal sufficient causes are merely a way to
describe the set of potential outcomes ${\mathcal{D}}(\mathbf{C},\Omega)$.
The representation $(\mathbf{A},{\mathfrak{B}})$ may be more compact;
compare Table~\ref{taboutcomes} and (\ref{eq21}). Extending this to
a population $\Omega$, the variables $%
\mathbf{A}$ in a representation $(\mathbf{A},{\mathfrak{B}})$ merely
describe subpopulations with particular patterns of potential outcomes.

Knowing that there
exists an individual for whom all representations ${\mathfrak{B}}$ have
certain features in common provides qualitative information about the
set of
potential outcomes.

For two binary causes
$\{X_{1},X_{2}\}$, Theorem~\ref{thmcharsci}
implies that $\{X_{1},X_{2}\}$ is irreducible relative
to $\mathbf{C}$ for $\omega^{\ast}$ if $D_{11}(\omega^{\ast})=1$
and $%
D_{10}(\omega^{\ast})=D_{01}(\omega^{\ast})=0$. Such a pattern is
of interest insofar as it indicates
that the causal process resulting in this individual's potential outcomes
${\mathcal{D}}(\mathbf{C},\omega^{\ast} )$ is such that (for some
setting of the variables in $\mathbf{C}\setminus\{X_1, X_2\}$), $D=1$
if both $X_1=1$ and $X_2=1$, but not when $X_1=1$ and $X_2=0$ or vice versa.

Similarly it follows from Theorem~\ref{thmsingularalt} that if
$\{X_{1},X_{2}\}$ is
singular relative to $\mathbf{C}$ for $\omega^{\ast}$ then
$D_{11}(\omega^{\ast})=1$ and $D_{10}(\omega^{\ast})=D_{01}(\omega^{\ast
})=D_{00}(\omega^{\ast})=0$. Hence
the causal process producing ${\mathcal{D}}(\mathbf{C},\omega^{\ast} )$
is such that, for some setting of the variables in $\mathbf{C}\setminus
\{X_1, X_2\}$, $D=1$ if both $X_1=1$ and $X_2=1$, but not when either
$X_1=0$ or $X_2=0$.

In contrast to the classical notions of interaction arising in linear
models (see Section~\ref{seclinearlink}),
irreducibility and singularity are causal in that they relate to the
potential outcomes.
Sections~\ref{sectestsnway} and~\ref{secsingular} contain
empirical tests
for the presence of irreducible or singular interactions.

%Such a pattern of potential outcomes
%These patterns of potential outcomes may be
%of interest because they represent more causal notions of interaction
%than
%the notion that arises from the use of statistical models as in .
% We have derived empirical tests for these causal
%interaction notions of irreducibilty and singularity in

%s7.2 #&#
%s7.2 ###
\subsection{Generative mechanism interpretations}
A minimal sufficient cause representation may be interpreted in terms
of a ``generative mechanism'':

%
%de7.1 #&#
\begin{defin}\label{defmech}
A \textit{ mechanism $M(\omega)$ relative to $\mathbf C$} takes as input
an assignment
$\mathbf c$ to $\mathbf C$,
and outputs a ``state'' $M_\mathbf{c}(\omega)$ which is either
``active'' ($1$) or ``inactive''~($0$). A mechanism is said to be
\textit{generative} for $D$ if whenever it is active, the event $D=1$
is caused, so that $M_\mathbf{c}(\omega)=1$ implies $D_\mathbf{c}(\omega
) = 1$. Conversely, a mechanism is said to be \textit{preventive} for
$D$ if whenever $M_\mathbf{c}(\omega)=1$, $D_\mathbf{c}(\omega) = 0$ is caused.
\end{defin}

Though this definition refers to a mechanism ``causing'' $D=1$ or
$D=0$, we
abstain from defining this more formally in terms of potential outcomes
until the next section. Our reason for proceeding in this way is that
there may be circumstances in which an investigator is able to posit the
existence of a causal mechanism causing $D=1$ or $D=0$, for example,
based on
experiments manipulating the inputs $\mathbf C$ and output $D$, but lacks
sufficiently detailed information to posit well-defined counterfactual
outcomes involving interventions on these (hypothesized) mechanisms.

%
%de7.2 #&#
\begin{defin}\label{defexhaust}
A set of generative mechanisms $\mathbf{M}=\langle
M^{1},\ldots, M^{p}\rangle$ will be said to be \textit{exhaustive} for
a given
set of potential outcomes ${\mathcal{D}}(%
\mathbf{C},\Omega)$ if for all $\omega\in\Omega$, and all $\mathbf c$,
if $D_\mathbf{c}(\omega) = 1$, then for some $M^i \in\mathbf{M}$,
$M^i_{\mathbf{c}}(\omega) = 1$.
\end{defin}

Note that if $\mathbf M$ forms an exhaustive set of mechanisms for
${\mathcal{D}}(%
\mathbf{C},\Omega)$,
then it follows that in a context in which no mechanism $M^i$ is
active, $D=0$.

%
%pr7.3 #&#
\begin{proposition} If $\mathbf{M}=\langle
M^{1},\ldots, M^{p}\rangle$
forms an exhaustive set of generative mechanisms
for ${\mathcal{D}}(%
\mathbf{C},\Omega)$, then
$D=\bigvee  (\mathbf{M})$ and
$
D_{\mathbf{c}}(\omega) =\bigvee  (\mathbf{M}_{\mathbf{c}%
}(\omega)).
$
\end{proposition}

\begin{pf} The proof follows from Definitions
\ref{defmech} and~\ref{defexhaust}.
\end{pf}

%
%pr7.4 #&#
\begin{proposition}\label{propmechanism-irreducible}
Suppose $\mathbf{M}$ forms an exhaustive set of generative mechanisms
for ${%
\mathcal{D}}(\mathbf{C},\Omega)$. If $\mathbf{B}\in\dot{\mathbb{P}}(%
\mathbb{L}(\mathbf{C}))$, $|\mathbf{B}|=|\mathbf{C}|$ and $\mathbf B$ is
irreducible for ${%
\mathcal{D}}(\mathbf{C};\Omega)$, then there exists an individual
$\omega^{\ast}$ and a mechanism $M^{i}$ such that
$M^i_{\mathbf{B}=\mathbf{1}}(\omega^{\ast})=1$
but for all $L\in{\mathbf{B}}$, $M^i_{\mathbf{B}\setminus\{L\}=\mathbf
{1},L=0}(\omega^{\ast})=0$.
\end{proposition}

Thus if there exists an exhaustive set of generative mechanisms for
${%
\mathcal{D}}(\mathbf{C},\Omega)$ and
$\mathbf B$ is irreducible, then there is an individual $\omega^*$ and a
mechanism $M^i$ such that $M^i$ is active when
all the literals in $\mathbf B$ take the value $1$, and is inactive when
any one
literal is $0$, and the rest continue to take the value $1$.

\begin{pf} By Theorem~\ref{thmcharsci}, since $\mathbf
{B}$ is irreducible for ${\mathcal{D}}(%
\mathbf{C};\Omega)$, there exists $\omega^{\ast}\in\Omega$ such
that $%
D_{\mathbf{B}=\mathbf{1}}(\omega^{\ast})=1$ and for all $L\in\mathbf{B}$,
$D_{\mathbf{B}\setminus\{L\}=\mathbf{1},L=0}(\omega^{\ast})$. Since
$%
\mathbf{M}$ is an exhaustive set of generative
mechanisms for ${\mathcal{D}}(\mathbf{C}%
,\Omega)$, we have that for all $\mathbf{c}$, $D_{\mathbf{c}}(\omega^{\ast
})=\bigvee  (\mathbf{M}_{\mathbf{c}}(\omega)^{\ast})$. Since
$D_{\mathbf{B}=%
\mathbf{1}}(\omega^{\ast})=1$, for some $M_{i} \in\mathbf{M}$,
$M^i_{\mathbf{B}=\mathbf{1}}(\omega^{\ast})=1$. Since for
all $L\in\mathbf{B}$, $D_{\mathbf{B}\setminus\{L\}=\mathbf
{1},L=0}(\omega^{\ast})=0$
we have that $M^i_{\mathbf{B}\setminus\{L\}=\mathbf{1},L=0}(\omega^{\ast})=0$.
\end{pf}

%
%pr7.5 #&#
\begin{proposition}\label{propmechanism-singular} Suppose $\mathbf{M}$
forms an exhaustive set of
generative mechanisms for ${%
\mathcal{D}}(\mathbf{C},\Omega)$. If $\mathbf{B}\in\dot{\mathbb{P}}(%
\mathbb{L}(\mathbf{C}))$, $|\mathbf{B}|=|\mathbf{C}|$ and $\mathbf B$ is
singular for ${%
\mathcal{D}}(\mathbf{C};\Omega)$, then there exists an individual
$\omega^{\ast}$ and a mechanism $M^{i}$ such that $M^{i}_{\mathbf{B}=\mathbf
{b}}(\omega^{\ast})=1$
if and only if $\mathbf{b}=\mathbf{1}$.
\end{proposition}

Hence under the conditions of Proposition~\ref
{propmechanism-singular}, if
$\mathbf B$ is singular, then there is an individual $\omega^*$ and a
mechanism $M^i$ such that $M^i$ is active if and only if
all the literals in $\mathbf B$ take the value $1$.

\begin{pf} The proof is similar to the proof of
Proposition~\ref{propmechanism-irreducible},
replacing Theorem~\ref{thmcharsci} by Theorem~\ref{thmsingularalt}.
\end{pf}

As the next example shows, the assumption that there
exists an exhaustive set of generative mechanisms
is substantive, and does not hold in all cases.

%
%ex6 #&#
\begin{example}\label{exprevent}
Suppose $\mathbf{C}=\{X_{1},X_{2}\}$ where $X_{1}$ and $X_{2}$
denote the presence of a variant allele at two particular loci. Let
$M^1$ and $M^2$ denote two different proteins such that $M^i$ is
produced if and only if $X_{i}=0$, that is, the associated allele is
not present. Finally, let $D$
denote some characteristic whose occurrence is blocked by the presence
of either $M^1$ or $M^2$ (or both). In this example,
\begin{eqnarray*}
M^i_{x_{1}x_{2}}(\omega) &=&(1-x_{i}),
\\
D_{x_{1}x_{2}}(\omega) &=&\bigl(1-M^1_{x_{1}x_{2}}\bigr)\vee
\bigl(1-M^2_{x_{1}x_{2}}\bigr) = x_{1}x_{2}.
\end{eqnarray*}
By De Morgan's law, the second equation here may also be expressed as
\[
1-D_{x_1x_2}(\omega) = M^1_{x_1x_2}(
\omega)M^2_{x_1x_2}(\omega) = 1-x_1x_2.
\]
The mechanisms $M^1$ and $M^2$ are preventive for $D$, so that $D=1$
only occurs when both mechanisms are inactive.
An exhaustive set of generative mechanisms does not exist because
in this example there are no generative mechanisms (all
mechanisms are preventive).
\end{example}

%

%
%causes
%$\bf C$ is a
%Suppose now that under the weak interpretation, one further commits to
%an
%ontological position that `Whenever some event $D=1$ occurs then there
%must
%be some mechanism that brings about $D=1$' (whatever may be meant by
%`mechanism'; we give one possible formalization of mechanism below; cf.
%Machamer et al., 2000). \ We would then have that $\{X_{1},X_{2}\}$
%being
%irreducible relative to $\mathbf{C}$ would imply that `There is some
%mechanism for $D=1$ that operates when both $X_{1}=1$ and $X_{2}=1$
%but that
%does not operate when just one of $X_{1}=1$ or $X_{2}=1$' (since
%irreducibility implies $D_{11}(\omega^{\ast})=1$ and $D_{10}(\omega
%^{\ast
%})=D_{01}(\omega^{\ast})=0$ for some $\omega^{\ast}$). \ We would
%also
%have that $\{X_{1},X_{2}\}$ being singular relative to $\mathbf{C}$
%would
%imply that `There is some mechanism for $D=1$ that operates iff both
%$X_{1}=1
%$ and $X_{2}=1$' (since singularity implies $D_{11}(\omega^{\ast
%})=1$ and $%
%D_{10}(\omega^{\ast})=D_{01}(\omega^{\ast})=D_{00}(\omega^{\ast
%})=0$
%for some $\omega^{\ast}$). \ An attempt to formalize the notion of
%`mechanism' motivates the following strong interpretation.

It is natural to suppose that mechanisms are ``modular'' and thus
may be isolated or rendered inactive without affecting other
such mechanisms. This may be formalized via potential outcomes:

%
%de7.6 #&#
\begin{defin} An exhaustive set of generative mechanisms $\mathbf{M}$
for ${\mathcal{D}}(%
\mathbf{C},\Omega)$ are said to \textit{support counterfactuals} if
there exist well-defined potential outcomes $D_{\mathbf{C}=%
\mathbf{c},\mathbf{M}=\mathbf{m}}(\omega)$ and $D_{\mathbf{M}=\mathbf
{m}%
}(\omega)$ such that
\[
D_{\mathbf{C}=\mathbf{c},\mathbf{M}=\mathbf{m}}(\omega)=D_{\mathbf{M}=%
\mathbf{m}}(\omega) =\Bigl(\bigvee  (\mathbf{M})\Bigr)_\mathbf{m}.
\]
\end{defin}

The important assumption here is the existence of the potential outcomes
$D_{
\mathbf{m}}(\omega)$ and $D_{\mathbf{c},\mathbf{m}}(\omega)$.
Note that if $\mathbf{M}$ supports counterfactuals then
interventions on $\mathbf{C}$
do not affect $D$ if interventions are also made on $%
\mathbf{M}$.

%
%pr7.7 #&#
\begin{proposition}
If the exhaustive set of generative mechanisms $\mathbf{M}$ support
counterfactuals, then
\[
D_{
\mathbf{M}=\mathbf{M}(\omega)}(\omega) = D_{\mathbf{C}=\mathbf
{C}(\omega),\mathbf{M}=\mathbf{M}(\omega
)}(\omega)=D(\omega)
\]
so that
consistency holds for the potential outcomes $D_{
\mathbf{m}}(\omega)$ and $D_{\mathbf{c},\mathbf{m}}(\omega)$.\vadjust{\goodbreak}
\end{proposition}

\begin{pf} This follows because
\[
D_{\mathbf{C}=\mathbf{C}(\omega),\mathbf{M}=%
\mathbf{M}(\omega)}(\omega)=D_{\mathbf{M}=\mathbf{M}(\omega)}(\omega )=\bigvee
\bigl(\mathbf{M}(\omega)\bigr)=D(\omega).
\]
\upqed\end{pf}

%s7.3 #&#
%s7.3 ###
\subsection{Counterfactual interpretation of a
sufficient cause representation}

If we have an exhaustive set of generative mechanisms which supports
counterfactuals, and further
each mechanism is a conjunction of literals, then there will be a
sufficient cause representation that itself supports counterfactuals.

%
%de7.8 #&#
\begin{defin}\label{defsupport-counterfact-msc}
A representation $(\mathbf{A},{\mathfrak{B}})$
for ${\mathcal{D}}(\mathbf{C},\Omega)$ will be said to \textit{be
structural} if for each
pair $(A_{i},\mathbf{B}_{i})$, $A_{i}\mathbf{\in A}$, $\mathbf
{B}_{i}\mathbf{%
\in{\mathfrak{B}}}$ there exists a generative mechanism (or
mechanisms) $M^i$ such that $M^{i}=A_{i}\wedge(\bigwedge
(\mathbf{B}_{i}))$ and
\[
M^{i}{}_{\mathbf{C}=\mathbf{c}}(\omega)=A_{i}(\omega)\wedge\Bigl(
\bigwedge  (\mathbf{B}%
_{i})
\Bigr){}_{\mathbf{c}}.
\]
\end{defin}
Thus if
$(\mathbf{A},{\mathfrak{B}})$ is structural, then
each pair $(A_{i},\mathbf{B}_{i})$, $A_{i}\mathbf{\in A}$, $%
\mathbf{B}_{i}\mathbf{\in{\mathfrak{B}}}$ corresponds to a mechanism
$M^{i}$.
Thus in this case the variables $A_i(\omega)$ may be interpreted as
indicating whether the corresponding
mechanism(s) $M^i$ is ``present'' in individual~$\omega$. We may thus associate
potential outcomes with the $A_i$, corresponding to removing (or
inserting) the
corresponding mechanism(s).
This interpretation of the $A_i$'s is consistent with the notion of
``co-cause'' which arises in the literature on minimal sufficient causes.

We note that ``structural'' is often used as a synonym for ``causal.''
However, even under the weak interpretation, a sufficient cause
representation is causal in that it represents a set of potential
outcomes.
The word is used in Definition~\ref{defsupport-counterfact-msc} to
connote that the \textit{structure} of the
representation itself represents (additional) potential outcomes for a set of mechanisms $\mathbf{M}$ that correspond with
the pairs $(A_{i},\mathbf{B}_{i})$, $A_{i}\mathbf{\in A}$, $\mathbf{B}_{i}
\mathbf{\in{\mathfrak{B}}}$. Note that there need not be a unique
structural representation $(\mathbf{A},{\mathfrak{B}})$ for
${\mathcal{D}}(\mathbf{C},\Omega)$.
There might be several functionally equivalent, yet substantively
different, generative mechanisms corresponding to a given pair $(A_i,
\mathbf{B}_i)$; see Example~\ref{exprevent}.

%
%pr7.9 #&#
\begin{proposition} \label{propstructural-implies-exhaustive}
If a representation $(\mathbf{A},{\mathfrak{B}})$ for
${\mathcal{D}}(\mathbf{C},\Omega)$, where $\mathbf{A}= \langle
A_1,\ldots, A_p\rangle$, is structural, then the associated set of
generative mechanisms
$\langle M^1,\ldots, M^p\rangle$ is exhaustive.
\end{proposition}

\begin{pf}This follows from Definitions~\ref{defpopsuffcause}
and~\ref{defexhaust}.
\end{pf}

%
%pr7.10 #&#
\begin{proposition} \label{propexists-structural-rep}
Suppose that $\mathbf{M}$ forms an exhaustive
set of generative mechanisms for ${%
\mathcal{D}}(\mathbf{C},\Omega)$, and $\mathbf{M}$ supports
counterfactuals. If for all $M^i \in\mathbf{M}$ there exists $%
\mathbf{B}_i \in\dot{\mathbb{P}}({\mathbb{L}}(\mathbf{C}))$,
and an $A_i$
such that for all $\mathbf{c}$, and $\omega\in\Omega$, if $A_i(\omega
) =1$, then $(M^i)_\mathbf{c}(\omega) = \bigwedge  (\mathbf
{B}_i)_\mathbf{c}$, then $(\mathbf{A} = \langle A_1,\ldots, A_p\rangle
, \mathfrak{B} = \langle\mathbf{B}_1,\ldots, \mathbf{B}_p\rangle)$
forms a representation for ${%
\mathcal{D}}(\mathbf{C},\Omega)$ that is structural.
\end{proposition}

\begin{pf} This follows from Definitions~\ref{defpopsuffcause}
and~\ref{defsupport-counterfact-msc}.\vadjust{\goodbreak}
\end{pf}

%p

%If $(\mathbf{A},{\mathfrak{B}})$ is structural for ${\mathcal{D}}(
%,\Omega)$,We note that `structural' is often used as a synonym for
%`causal'.
%However, even under the weak interpretation, a sufficient cause
%representation is causal in that it represents a set of potential
%outcomes.
%The word is used here to connote that the \textit{structure} of the
%representation itself represents (additional) potential outcomes
%arising
%from intervention on a set of mechanisms $\mathbf{M}$ that correspond
%with
%the pairs $(A_{i},\mathbf{B}_{i})$, $A_{i}\mathbf{\in A}$, $
%structural representation $(\mathbf{A},{\mathfrak{B}})$ for ${

%
%pr7.11 #&#
\begin{proposition}\label{propmechanism-exists}
If there is some representation $(\mathbf{A},{\mathfrak{B}})$ that is
structural, and $\mathbf{B}\in\dot{\mathbb{P}}(\mathbb{L}(%
\mathbf{C}))$ is irreducible for ${\mathcal{D}}(\mathbf{C};\Omega)$, then
there exists a mechanism $M_{i}$ that is active only if $\mathbf
{B}=\mathbf{1}
$.
\end{proposition}

\begin{pf} If $\mathbf{B}$ is irreducible for ${\mathcal
{D}}(\mathbf{C};\Omega)$,
then there exists $\mathbf{B}_{i}\in{\mathfrak{B}}$ with $\mathbf{%
B\subseteq B}_{i}$; the mechanism $M^{i}=A_{i}\wedge(\bigwedge  (\mathbf{B}_{i}))$ is
such that $M^{i}=1$ only if $\mathbf{B}=\mathbf{1}$.
\end{pf}

Note that the conclusion of Proposition~\ref{propmechanism-exists},
unlike those of Propositions~\ref{propmechanism-irreducible} and \ref
{propmechanism-singular}, %are
%global in that they
does not make reference to an individual $\omega^*$.
This is because Proposition~\ref{propmechanism-exists}
assumes that there is a representation $(\mathbf{A},{\mathfrak{B}})$
that is structural: in
this representation
the $A_i$ variables may be seen as a constituent part of the
corresponding mechanism $M_i$.
%, not merely an indicator of a sub-population of units in which $M_i$
%is present.

Note that there may exist a set of exhaustive
generative mechanisms, but these mechanisms may not themselves be
conjunction of literals so that there
is no sufficient cause representation for ${\mathcal{D}}(\mathbf
{C};\Omega)$ that is structural:

%
%ex7 #&#
\begin{example}Suppose $\mathbf{C}=\{X_{1},X_{2}\}$ where $X_{1}$ and $X_{2}$
again denote the presence of variant alleles, acquired by
maternal and paternal inheritance, respectively, at a particular locus.
Let $M$ denote a protein
that is produced if and only if either $X_{1}=X_{2}=1$ or
$X_{1}=X_{2}=0$ and let $D$
denote some characteristic that occurs if and only if $M=1$.  Suppose
we can intervene
to remove or add the protein. We then have that
\begin{eqnarray*}
M_{x_{1}x_{2}}(\omega) &=&x_{1}x_{2}
\vee(1-x_{1}) (1-x_{2}),
\\
D_{x_{1}x_{2}}(\omega) &=&M_{x_{1}x_{2}}(\omega),
\\
D_{x_{1}x_{2}m}(\omega)=D_{m}(\omega) &=&m.
\end{eqnarray*}
Thus $\{M \}$ constitutes an exhaustive set of generative mechanisms
for ${\mathcal{D}}(\mathbf{C}%
,\Omega)$. We have the following representation for ${\mathcal
{D}}(\mathbf{C},\Omega)$:%
\[
D_{x_{1}x_{2}}(\omega)=\bigl(A_{1}(\omega)X_{1}X_{2}
\vee A_{2}(\omega ) (1-X_{1}) (1-X_{2})
\bigr)_{x_{1}x_{2}}
\]
with $A_{1}(\omega)=A_{2}(\omega)=1$ for all $\omega\in\Omega$.
However, this representation is not structural because $A_{1}(\omega
)X_{1}X_{2}$ and $A_{2}(\omega)(1-X_{1})(1-X_{2})$ do not constitute
separate mechanisms for which interventions are conceivable; there is only
one mechanism $M$, the protein. Since for any $\omega\in\Omega$, $%
D_{11}(\omega)=1$ and $D_{10}(\omega)=D_{01}(\omega)=0$, $\{
X_{1},X_{2}\}$
is irreducible relative to $\mathbf{C}$; however it is not the case that
there is a mechanism $M_{i}$ that is active only if $X_{1}X_{2}=1$ since
for the only mechanism $M$ it is the case that $M=1$ if
$X_{1}=X_{2}=0$.  Note, however, in this example there is still a mechanism, namely $M$, that
will be ``active'' if $X_{1}=X_{2}=1$ but will be ``inactive'' if only
one of $X_{1}$
or $X_{2}$ is~$1$.\looseness=1
\end{example}

\begin{example}
To illustrate the results in the paper we consider again the data presented
in Table~\ref{tabdata}. We let $D$ denote bladder cancer, $X_{1}$
smoking, $X_{2}$ the
S NAT2 genotype and $X_{3}$ the *10 allele on NAT1. As discussed in Example~\ref{exthree-way}, if the effect of ${\mathbf C}=\{X_{1},X_{2},X_{3}\}
$ is unconfounded for $D$ and
we fit the model
%
%e26 #&#
%e7.2 ###
\begin{eqnarray}\label{eqm1}
&& E[D = 1\mid X_{1} = x_{1},X_{2} =
x_{2},X_{3} = x_{3}]
\nonumber
\\
&&\qquad=\beta_{\varnothing}+
\beta_{1}x_{1}+\beta_{2}x_{2}+
\beta_{3}x_{3}+\beta_{12}x_{1}x_{2}\\
&&\qquad\quad{}+\beta_{13}x_{1}x_{3}+
\beta_{23}x_{2}x_{3}+\beta_{123}x_{1}x_{2}x_{2},\nonumber
\end{eqnarray}
then if $X_{1}$, $X_{2}$ and $X_{3}$ have positive monotonic effects on $D$
(relative to $\mathbf{C}$), then $\{X_{1},X_{2},X_{3}\}$ is irreducible
relative to ${\mathbf C}$
if any of the following hold:
\[
\beta_{123} > \beta_{1}, \qquad\beta_{123} >
\beta_{2}, \qquad\beta_{123} > \beta_{3}.
\]
We cannot fit model (\ref{eqm1}) directly with case control data.  However, under
the assumption that the outcome is rare (reasonable with bladder
cancer) so
that odds ratios approximate risk ratios, we can fit the model%
%
%e27 #&#
%e7.3 ###
\begin{eqnarray}\label{eqm2}
\qquad&&E[D = 1\mid X_{1} = x_{1},X_{2} =
x_{2},X_{3} = x_{3}]/E[D = 1\mid
X_{1} = 0,X_{2} = 0,X_{3} = 0]
\nonumber
\\[-8pt]
\\[-8pt]
\nonumber
&&\qquad=\theta_{1}x_{1}+\theta_{2}x_{2}+
\theta_{3}x_{3}+\theta_{12}x_{1}x_{2}+
\theta_{13}x_{1}x_{3}+\theta_{23}x_{2}x_{3}+
\theta_{123}x_{1}x_{2}x_{2},
\end{eqnarray}
and the conditions for the irreducibility of $\{X_{1},X_{2},X_{3}\}$
relative to ${\mathbf C}$ under monotonicity of $\{X_{1},X_{2},X_{3}\}
$ become%
\[
\theta_{123} > \theta_{1}, \qquad\theta_{123} >
\theta_{2},\qquad \theta_{123} > \theta_{3}.
\]
If we fit model (\ref{eqm2}) using maximum likelihood, we find that
\begin{eqnarray*}
\theta_{123}-\theta_{1}&=&1.21\ (95\% \mbox{ CI:}
-3.83,6.26),
\\
\theta_{123}-\theta_{2}&=&2.93\ (95\% \mbox{ CI:}
-2.85,8.72),
\\
\theta_{123}-\theta_{3}&=&2.97\ (95\% \mbox {CI:}
-2.80,8.74).
\end{eqnarray*}
In each case, under our assumption of no confounding, the point
estimate suggests evidence of irreducibility, under monotonicity of $%
\{X_{1},X_{2},X_{3}\}$, but the sample size is not sufficiently large to
draw this conclusion confidently. With monotonicity of $\{
X_{1},X_{2},X_{3}\}
$, irreducibility also implies a singular interaction for $%
\{X_{1},X_{2},X_{3}\}$. If we assume that only $\{X_{1},X_{2}\}$ or $%
\{X_{1},X_{3}\}$ or $\{X_{2},X_{3}\}$ are monotonic relative to
${\mathbf C}$, then the
conditions for irreducibility in Example~\ref{exthree-way} can be
expressed, respectively, as
\[
\theta_{123} > 1+\theta_{1}+\theta_{2},\qquad
\theta_{123} >1+\theta_{1}+\theta_{3},\qquad
\theta_{123} >1+\theta_{2}+\theta_{3}.
\]
From model (\ref{eqm2}) we have that
\begin{eqnarray*}
\theta_{123} -(1+\theta_{1}+\theta_{2})&=&0.09\
(95\% \mbox{ CI:} -4.77,4.96),
\\
\theta_{123} -(1+\theta_{1}+\theta_{3})&=&0.13\
(95\% \mbox{ CI:} -4.69,4.95),
\\
\theta_{123} -(1+\theta_{2}+\theta_{3})&=&1.86\
(95\% \mbox{ CI:} -3.41,7.12).
\end{eqnarray*}
Again, under our assumption of no confounding, in each case the point
estimate suggests
evidence of irreducibility, under monotonicity of just two of the three
exposures, but the sample size is not sufficiently large to draw this
conclusion confidently. With monotonicity of two of the three exposures,
irreducibility also implies a singular interaction for $\{
X_{1},X_{2},X_{3}\}
$. The test for irreducibility in Example~\ref{exthree-way} without
assumptions about
monotonicity can be expressed as
\[
\theta_{123} > 2+\theta_{1}+\theta_{2}+
\theta_{3}.
\]
From model (\ref{eqm2}) we have that
\[
\theta_{123} -(2+\theta_{1}+\theta_{2}+
\theta_{3}) = -0.99\ (95\% \mbox{ CI:} -5.86,3.88).
\]
In this case, not even the
point estimate is positive.

If $\{X_{1},X_{2},X_{3}\}$ is in fact irreducible and if there exists a
representation that is structural, then it follows by Proposition~\ref{propmechanism-exists}
that there exists a mechanism that is active only if $X_{1}=1,X_{2}=1,X_{3}=1
$.
\end{example}

%s8 #&#
%s8 ###
\section{Concluding remarks}
\label{secconclude}

In this paper we have developed general results for notions of interaction
that we referred to as ``irreducibility'' (aka ``a sufficient cause
interaction'') and ``singularity'' (aka ``a singular interaction'') for
sufficient cause models with an
arbitrary number of dichotomous causes. The theory
could be extended by developing notions of sufficient cause,
irreducibility and singularity for causes and outcomes that are
categorical and/or ordinal in nature; see~\cite{vanderweelesufficient2010}.

% for categorical and ordinal causes and/or for categorical or
%ordinal outcomes
% the extent to which the sufficient cause
%framework and notions of irreducibility and singularity can be
%extended so
%as to .

%The present work could be extended in a number of directions.
%Efficiency,
%bias and power properties of different methods of testing for
%sufficient
%cause interactions could be explored. \ Work could be done on
%extending the
%theory developed in this paper to stochastic counterfactuals. \
%Extensions
%concerning time-varying exposures could also be explored. \ The theory
%developed in this paper might also be able to be related to continuous
%variables. \ Continuous variables can be re-coded as a infinite number
%of
%dichotomous variables of the form ${\mathbb I}(X>x)$ for different
%values of $x$. \ It
%may be of interest to explore whether progress with the
%sufficient-component
%cause framework could be made with these re-coded continuous
%variables. \
%Finally, it can be shown that there are cases in which any sufficient
%cause
%representation must have a sufficient cause with two or more causes in
%its
%conjunction but for which no two particular causes ever need to be
%present
%in the same sufficient cause. \ Further work could be done on
%characterizing
%these cases and deriving further counterfactual and empirical
%conditions for
%such cases.

\section*{Acknowledgments}
We thank Stephen Stigler for pointing us to earlier work by Cayley on minimal
sufficient cause models. We thank Mathias Drton, James Robins and Rekha
Thomas for helpful conversations. We are grateful to the Associate
Editor and Referees for helpful
suggestions which improved the manuscript.

% imsref loaded by akundreckaite, 2012-09-11 11:10:09
%
% imsref loaded by akundreckaite, 2012-09-18 11:21:44

%

%suskaldyti doi

\printaddresses


\begin{thebibliography}{45}
% BibTex style file: ims.bst, 2012-08-21
% Default style options (sort=0,type=number).
% Used options (sort=1,type=number).

%b1 ###
\bibitem{Aickin2002}
\begin{bbook}[author]
\bauthor{\bsnm{Aickin},~\bfnm{M.}\binits{M.}}
(\byear{2002}).
\btitle{Causal Analysis in Biomedicine and Epidemiology Based on Minimal
  Sufficient Causation}.
\bpublisher{Dekker}, \blocation{New York}.
\bptok{imsref}%
\end{bbook}
\endbibitem

%b2 ###
\bibitem{Bateson1909}
\begin{bbook}[author]
\bauthor{\bsnm{Bateson},~\bfnm{W.}\binits{W.}}
(\byear{1909}).
\btitle{Mendel's Principles of Heredity}.
\bpublisher{Cambridge Univ. Press}, \blocation{London}.
\bptok{imsref}%
\end{bbook}
\endbibitem

%b3 ###
\bibitem{blisstoxicity1939}
\begin{barticle}[author]
\bauthor{\bsnm{Bliss},~\bfnm{C.~I.}\binits{C.~I.}}
(\byear{1939}).
\btitle{The toxicity of poisons applied jointly}.
\bjournal{Annals of Applied Biology}
\bvolume{26}
\bpages{585--615}.
\bptok{imsref}%
\end{barticle}
\endbibitem

%b4 ###
\bibitem{cayley1853}
\begin{barticle}[author]
\bauthor{\bsnm{Cayley},~\bfnm{A.}\binits{A.}}
(\byear{1853}).
\btitle{Note on a question in the theory of probabilities}.
\bjournal{London, Edinburgh and Dublin Philosophical Magazine}
\bvolume{VI}
\bpages{259}.
\bptok{imsref}%
\end{barticle}
\endbibitem

%b5 ###
\bibitem{cayley1889}
\begin{barticle}[author]
\bauthor{\bsnm{Cayley},~\bfnm{A.}\binits{A.}}
(\byear{1889}).
\btitle{A theorem on trees}.
\bjournal{Quart. J. Math.}
\bvolume{23}
\bpages{376--378}.
\bptok{imsref}%
\end{barticle}
\endbibitem

%b6 ###
\bibitem{Cordell2002}
\begin{barticle}[pbm]
\bauthor{\bsnm{Cordell},~\bfnm{Heather~J.}\binits{H.~J.}}
(\byear{2002}).
\btitle{Epistasis: What it means, what it doesn't mean, and statistical methods
  to detect it in humans}.
\bjournal{Hum. Mol. Genet.}
\bvolume{11}
\bpages{2463--2468}.
\bid{issn={0964-6906}, pmid={12351582}}
\bptok{imsref}%
\end{barticle}
\endbibitem

%b7 ###
\bibitem{Cox1958}
\begin{bbook}[mr]
\bauthor{\bsnm{Cox},~\bfnm{D.~R.}\binits{D.~R.}}
(\byear{1958}).
\btitle{Planning of Experiments}.
\bpublisher{Wiley}, \blocation{New York}.
\bid{mr={0095561}}
\bptok{imsref}%
\end{bbook}
\endbibitem

%b8 ###
\bibitem{dedekind1897}
\begin{barticle}[author]
\bauthor{\bsnm{Dedekind},~\bfnm{Richard}\binits{R.}}
(\byear{1897}).
\btitle{{\"U}ber {Z}erlegungen von {Z}ahlen durch ihre gr{\"{o}}{\ss}ten
  gemeinsamen {T}eiler}.
\bjournal{Gesammelte Werke}
\bvolume{2}
\bpages{103--148}.
\bptok{imsref}%
\end{barticle}\vadjust{\goodbreak}
\endbibitem

%b9 ###
\bibitem{Flanders2006}
\begin{barticle}[author]
\bauthor{\bsnm{Flanders},~\bfnm{D.}\binits{D.}}
(\byear{2006}).
\btitle{Sufficient-component cause and potential outcome models}.
\bjournal{Eur. J. Epidemiol.}
\bvolume{21}
\bpages{847--853}.
\bptok{imsref}%
\end{barticle}
\endbibitem

%b10 ###
\bibitem{fukuda2005}
\begin{bmisc}[author]
\bauthor{\bsnm{Fukuda},~\bfnm{Komei}\binits{K.}}
(\byear{2005}).
\bhowpublished{\texttt{cddlib} reference manual. Technical report, EPFL
  Lausanne and ETH Z{\"{u}}rich. Available at
  \href{http://www.ifor.math.ethz.ch/\textasciitilde fukuda/cdd\_home/cdd.html}{www.ifor.math.ethz.ch/\textasciitilde fukuda/cdd\_home/cdd.html}.}
\bptok{imsref}%
\end{bmisc}
\endbibitem

%b11 ###
\bibitem{Greenland2002}
\begin{barticle}[pbm]
\bauthor{\bsnm{Greenland},~\bfnm{Sander}\binits{S.}} \AND
  \bauthor{\bsnm{Brumback},~\bfnm{Babette}\binits{B.}}
(\byear{2002}).
\btitle{An overview of relations among causal modelling methods}.
\bjournal{Int. J. Epidemiol.}
\bvolume{31}
\bpages{1030--1037}.
\bid{issn={0300-5771}, pmid={12435780}}
\bptok{imsref}%
\end{barticle}
\endbibitem

%b12 ###
\bibitem{Greenland1988}
\begin{barticle}[pbm]
\bauthor{\bsnm{Greenland},~\bfnm{S.}\binits{S.}} \AND
  \bauthor{\bsnm{Poole},~\bfnm{C.}\binits{C.}}
(\byear{1988}).
\btitle{Invariants and noninvariants in the concept of interdependent effects}.
\bjournal{Scand. J. Work Environ. Health}
\bvolume{14}
\bpages{125--129}.
\bid{issn={0355-3140}, pii={1945}, pmid={3387960}}
\bptok{imsref}%
\end{barticle}
\endbibitem

%b13 ###
\bibitem{Koopman1981}
\begin{barticle}[pbm]
\bauthor{\bsnm{Koopman},~\bfnm{J.~S.}\binits{J.~S.}}
(\byear{1981}).
\btitle{Interaction between discrete causes}.
\bjournal{Am. J. Epidemiol.}
\bvolume{113}
\bpages{716--724}.
\bid{issn={0002-9262}, pmid={7234861}}
\bptok{imsref}%
\end{barticle}
\endbibitem

%b14 ###
\bibitem{Mackie1965}
\begin{barticle}[author]
\bauthor{\bsnm{Mackie},~\bfnm{J.~L.}\binits{J.~L.}}
(\byear{1965}).
\btitle{Causes and conditions}.
\bjournal{American Philosophical Quarterly}
\bvolume{2}
\bpages{245--255}.
\bptok{imsref}%
\end{barticle}
\endbibitem

%b15 ###
\bibitem{marcovitz2001}
\begin{bbook}[author]
\bauthor{\bsnm{Marcovitz},~\bfnm{A.~B.}\binits{A.~B.}}
(\byear{2001}).
\btitle{Introduction to Logic Design}.
\bpublisher{McGraw-Hill}, \blocation{New York}.
\bptok{imsref}%
\end{bbook}
\endbibitem

%b16 ###
\bibitem{mccluskey1956}
\begin{barticle}[mr]
\bauthor{\bsnm{McCluskey},~\bfnm{E.~J.}\binits{E.~J.} \bsuffix{Jr.}}
(\byear{1956}).
\btitle{Minimization of {B}oolean functions}.
\bjournal{Bell System Tech. J.}
\bvolume{35}
\bpages{1417--1444}.
\bid{issn={0005-8580}, mr={0082876}}
\bptok{imsref}%
\end{barticle}
\endbibitem

%(\byear{1923}).
%  Agaricales: Essay des principle}.
%In \bbooktitle{Statistical Science 5}
%(\beditor{\bfnm{D.}\binits{D.}~\bsnm{Dabrowska}} \AND
%  \beditor{\bfnm{T.}\binits{T.}~\bsnm{Speed}}, eds.)

%b17 ###
\bibitem{Novick2004}
\begin{barticle}[pbm]
\bauthor{\bsnm{Novick},~\bfnm{Laura~R.}\binits{L.~R.}} \AND
  \bauthor{\bsnm{Cheng},~\bfnm{Patricia~W.}\binits{P.~W.}}
(\byear{2004}).
\btitle{Assessing interactive causal influence}.
\bjournal{Psychol. Rev.}
\bvolume{111}
\bpages{455--485}.
\bid{doi={10.1037/0033-295X.111.2.455}, issn={0033-295X}, pii={2004-12248-008},
  pmid={15065918}}
\bptok{imsref}%
\end{barticle}
\endbibitem

%b18 ###
\bibitem{Pearl2000}
\begin{bbook}[mr]
\bauthor{\bsnm{Pearl},~\bfnm{Judea}\binits{J.}}
(\byear{2000}).
\btitle{Causality: Models, Reasoning, and Inference}.
\bpublisher{Cambridge Univ. Press}, \blocation{Cambridge}.
\bid{mr={1744773}}
\bptok{imsref}%
\end{bbook}
\endbibitem

%b19 ###
\bibitem{Phillips2008}
\begin{barticle}[pbm]
\bauthor{\bsnm{Phillips},~\bfnm{Patrick~C.}\binits{P.~C.}}
(\byear{2008}).
\btitle{Epistasis---the essential role of gene interactions in the structure
  and evolution of genetic systems}.
\bjournal{Nat. Rev. Genet.}
\bvolume{9}
\bpages{855--867}.
\bid{doi={10.1038/nrg2452}, issn={1471-0064}, mid={NIHMS116882}, pii={nrg2452},
  pmcid={2689140}, pmid={18852697}}
\bptok{imsref}%
\end{barticle}
\endbibitem

%b20 ###
\bibitem{quine1952}
\begin{barticle}[mr]
\bauthor{\bsnm{Quine},~\bfnm{W.~V.}\binits{W.~V.}}
(\byear{1952}).
\btitle{The problem of simplifying truth functions}.
\bjournal{Amer. Math. Monthly}
\bvolume{59}
\bpages{521--531}.
\bid{issn={0002-9890}, mr={0051191}}
\bptok{imsref}%
\end{barticle}
\endbibitem

%b21 ###
\bibitem{quine1955}
\begin{barticle}[mr]
\bauthor{\bsnm{Quine},~\bfnm{W.~V.}\binits{W.~V.}}
(\byear{1955}).
\btitle{A way to simplify truth functions}.
\bjournal{Amer. Math. Monthly}
\bvolume{62}
\bpages{627--631}.
\bid{issn={0002-9890}, mr={0075886}}
\bptok{imsref}%
\end{barticle}
\endbibitem

%b22 ###
\bibitem{Robins1986}
\begin{barticle}[mr]
\bauthor{\bsnm{Robins},~\bfnm{James}\binits{J.}}
(\byear{1986}).
\btitle{A new approach to causal inference in mortality studies with a
  sustained exposure period---application to control of the healthy worker
  survivor effect}.
\bjournal{Math. Modelling}
\bvolume{7}
\bpages{1393--1512}.
\bid{doi={10.1016/0270-0255(86)90088-6}, issn={0270-0255}, mr={0877758}}
\bptok{imsref}%
\end{barticle}
\endbibitem

%b23 ###
\bibitem{Robins1999}
\begin{bincollection}[mr]
\bauthor{\bsnm{Robins},~\bfnm{James~M.}\binits{J.~M.}}
(\byear{2000}).
\btitle{Marginal structural models versus structural nested models as tools for
  causal inference}.
In \bbooktitle{Statistical Models in Epidemiology, the Environment, and
  Clinical Trials ({M}inneapolis, {MN}, 1997)}.
\bseries{IMA Vol. Math. Appl.}
\bvolume{116}
\bpages{95--133}.
\bpublisher{Springer}, \blocation{New York}.
\bid{doi={10.1007/978-1-4612-1284-3_2}, mr={1731682}}
\bptnote{check year}%
\bptok{imsref}%
\end{bincollection}
\endbibitem

%b24 ###
\bibitem{Rosenbaum1983}
\begin{barticle}[mr]
\bauthor{\bsnm{Rosenbaum},~\bfnm{Paul~R.}\binits{P.~R.}} \AND
  \bauthor{\bsnm{Rubin},~\bfnm{Donald~B.}\binits{D.~B.}}
(\byear{1983}).
\btitle{The central role of the propensity score in observational studies for
  causal effects}.
\bjournal{Biometrika}
\bvolume{70}
\bpages{41--55}.
\bid{doi={10.1093/biomet/70.1.41}, issn={0006-3444}, mr={0742974}}
\bptok{imsref}%
\end{barticle}
\endbibitem

%b25 ###
\bibitem{Rothman1976}
\begin{barticle}[pbm]
\bauthor{\bsnm{Rothman},~\bfnm{K.~J.}\binits{K.~J.}}
(\byear{1976}).
\btitle{Causes}.
\bjournal{Am. J. Epidemiol.}
\bvolume{104}
\bpages{587--592}.
\bid{issn={0002-9262}, pmid={998606}}
\bptok{imsref}%
\end{barticle}
\endbibitem

%b26 ###
\bibitem{Rothman1998}
\begin{bbook}[author]
\bauthor{\bsnm{Rothman},~\bfnm{K.~J.}\binits{K.~J.}} \AND
  \bauthor{\bsnm{Greenland},~\bfnm{S.}\binits{S.}}
(\byear{1998}).
\btitle{Modern Epidemiology}.
\bpublisher{Lippincott-Raven}, \blocation{Philadelphia}.
\bptok{imsref}%
\end{bbook}
\endbibitem

%b27 ###
\bibitem{Rubin1974}
\begin{barticle}[author]
\bauthor{\bsnm{Rubin},~\bfnm{D.~B.}\binits{D.~B.}}
(\byear{1974}).
\btitle{Estimating causal effects of treatments in randomized and nonrandomized
  studies}.
\bjournal{J. Educ. Psychol.}
\bvolume{66}
\bpages{688--701}.
\bptok{imsref}%
\end{barticle}
\endbibitem

%b28 ###
\bibitem{Rubin1978}
\begin{barticle}[mr]
\bauthor{\bsnm{Rubin},~\bfnm{Donald~B.}\binits{D.~B.}}
(\byear{1978}).
\btitle{Bayesian inference for causal effects: The role of randomization}.
\bjournal{Ann. Statist.}
\bvolume{6}
\bpages{34--58}.
\bid{issn={0090-5364}, mr={0472152}}
\bptok{imsref}%
\end{barticle}
\endbibitem

%b29 ###
\bibitem{Rubin1990}
\begin{barticle}[mr]
\bauthor{\bsnm{Rubin},~\bfnm{Donald~B.}\binits{D.~B.}}
(\byear{1990}).
\btitle{Comment on {J}. {N}eyman and causal inference in experiments and
  observational studies: ``{O}n the application of probability theory to
  agricultural experiments. {E}ssay on principles. {S}ection 9''
  [\textit{{A}nn. {A}gric. {S}ci.} \textbf{10} (1923) 1--51]}.
\bjournal{Statist. Sci.}
\bvolume{5}
\bpages{472--480}.
\bid{issn={0883-4237}, mr={1092987}}
\bptok{imsref}%
\end{barticle}
\endbibitem


%b30 ###
\bibitem{Neyman1923}
\begin{barticle}[mr]
\bauthor{\bsnm{Splawa-Neyman},~\bfnm{Jerzy}\binits{J.}}
(\byear{1990}).
\btitle{On the application of probability theory to agricultural experiments.
  {E}ssay on principles. {S}ection 9 [\textit{{A}nn. {A}gric. {S}ci.} \textbf{10} (1923) 1--51]}.
\bjournal{Statist. Sci.}
\bvolume{5}
\bpages{465--472}.
\bnote{Translated from the Polish and edited by D. M. D\c{a}browska and
  T. P. Speed}.
\bid{issn={0883-4237}, mr={1092986}}
\bptok{imsref}%
\end{barticle}
\endbibitem


%b31 ###
\bibitem{Taylor1998}
\begin{barticle}[author]
\bauthor{\bsnm{Taylor},~\bfnm{Jack~A.}\binits{J.~A.}},
  \bauthor{\bsnm{Umbach},~\bfnm{David~M.}\binits{D.~M.}},
  \bauthor{\bsnm{Stephens},~\bfnm{Elizabeth}\binits{E.}},
  \bauthor{\bsnm{Castranio},~\bfnm{Trisha}\binits{T.}},
  \bauthor{\bsnm{Paulson},~\bfnm{David}\binits{D.}},
  \bauthor{\bsnm{Robertson},~\bfnm{Gary}\binits{G.}},
  \bauthor{\bsnm{Mohler},~\bfnm{James~L.}\binits{J.~L.}} \AND
  \bauthor{\bsnm{Bell},~\bfnm{Douglas~A.}\binits{D.~A.}}
(\byear{1998}).
\btitle{The role of {N}-acetylation polymorphisms in smoking-associated bladder
  cancer: Evidence of a gene-gene-exposure three-way interaction}.
\bjournal{Cancer Research}
\bvolume{58}
\bpages{3603--3610}.
\bptok{imsref}%
\end{barticle}
\endbibitem

%b32 ###
\bibitem{VanderWeele2010a}
\begin{barticle}[pbm]
\bauthor{\bsnm{VanderWeele},~\bfnm{Tyler~J.}\binits{T.~J.}}
(\byear{2010}).
\btitle{Empirical tests for compositional epistasis}.
\bjournal{Nat. Rev. Genet.}
\bvolume{11}
\bpages{166}.
\bid{doi={10.1038/nrg2579-c1}, issn={1471-0064}, pii={nrg2579-c1},
  pmid={20084088}}
\bptok{imsref}%
\end{barticle}
\endbibitem

%b33 ###
\bibitem{VanderWeele2010b}
\begin{barticle}[mr]
\bauthor{\bsnm{VanderWeele},~\bfnm{Tyler~J.}\binits{T.~J.}}
(\byear{2010}).
\btitle{Epistatic interactions}.
\bjournal{Stat. Appl. Genet. Mol. Biol.}
\bvolume{9}
\bpages{24}.
\bid{doi={10.2202/1544-6115.1517}, issn={1544-6115}, mr={2594940}}
\bptok{imsref}%
\end{barticle}
\endbibitem

%b34 ###
\bibitem{vanderweelesufficient2010}
\begin{barticle}[mr]
\bauthor{\bsnm{VanderWeele},~\bfnm{Tyler~J.}\binits{T.~J.}}
(\byear{2010}).
\btitle{Sufficient cause interactions for categorical and ordinal exposures
  with three levels}.
\bjournal{Biometrika}
\bvolume{97}
\bpages{647--659}.
\bid{doi={10.1093/biomet/asq030}, issn={0006-3444}, mr={2672489}}
\bptok{imsref}%
\end{barticle}
\endbibitem

%b35 ###
\bibitem{VanderWeele2006}
\begin{barticle}[pbm]
\bauthor{\bsnm{VanderWeele},~\bfnm{Tyler~J.}\binits{T.~J.}} \AND
  \bauthor{\bsnm{Hern{\'{a}}n},~\bfnm{Miguel~A.}\binits{M.~A.}}
(\byear{2006}).
\btitle{From counterfactuals to sufficient component causes and vice versa}.
\bjournal{Eur. J. Epidemiol.}
\bvolume{21}
\bpages{855--858}.
\bid{doi={10.1007/s10654-006-9075-0}, issn={0393-2990}, pmid={17225959}}
\bptok{imsref}%
\end{barticle}
\endbibitem

%b36 ###
\bibitem{vanderweeleremarks2011}
\begin{barticle}[pbm]
\bauthor{\bsnm{VanderWeele},~\bfnm{Tyler~J.}\binits{T.~J.}} \AND
  \bauthor{\bsnm{Knol},~\bfnm{Mirjam~J.}\binits{M.~J.}}
(\byear{2011}).
\btitle{Remarks on antagonism}.
\bjournal{Am. J. Epidemiol.}
\bvolume{173}
\bpages{1140--1147}.
\bid{doi={10.1093/aje/kwr009}, issn={1476-6256}, pii={kwr009}, pmcid={3121324},
  pmid={21490044}}
\bptok{imsref}%
\end{barticle}
\endbibitem

%b37 ###
\bibitem{VanderWeele2007a}
\begin{barticle}[author]
\bauthor{\bsnm{VanderWeele},~\bfnm{T.~J.}\binits{T.~J.}} \AND
  \bauthor{\bsnm{Robins},~\bfnm{J.~M.}\binits{J.~M.}}
(\byear{2007}).
\btitle{The identification of synergism in the sufficient-component cause
  framework}.
\bjournal{Epidemiol.}
\bvolume{18}
\bpages{329--339}.
\bptok{imsref}%
\end{barticle}
\endbibitem

%b38 ###
\bibitem{VanderWeele2007c}
\begin{barticle}[pbm]
\bauthor{\bsnm{VanderWeele},~\bfnm{Tyler~J.}\binits{T.~J.}} \AND
  \bauthor{\bsnm{Robins},~\bfnm{James~M.}\binits{J.~M.}}
(\byear{2007}).
\btitle{Directed acyclic graphs, sufficient causes, and the properties of
  conditioning on a common effect}.
\bjournal{Am. J. Epidemiol.}
\bvolume{166}
\bpages{1096--1104}.
\bid{doi={10.1093/aje/kwm179}, issn={0002-9262}, pii={kwm179}, pmid={17702973}}
\bptok{imsref}%
\end{barticle}
\endbibitem

%b39 ###
\bibitem{VanderWeele2008a}
\begin{barticle}[mr]
\bauthor{\bsnm{VanderWeele},~\bfnm{Tyler~J.}\binits{T.~J.}} \AND
  \bauthor{\bsnm{Robins},~\bfnm{James~M.}\binits{J.~M.}}
(\byear{2008}).
\btitle{Empirical and counterfactual conditions for sufficient cause
  interactions}.
\bjournal{Biometrika}
\bvolume{95}
\bpages{49--61}.
\bid{doi={10.1093/biomet/asm090}, issn={0006-3444}, mr={2409714}}
\bptok{imsref}%
\end{barticle}
\endbibitem

%b40 ###
\bibitem{VanderWeele2008c}
\begin{barticle}[mr]
\bauthor{\bsnm{VanderWeele},~\bfnm{Tyler~J.}\binits{T.~J.}} \AND
  \bauthor{\bsnm{Robins},~\bfnm{James~M.}\binits{J.~M.}}
(\byear{2010}).
\btitle{Signed directed acyclic graphs for causal inference}.
\bjournal{J. R. Stat. Soc. Ser. B Stat. Methodol.}
\bvolume{72}
\bpages{111--127}.
\bid{doi={10.1111/j.1467-9868.2009.00728.x}, issn={1369-7412}, mr={2751246}}
\bptok{imsref}%
\end{barticle}
\endbibitem

%b41 ###
\bibitem{vanderweelemarginal2010}
\begin{barticle}[pbm]
\bauthor{\bsnm{VanderWeele},~\bfnm{Tyler~J.}\binits{T.~J.}},
  \bauthor{\bsnm{Vansteelandt},~\bfnm{Stijn}\binits{S.}} \AND
  \bauthor{\bsnm{Robins},~\bfnm{James~M.}\binits{J.~M.}}
(\byear{2010}).
\btitle{Marginal structural models for sufficient cause interactions}.
\bjournal{Am. J. Epidemiol.}
\bvolume{171}
\bpages{506--514}.
\bid{doi={10.1093/aje/kwp396}, issn={1476-6256}, pii={kwp396}, pmcid={2877448},
  pmid={20067916}}
\bptok{imsref}%
\end{barticle}
\endbibitem

%b42 ###
\bibitem{Vansteelandt2003}
\begin{barticle}[mr]
\bauthor{\bsnm{Vansteelandt},~\bfnm{S.}\binits{S.}} \AND
  \bauthor{\bsnm{Goetghebeur},~\bfnm{E.}\binits{E.}}
(\byear{2003}).
\btitle{Causal inference with generalized structural mean models}.
\bjournal{J. R. Stat. Soc. Ser. B Stat. Methodol.}
\bvolume{65}
\bpages{817--835}.
\bid{doi={10.1046/j.1369-7412.2003.00417.x}, issn={1369-7412}, mr={2017872}}
\bptok{imsref}%
\end{barticle}
\endbibitem

%b43 ###
\bibitem{vansteelandtmarginal2012}
\begin{barticle}[mr]
\bauthor{\bsnm{Vansteelandt},~\bfnm{Stijn}\binits{S.}},
  \bauthor{\bsnm{VanderWeele},~\bfnm{Tyler~J.}\binits{T.~J.}} \AND
  \bauthor{\bsnm{Robins},~\bfnm{James~M.}\binits{J.~M.}}
(\byear{2012}).
\btitle{Semiparametric tests for sufficient cause interaction}.
\bjournal{J. R. Stat. Soc. Ser. B Stat. Methodol.}
\bvolume{74}
\bpages{223--244}.
\bid{doi={10.1111/j.1467-9868.2011.01011.x}, issn={1369-7412}, mr={2899861}}
\bptok{imsref}%
\end{barticle}
\endbibitem

%b44 ###
\bibitem{wiedemann1991}
\begin{barticle}[mr]
\bauthor{\bsnm{Wiedemann},~\bfnm{Doug}\binits{D.}}
(\byear{1991}).
\btitle{A computation of the eighth {D}edekind number}.
\bjournal{Order}
\bvolume{8}
\bpages{5--6}.
\bid{doi={10.1007/BF00385808}, issn={0167-8094}, mr={1129608}}
\bptok{imsref}%
\end{barticle}
\endbibitem

\end{thebibliography}
\end{document}